\documentclass[11pt]{article}
\usepackage{amsmath,amsfonts,amssymb,theorem}
\usepackage{amsmath}
\usepackage{amsfonts}
\usepackage{amssymb}
\usepackage[usenames,dvipsnames]{pstricks}
\usepackage{epsfig}
\usepackage{pst-grad}
\usepackage{pst-plot}
\usepackage{multicol}
\usepackage{mathdots}

\usepackage{pst-node}
\usepackage{graphicx}
\usepackage{xcolor}
\usepackage{pgf}
\usepackage{soul}

\def\eref#1{(\ref{#1})}

\newtheorem{theorem}{Theorem}[section]
\newtheorem{proposition}[theorem]{Proposition}
\newtheorem{lemma}[theorem]{Lemma}
\newtheorem{corollary}[theorem]{Corollary}
\newtheorem{definition}[theorem]{Definition}
\newtheorem{conjecture}{Conjecture}
\newtheorem{example}[theorem]{Example}
\newtheorem{problem}{Problem}
\newtheorem{remark}[theorem]{Remark}

\begin{document}

\author{Gi-Sang Cheon\thanks{Applied Algebra and Optimization Research Center, Department of Mathematics, Sungkyunkwan University, Suwon 16419, Republic of Korea}, Ji-Hwan Jung\footnotemark[2], Sergey Kitaev
\thanks{Department of Computer and Information Sciences, University of Strathclyde, 26 Richmond Street, Glasgow, G1 1XH, United Kingdom}\  \ and Seyed Ahmad Mojallal\footnotemark[2]\\
{\footnotesize gscheon@skku.edu},\,{\footnotesize
jh56k@skku.edu},\,{\footnotesize sergey.kitaev@cis.strath.ac.uk,}\,{\footnotesize mojallal@skku.edu}}

\title{Riordan graphs I: Structural properties\thanks{This work was supported by the National Research Foundation of Korea (NRF) grant funded by the Korea government (MSIP) (2016R1A5A1008055) and the Ministry of Education of Korea (NRF-2016R1A6A3A11930452).}}
\date{}
\maketitle

\begin{abstract}
In this paper, we use the theory of Riordan matrices to introduce
the notion of a Riordan graph. The Riordan graphs are a far-reaching
generalization of the well known and well studied Pascal graphs and
Toeplitz graphs, and also some other families of graphs. The Riordan
graphs are proved to have a number of interesting (fractal)
properties, which can be useful in creating computer networks with
certain desirable features, or in obtaining useful information when
designing algorithms to compute values of graph invariants. The main
focus in this paper is the study of structural properties of
families of Riordan graphs obtained from infinite Riordan graphs,
which includes a fundamental decomposition theorem and certain
conditions on Riordan graphs to have an Eulerian trail/cycle or a
Hamiltonian cycle. We will study spectral properties of the Riordan
graphs in a follow up paper. \\



\end{abstract}

\section{Introduction}

{\em Pascal's triangle} is a classical combinatorial object, and its
roots can be traced back to the 2nd century BC. In 1991, Shapiro,
Getu, Woan and Woodson \cite{SGWW} have introduced the notion of a
{\em Riordan array}, also known as a {\em Riordan matrix}, in order
to define a class of infinite lower triangular matrices with
properties analogous to those of the Pascal triangle
(matrices)~\cite{ES}:
\begin{eqnarray}\label{e:ppascal}
P=\left[{i\choose j}\right]_{i,j\ge 0}=\left(\begin{array}{cccccc}
       1 & &  &   &  \\
       1 & 1 &  &  &O&  \\
       1 & 2 & 1 &    &&  \\
       1 & 3 &3 & 1 & &  \\
       1 & 4 &6 &4 & 1& \\
       \vdots &\vdots&\vdots & \vdots&\vdots& \ddots \end{array}\right).
\end{eqnarray}
The idea here is to express each column of the Pascal matrix in terms of
generating functions. Since the $j$th column generating function is
 \begin{eqnarray*}
 \sum_{i=0}^\infty{i\choose j}t^i={t^{j}\over (1-t)^{j+1}}={1\over 1-t}\left({t\over
 1-t}\right)^j,\;j\ge0,
\end{eqnarray*}
we move from one column to the next by multiplying by ${t\over
1-t}$. It is convenient to start our numbering at 0 so that the leftmost column is the 0th column. In general, we start with
$g(t)=g_0+g_{1}t+g_{2}t^{2}+\cdots$ as the 0th column generating
function and then multiply by $f(t)=f_{1}t+f_{2}t^{2}+\cdots$ so
that the next column has the generating function
$g(t)f(t)=g_0f_1t+(g_0f_2+g_1f_1)t^2+\cdots$. Continuing, we
multiply by $f(t)$ again to obtain the generating function
$g(t)f(t)^2$ for the 2nd column, and so on. We usually abbreviate $g(t)$ by $g$,
$f(t)$ by $f$, etc, and the resulting matrix is called the Riordan matrix $(g,f)$. Thus, the generating function for the $j$th
column is $gf^j$ for $j=0,1,2,\ldots$. In the Pascal matrix we have
$g={1\over 1-t}$ and $f={t\over 1-t}$.

Since their introduction,  Riordan matrices became an active subject of
research; see \cite{CH,CJ,He, He1, Alg, MerSpr, S2, Sprugnoli} and
references there in, for examples of results in this direction. In
particular, Riordan matrices found applications in the context of
the computation of combinatorial sums \cite{Sprugnoli}. Also, see~\cite{CLMPS} for a recent paper about Lie theory on the
Riordan group, which is the set of invertible Riordan matrices.

The notion of Pascal's triangle was also influential in graph theory
and its applications to computer networks. Indeed, in 1983, Deo and
Quinn \cite{DQ} have introduced the {\it Pascal graphs} in their
searching for a class of optimal graphs for computer networks with
certain desirable properties, such as
\begin{itemize}
\item the design is to be simple and recursive;
\item there must be a universal vertex, i.e.\ a vertex adjacent to all others;
\item there must exist several paths between each pair of vertices.
\end{itemize}
The adjacency matrix of the Pascal graph with $n$
vertices is denoted by $PG_n$  and is defined to be an $n\times n$ symmetric $(0,1)$-matrix
whose main diagonal entries are all 0s and the lower triangular
part of the matrix consists of the first $n-1$ rows of Pascal's
triangle modulo 2. These graphs have attracted much attention in the
literature (see \cite{BD,CKM,DBD} and references there in). For
example, the Pascal graph with 6 vertices is given by
\begin{center}
 \begin{tabular}{cc}\label{e:Pascalgraph}
${\cal A}(PG_6)=\left(\begin{array}{cccccc}
       0 & 1 & 1 & 1 & 1 & 1 \\
       1 & 0 & 1 & 0 & 1 & 0 \\
       1 & 1 & 0 & 1 & 1 & 0 \\
       1 &0 &1 & 0 & 1 & 0 \\
       1 &1 &1 &1 & 0 & 1 \\
       1 &0 &0 & 0 &1 & 0
     \end{array}\right),$& \scalebox{0.5} 
{
\begin{pspicture}(0,0)(7.02625,0.0939063)
\psdots[dotsize=0.2](1.0253808,1.6071479)
\usefont{T1}{ptm}{m}{n}
\rput(0.7522559,1.7371479){\Large$1$}
\psdots[dotsize=0.2](3.8253808,1.6071479)
\usefont{T1}{ptm}{m}{n}
\rput(4.1,1.7571479){\Large$2$}
\psline[linewidth=0.04cm](1.0453808,1.6271479)(3.8253808,1.6271479)
\usefont{T1}{ptm}{m}{n}
\rput(4.1,-1.0628521){\Large$3$}
\psline[linewidth=0.04cm](3.8253808,1.6071479)(3.8253808,-0.9928521)
\psline[linewidth=0.04cm](1.0253808,1.6071479)(3.8253808,-0.9928521)
\psdots[dotsize=0.2](3.8253808,-0.9928521)
\usefont{T1}{ptm}{m}{n}
\psdots[dotsize=0.2](1.0253808,-0.9928521)
\usefont{T1}{ptm}{m}{n}
\rput(0.72631836,-1.0628521){\Large$4$}
\psline[linewidth=0.04cm](1.0053809,1.6471479)(1.0253808,-0.9328521)
\psline[linewidth=0.04cm](1.0253808,-0.9928521)(3.8253808,-0.9928521)
\psdots[dotsize=0.2](2.425381,-2.7928522)
\usefont{T1}{ptm}{m}{n}
\rput(2.4149122,-3.1){\Large$5$}
\psarc[linewidth=0.04](2.4226904,-0.36202875){2.4226904}{125.74507}{269.2245}
\rput{-180.22153}(4.9239626,-0.73357636){\psarc[linewidth=0.04](2.4626904,-0.36202875){2.4226904}{90.34283}{235.15048}}
\psline[linewidth=0.04cm](1.0453808,-1.0128521)(2.405381,-2.7728522)
\psline[linewidth=0.04cm](3.8053808,-0.9928521)(2.445381,-2.7728522)
\psdots[dotsize=0.2](6.165381,2.0071478)
\usefont{T1}{ptm}{m}{n}
\rput(6.4606935,2.057148){\Large$6$}
\rput{-0.18480174}(0.0,0.012048266){\psarc[linewidth=0.04](3.735381,0.0371479){3.13}{245.90797}{150.00491}}
\end{pspicture}
}\\
 &\\
 \end{tabular}
  {\centerline{ Fig.\ 1: The Pascal graph $PG_6$ with 6 vertices}}
 \end{center}

Another important object of interest to us is the well studied
notion of a {\em Toeplitz graph} with $n$ vertices, which is denoted by $TG_n$ (see \cite{DTTVZZ,Eul,Eul3,
EZ,Ghorban,IB,NP} and references there in for results on Toeplitz graphs).
The adjacency matrix ${\cal A}(TG_n)=[t_{ij}]$ is an $n\times n$
symmetric $(0,1)$-Toeplitz matrix defined by
\begin{align*}
t_{ij}=\left\{
\begin{array}{ll}
0 & \text{if $i=j$} \\
a_{|i-j|}\in\{0,1\} & \text{if $i\ne j$.}%
\end{array}
\right.
\end{align*}
Using Riordan language, ${\cal A}(TG_n)$ can be expressed as a
Riordan matrix of the Appell type, $(g,t)$ modulo 2 where
$g=\sum_{n\ge0}a_nt^n\in{\mathbb Z}_2[[t]]$. See
Section~\ref{families-sec} for more information about Toeplitz
graphs. In this context, if $g=t^{t_1-1}+t^{t_2-1}+\cdots+t^{t_k-1}$, where $1\le t_1<t_2<\cdots<t_k< n$,
then $TG_n$ is denoted by $T_n\langle t_1, t_2,\ldots,
t_k\rangle$. For example, the
Toeplitz graph $T_6\langle 1,3,5\rangle$
is the bipartite graph given by
 \begin{center}
 \begin{tabular}{cc}
${\cal A}(TG_6)=\left(\begin{array}{cccccc}
       0 & 1 & 0 & 1 & 0 & 1 \\
       1 & 0 & 1 & 0 & 1 & 0 \\
       0 & 1 & 0 & 1 & 0 & 1 \\
       1 & 0 & 1 & 0 & 1 & 0 \\
       0 & 1 & 0 & 1 & 0 & 1 \\
       1 & 0 & 1 & 0 & 1 & 0
     \end{array}\right)$,& \scalebox{0.7} 
{
\begin{pspicture}(0.1,0)(6.62,1.8639063)
\psdots[dotsize=0.2](0.1,1.2804687)
\psdots[dotsize=0.2](3.3,1.2804687)
\psdots[dotsize=0.2](6.5,1.2804687)
\psdots[dotsize=0.2](0.1,-1.3195312)
\psdots[dotsize=0.2](3.3,-1.3195312)
\psdots[dotsize=0.2](6.5,-1.3195312)
\usefont{T1}{ptm}{m}{n}
\rput(0.046875,1.6904688){1}
\usefont{T1}{ptm}{m}{n}
\rput(3.2676563,1.6904688){3}
\usefont{T1}{ptm}{m}{n}
\rput(6.469531,1.6904688){5}
\usefont{T1}{ptm}{m}{n}
\rput(0.07859375,-1.7095313){2}
\usefont{T1}{ptm}{m}{n}
\rput(3.2809374,-1.7095313){4}
\usefont{T1}{ptm}{m}{n}
\rput(6.4753127,-1.7095313){6}
\psline[linewidth=0.04cm](0.1,1.2604687)(0.1,-1.2595313)
\psline[linewidth=0.04cm](0.1,1.2804687)(3.28,-1.2995312)
\psline[linewidth=0.04cm](0.1,1.3204688)(6.48,-1.3195312)
\psline[linewidth=0.04cm](3.26,1.3004688)(0.16,-1.2595313)
\psline[linewidth=0.04cm](3.3,1.2604687)(3.3,-1.2395313)
\psline[linewidth=0.04cm](3.32,1.2804687)(6.46,-1.2595313)
\psline[linewidth=0.04cm](6.5,1.3004688)(0.12,-1.3195312)
\psline[linewidth=0.04cm](6.48,1.2604687)(3.32,-1.2995312)
\psline[linewidth=0.04cm](6.5,1.2804687)(6.5,-1.2795312)
\end{pspicture}
}\\
 \end{tabular}
 \vskip.5pc
 \centerline{Fig.\ 2: The Toeplitz graph $TG_6$ with 6 vertices}
 \end{center}

In this paper we introduce the notion of a {\em Riordan graph}. This notion not only provides a far-reaching generalizaiton of the notions of a Pascal graph and a Toeplitz graph, but also extends the theory of Riordan matrices to the domain of graph theory, similarly to the introduction of the Pascal graphs based on Pascal's triangles.  The Riordan graphs are proved to have a number of interesting (fractal) properties, which can be useful in creating computer networks with
certain desired features, or in obtaining useful information when designing algorithms to compute values of graph invariants. Indeed, the Pascal graphs are just one instance of a family of Riordan graphs having a universal vertex and a simple recursive structure. Thus, other members of that family could be used instead of the Pascal graphs in designing computer networks, and depending on the context, these could be a better choice than the Pascal graphs.

Our basic idea here is in building the infinite adjacency matrix
$\mathcal{A}$ based on an infinite Riordan matrix over ${\mathbb Z}$
modulo 2, and considering the leading principal submatrices giving
(finite) Riordan graphs. See Section~\ref{Riordan-graphs-subsec} for
the precise definitions. We introduce various families of Riordan
graphs based on the choice of the generating functions defining
these graphs via Riordan matrices. For example, Riordan graphs
of the {\em Appell type} are precisely the class of Toeplitz graphs
including the complete graphs $K_n$ and bipartite graphs, while
Riordan graphs of the {\em Bell type} include Pascal graphs,
{\em Catalan graphs} and {\em Motzkin graphs}.

 One of the basic questions in our context is
whether or not a given labeled or unlabeled graph is a Riordan
graph. It turns out that all unlabeled graphs on at most four
vertices are Riordan graphs (see Theorem~\ref{unlabeled}), while
non-Riordan unlabeled graphs always exist for larger graphs on any
number of vertices. However, the main focus in this paper is the
study of {\em labeled} Riordan graphs, and we give structural
properties of certain families of graphs obtained from infinite
Riordan graphs.

Throughout this paper, we normally label graphs on $n$ vertices by
the elements of the $n$-set $[n]:=\{1,2,\ldots,n\}$. However, we
also meet graphs labeled by odd numbers, or even numbers, or
consecutive subintervals in $[n]$. For two isomorphic graphs, $G$
and $H$, we write $G\cong H$. For a graph $G$, $V$ (resp., $E$)
denotes the set of vertices (resp., edges) in $G$. Also, for a
subset $S$ of vertices $V$, we let $\left< S\right>$ denote the
graph induced by the vertices in $S\subset V$. Finally, we let
$\mathbb{N}_0=\{0,1,\ldots\}$.

This paper is organized as follows. In Section~\ref{sec2} we review
the notion of a Riordan matrix and use it to introduce the notion of
an (infinite) Riordan graph in the {\em labeled} and {\em unlabeled}
cases. However, the main focus in this paper is the {\em labeled}
case, so unless we use the word ``unlabeled'' explicitely, our
Riordan graphs are {\em labeled}. A number of basic results on
Riordan graphs are established in
Section~\ref{Riordan-graphs-subsec}, and this includes the number of
Riordan graphs on $n$ vertices (see Proposition~\ref{num-graphs}),
and the necessary conditions on Riordan graphs (see
Theorem~\ref{thm1}). In Section~\ref{product-sec} we define the {\em
R-product} $\otimes_R$ of two Riordan graphs and then give its
combinatorial interpretation in terms of directed walks in certain
graphs (see Theorem~\ref{e:th8}). We also discuss the {\em ring sum}
$\oplus$ of two graphs in Section~\ref{product-sec} that can be used
to define certain classes of Riordan graphs in Section~\ref{families-sec}. These include, but are not
limited to Riordan graphs of the {\em Appell type}, {\em Bell type},
{\em Lagrange type}, {\em checkerboard type}, {\em derivative type},
and {\em hitting time type}.

In Section~\ref{struct-sec} we give structural results applicable to any Riordan graphs. In particular, in Section~\ref{A-seq-sec}
we show that every  Riordan graph is a fractal (see Theorem~\ref{e:coro2}). Also, the {\em reverse relabelling} of proper Riordan graphs
is defined and studied in Section~\ref{revRio}. Further, in Section~\ref{decomp-sec} we prove the {\em Riordan Graph Decompostion} theorem (see Theorem~\ref{e:th}).
In addition, in Section~\ref{eh-sec} we give certain conditions on Riordan graphs to have an Eulerian trail/cycle or a Hamiltonian cycle.

In Section~\ref{ioie-sec}, we consider {\em io-decomposable} and
{\em ie-decomposable} proper Riordan graphs, and provide a
characterization result for these graphs (see Theorem~\ref{e:th12}).
One of the main focuses in this paper is the study of Riordan graphs
of the {\em Bell type} conducted in Section~\ref{Bell-type-sec}. In
particular, we provide two characterization results for
io-decomposable Riordan graphs of the Bell type (see
Lemma~\ref{e:lem3} and Theorem~\ref{e:th11}) and use the results to
show that the {\em Pascal graphs} and {\em Catalan graphs} are
io-decomposable, while the {\em Motzkin graphs} are not
io-decomposable. Also, in Section~\ref{Bell-type-sec} we study basic
graph invariants of io-decomposable Riordan graphs of the Bell type:
number of edges, number of universal vertices, clique number,
chromatic number, diameter, and others. In
Section~\ref{derivative-sec} we provide two characterization results
(Lemma~\ref{e:th2} and Theorem~\ref{e:th3}) for ie-decomposable
Riordan graphs of the {\em derivative type}. Finally, in
Section~\ref{last-sec} we provide concluding remarks and state
directions for further research.

 We study spectral properties of the Riordan graphs in the
follow up paper~\cite{CJKM}.

\section{From Riordan matrices to Riordan graphs}\label{sec2}

After briefly reviewing the notion of a {\em Riordan matrix} in the
Riordan group, we will introduce the notion of a {\em Riordan
graph}. Then, in Section~\ref{families-sec} we introduce various
families of Riordan graphs.

\subsection{A brief introduction to Riordan group}\label{Riordan-matrix-sec}
Consider $\kappa[[t]]$ the ring of formal power series
$g=\sum_{n\geq 0}g_nt^n$ in the variable $t$ over an integral domain
$\kappa$. An important operation on formal power series is
coefficient extraction, $[t^n]g=g_n$.

Let $L=[\ell_{ij}]_{i,j\ge0}$ be an infinite matrix over $\kappa$.
If there exists a pair of generating functions $(g,f)\in
\kappa[[t]]\times \kappa[[t]]$, $f(0)=0$ such
 that
\begin{eqnarray*}
 gf^j=\sum_{i\ge0}\ell_{ij}t^i,\;j\ge0\;\;{\text{or equivalently}}\;\; \ell_{ij}=[t^i]gf^j
\end{eqnarray*}
then the matrix  $L$ is called a {\it Riordan matrix} (or, a {\it
Riordan array})  over $\kappa$ generated by $g$ and $f$. Usually, we
write $L=(g(t),f(t))$ or $(g,f)$. Since $f(0)=0$, every Riordan
matrix $(g,f)$ is an infinite lower triangular matrix. If a Riordan
matrix is invertible, it is called {\it proper}. Note that $(g,f)$
is invertible if and only if $g(0)\ne0$, $f(0)=0$ and
$f^\prime(0)\ne0$, where $f'$ denotes the first derivative of $f$.

It is well known \cite{MRSV} that an infinite lower triangular
matrix $L=[\ell_{i,j}]_{i,j\ge0}$ with $\ell_{0,0}\ne0$ is a proper
Riordan matrix given by $L=(g,f)$ if and only if there exists a
unique pair of sequences $(a_n)_{n\ge0}$ with $a_0\ne0$ and
$(z_n)_{n\ge0}$ such that for all $i\ge j\ge0$,
\begin{eqnarray}
\ell_{i+1,j+1}&=&a_0\ell_{i,j}+a_1\ell_{i,j+1}+\cdots+a_{i-j}\ell_{i,i};\label{e:A-sequence-relation}\\
\ell_{i+1,0}&=&z_0\ell_{i,0}+z_1\ell_{i,1}+\cdots+z_{i-j}\ell_{i,i}\nonumber.
\end{eqnarray}
Equivalently, there exists a unique pair of generating functions
$A=\sum_{n\ge0}a_nt^n$ and $Z=\sum_{n\ge0}z_nt^n$ such that
\begin{eqnarray*}\label{e:eq12}
g={1\over1-tZ(f)},\quad\hbox{and}\quad f=tA(f)\;{\rm or}\; A=t/\bar
f.
\end{eqnarray*}
The sequences $(a_n)_{n\ge0}$ and $(z_n)_{n\ge0}$ are called the
$A$-sequence and the $Z$-sequence of the Riordan matrix $L$,
respectively. For example, the Pascal matrix $P=\left({1\over
1-t},{t\over 1-t}\right)$ in \eref{e:ppascal} has the $A$-sequence
$(1,1,0,\ldots)$ and the $Z$-sequence $(1,0,\ldots)$.

\begin{theorem}{\rm \cite{SGWW}} Consider a matrix equation $Bx=c$ where $x=(x_0,x_1,\ldots)^T$. If
$B$ is a Riordan matrix $(g,f)$, then the resulting column vector
$c=(c_0,c_1,\ldots)^T$ has the generating function
$$c(t)=\sum_{n\ge0}c_nt^n=g(t)(x\circ f)(t)=gx(f)$$
where $x(t)=x_0+x_1t+x_2t^2+\cdots$ and $\circ$ is the composition operator.
\end{theorem}

This property is called the {\it fundamental theorem for Riordan
matrices} (FTRM) and we write it simply as $(g,f)x=gx(f)$ in terms
of generating functions. This leads to the multiplication of Riordan
matrices as
\begin{eqnarray*}\label{e:Rm}
(g,f)*(h,\ell)=(gh(f),\ell(f)).
\end{eqnarray*}
The set of all proper Riordan matrices under the above {\em Riordan
multiplication} in terms of generating functions forms a group
called the \textit{Riordan group} and denoted $({\cal R},*)$ where
$*$ is the usual matrix multiplication. The identity of the group is
$(1,t)$, the usual identity matrix and
$(g,f)^{-1}=({1/g(\overline{f})},\overline{f})$ where $\overline{f}$
is the {\em compositional inverse} of $f$, i.e.\
$\overline{f}(f(t))=f(\overline{f}(t))=t$. Let
\begin{eqnarray*}
{\cal F}_0&=&\{g\in\kappa[[t]]|g(0)\ne0\},\\
{\cal F}_1&=&\{f\in\kappa[[t]]|f(0)=0,f'(0)\ne0\}=t{\cal F}_0.
\end{eqnarray*}
Then $\left({\cal F}_0,\cdot\right)$ and $\left({\cal
F}_1,\circ\right)$ form groups under the convolution $\cdot$ and the
composition $\circ$ of two formal power series, respectively.
Moreover, it can be easily shown that the Riordan group ${\cal R}$
is isomorphic to the semidirect product $\rtimes$ of ${\cal F}_0$ and ${\cal
F}_1$:
$${\cal R}\cong
{\cal F}_0\rtimes{\cal F}_1.$$

It is also well known \cite{He,Alg,S2} that there are several
important subgroups of the Riordan group. These subgroups will be
used for classifying the types of Riordan graphs in
Section~\ref{families-sec}. \vskip.5pc

\noindent
{\bf Appell subgroup.} A Riordan matrix of the form
\begin{eqnarray*}\label{e:appell}
(g,t)=\left(\begin{array}{ccccc} g_0&&&&\\
g_1&g_0&&O&\\
g_2&g_1&g_0&&\\
g_3&g_2&g_1&g_0& \\
\vdots&\vdots&\vdots&\ddots&\ddots\end{array} \right)
\end{eqnarray*}
is called an {\it Appell matrix} which is a Toeplitz matrix. The set
of all invertible Appell matrices forms a normal subgroup of the
Riordan group which is called the Appell subgroup:
$$ {\cal A}:=\{(g,t)\in{\cal R}\mid g\in {\cal
F}_0\}\cong{\cal F}_0.
$$
\vskip.5pc

\noindent
 {\bf Lagrange subgroup.} A Riordan matrix of the form
$(1,f)$ is called a {\it Lagrange matrix}. The set of all invertible
Lagrange matrices forms a subgroup of the Riordan group which is
called the Lagrange subgroup:
$$ {\cal L}:=\{(1,f)\in{\cal R}\mid f\in {\cal
F}_1\}\cong{\cal F}_1.
$$
Lagrange matrices are sometimes called the {\it associated
matrices}.

Using Riordan multiplication we can see that every Riordan
matrix $(g,f)$ can be expressed as the product of the Appell matrix
$(g,t)$ and the Lagrange matrix $(1,f)$. Since ${\cal A}\cap{\cal
L}=\{(1,t)\}$ it proves that the Riordan group is the semidirect
product of ${\cal A}$ and ${\cal L}$:
 $$ {\cal R}=\{(g,f)=(g,t)(1,f)\in{\cal A}{\cal L}\}={\cal A}\rtimes{\cal L}.
 $$
\vskip.5pc

\noindent
{\bf Bell subgroup.} A Riordan matrix of the form $(g,
tg)$ or $(f/t, f)$ is called a {\it Bell matrix}. The set of all
invertible Bell matrices forms a subgroup of the Riordan group which
is called the Bell subgroup:
$$ {\cal B}:=\{(g,tg)\in{\cal R}\mid g\in {\cal
F}_0\}.
$$
Bell matrices are frequently arising in combinatorics. For example,
\begin{itemize}
\item Pascal matrix $\left({1\over 1-t},{t\over 1-t}\right)$ (see
\eref{e:ppascal} for the matrix expression);
\item Catalan matrix $\left(C,tC\right)=\left(\begin{array}{cccccc}
       1 & &  &   &  \\
       1 & 1 &  &  &O&  \\
       2 & 2 & 1 &    &&  \\
       5 & 5 &3 & 1 & &  \\
       14 & 14 &9 &4 & 1& \\
       \vdots &\vdots&\vdots & \vdots&\vdots& \ddots
       \end{array}\right),\;C={1-\sqrt{1-4t}\over
       2t};$
\item Motzkin matrix $\left(M,tM\right)=\left(\begin{array}{cccccc}
       1 & &  &   &  \\
       1 & 1 &  &  &O&  \\
       2 & 2 & 1 &    &&  \\
       4 & 5 &3 & 1 & &  \\
       9 & 12 &9 &4 & 1& \\
       \vdots &\vdots&\vdots & \vdots&\vdots& \ddots
       \end{array}\right),\;M={1-t-\sqrt{1-2t-3t^2}\over 2t^2}$.
       \end{itemize}
\vskip.5pc

\noindent
 {\bf Checkerboard subgroup.} A Riordan matrix of the form
$(g,f)$ for an even function $g$ and an odd function $f$ is called a
{\it checkerboard matrix}. The set of all invertible checkerboard
matrices forms a subgroup of the Riordan group which is called the
checkerboard subgroup. The centralizer of the element $(1,-t)$ of
the Riordan group is the checkerboard subgroup. \vskip.5pc

\noindent
 {\bf Derivative subgroup.} A Riordan matrix of the form
$(f',f)$ is called a {\it derivative matrix}. The set of all invertible derivative
matrices forms a subgroup of the Riordan group which is called the
derivative subgroup. \vskip.5pc

\noindent
{\bf Hitting time subgroup.} A
Riordan matrix of the form $\left(tf'/f,f\right)$ is called a {\it
hitting time matrix}. The set of all invertible hitting time
matrices forms a subgroup of the Riordan group which is called the
hitting time subgroup.

\subsection{(0,1)-Riordan matrices}\label{Binary-matrix-sec}

In this section, we introduce the notion of a $(0,1)$-Riordan matrix
of order $n$, which will play an important role in defining Riordan
graphs with $n$ vertices in the next section.

First note \cite{CLMPS} that the Riordan group $\cal R$ over
$\kappa={\mathbb R}$ or ${\mathbb C}$ can be described as an inverse
limit of an inverse sequence of groups $\left({\cal
R}_n\right)_{n\in{\mathbb N}}$ of {\it finite} matrices. So ${\cal
R}_n$ is a subgroup of the classical Lie group $GL(n,\kappa)$, and
elements in ${\cal R}_n$ are $n\times n$ Riordan matrices denoted by
$(g,f)_n$ which are obtained from (infinite) Riordan matrices
$(g,f)$ by taking their leading principal submatrix of order $n$.

 Now consider Riordan matrices over the
finite field $\kappa={\mathbb Z}_2$. We call a Riordan matrix over
${\mathbb Z}_2$ a {\it binary Riordan matrix} and denote it as ${\cal
B}(g,f)=[b_{ij}]$. So, ${\cal B}(g,f)$ is the $(0,1)$-matrix obtained from a Riordan
matrix $(g,f)$ over ${\mathbb Z}$, where $g,f\in{\mathbb Z}_2[[t]]$
or ${\mathbb Z}[[t]]$ with $f(0)=1$, by taking the coefficients modulo 2, i.e.
\begin{eqnarray*}
 b_{ij}\equiv[t^i]gf^j\;({\rm mod}\;2)\;\;{\text{or equivalently}}\;\; {\cal
B}(g,f)\equiv (g,f)\;({\rm mod}\;2).
\end{eqnarray*}

From now on, we simply write  $a\equiv b$ for $a\equiv b\;({\rm
mod\;2})$. Also, for $f,g\in{\mathbb Z}[[t]]$, $f\equiv g$ means
$[t^n]f\equiv [t^n]g\;({\rm mod\;2})$ for all $n\ge0$. Since many properties for Riordan
matrices are valid for any integral domain $\kappa$, we state, without proof, some
properties for Riordan matrices over ${\mathbb Z}_2$,
which will be useful in this paper.

It is obvious that $\mathcal{B}(g,f)$ is invertible if and only if
$g(0)=1$, $f(0)=0$ and $f'(0)=1$. By the properties of the modulo
operation, it follows from \eref{e:A-sequence-relation} that
 an infinite $(0,1)$-lower triangular
matrix $B=[b_{i,j}]_{i,j\ge0}$ with $b_{0,0}=1$ is an invertible
Riordan matrix given by $B={\cal B}(g,f)$ if and only if there
exists a unique pair of $(0,1)$-sequences called the {\it binary
$A$-sequence} $(1,\tilde{a}_1,\tilde{a}_2,\ldots)$ and the {\it
binary $Z$-sequence} $(\tilde{z}_0,\tilde{z}_1,\ldots)$ such that
for all $i\ge j\ge0$,
\begin{eqnarray*}\label{e:b-asequence}
b_{i+1,j+1}&\equiv&b_{i,j}+\tilde{a}_1b_{i,j+1}+\cdots+\tilde{a}_{i-j}b_{i,i};\\
b_{i+1,0}&\equiv&\tilde{z}_0b_{i,0}+\tilde{z}_1b_{i,1}+\cdots+\tilde{z}_{i-j}b_{i,i}.\nonumber
\end{eqnarray*}
Thus we have the following theorem.

\begin{theorem} Let $(g,f)$ be a proper Riordan matrix over ${\mathbb
Z}$ with the $A$-sequence $(a_n)_{n\ge0}$, $a_0=1$ and the
$Z$-sequence $(z_n)_{n\ge0}$ where $g(0)=1$. Then, $(g,f)\equiv {\cal
B}(g,f)$ if and only if ${a}_n\equiv \tilde{a}_n$ and
$z_n\equiv\tilde{z}_n$ for all $n\ge0$.
\end{theorem}

Applying the fundamental theorem for Riordan matrices over ${\mathbb
Z}_2$ yields
\begin{align}\label{e:ftbrm} \mathcal{B}(g,f)h\equiv
gh(f)
\end{align}
where $h\in{\mathbb Z}_2[[t]]$. By this property it can be
shown easily that the set of all invertible binary Riordan matrices forms a
group under the binary operation
$$
{\cal B}(g,f){\cal B}(h,\ell)\equiv {\cal B}(gh(f),\ell(f)).
$$

Now, let ${\cal B}(g, f)_n$ be the $n\times n$ $(0,1)$-Riordan matrix
obtained from the leading principal submatrix of order $n$ in ${\cal
B}(g,f)$. The group of $n\times n$ invertible $(0,1)$-Riordan
matrices will be denoted by ${\cal R}_n({\mathbb Z}_2)$. If
$\kappa={\mathbb R}$ or ${\mathbb C}$ then the Riordan group ${\cal
R}_n$ over $\kappa$ has infinite order. However, the Riordan group
${\cal R}_n({\mathbb Z}_2)$ is a finite group.

\begin{example}\label{3by3example}{\rm Consider all $3\times3$ Riordan matrices over
${\mathbb Z}_2$: if $g=g_0+g_1t+g_2t^2$ and
$f=f_1t+f_2t^2\in{\mathbb Z}_2[[t]]$ then
\begin{eqnarray*}
{\cal B}(g,f)_3=\left(\begin{array}{ccc}
       g_0 &0 &0 \\
       g_1 & g_0f_1&0\\
       g_2 & g_0f_2+g_1f_1& g_0f_1^2
       \end{array}\right).
\end{eqnarray*}
It can be easily shown that there are 22 binary Riordan matrices of
order 3 and exactly eight of them are invertible. Note that
$g_0=f_1=1$ in this case:
 {\small\begin{eqnarray*}
 &&{\cal
B}(1,t)_3=\left(\begin{array}{ccc}
       1 & 0 & 0\\
      0 & 1&0 \\
      0 & 0& 1
       \end{array}\right),\  {\cal B}(1,t+t^2)_3=\left(\begin{array}{ccc}
       1 &0 & 0\\
      0 & 1&0\\
      0 & 1& 1
       \end{array}\right),\  {\cal B}(1+t^2,t)_3=\left(\begin{array}{ccc}
       1 & 0& 0\\
      0 & 1&0\\
      1 & 0& 1
       \end{array}\right),\\
&& {\cal B}(1+t^2,t+t^2)_3=\left(\begin{array}{ccc}
       1 & 0& 0\\
      0 & 1&0\\
      1 & 1& 1
       \end{array}\right),\  {\cal B}(1+t,t)_3=\left(\begin{array}{ccc}
       1 & 0& 0\\
      1 & 1&0\\
      0 & 1& 1
       \end{array}\right),\  {\cal B}(1+t,t+t^2)_3=\left(\begin{array}{ccc}
       1 & 0&0 \\
      1 & 1&0\\
      0 & 0& 1
       \end{array}\right),\\
&&{\cal B}(1+t+t^2,t)_3=\left(\begin{array}{ccc}
       1 & 0& 0\\
      1 & 1&0\\
      1 & 1& 1
       \end{array}\right),\  {\cal B}(1+t+t^2,t+t^2)_3=\left(\begin{array}{ccc}
       1 & 0&0 \\
      1 & 1&0\\
     1 & 0& 1
       \end{array}\right).
\end{eqnarray*}}Thus, ${\cal R}_3({\mathbb
Z}_2)$ is an eight-element group, and it is isomorphic to the
Dihedral group $D_4$ (also see \cite{CH})}.
\end{example}

\subsection{Riordan graphs}\label{Riordan-graphs-subsec}

The following definition gives the notion of a Riordan graph in both
{\em labeled} and {\em unlabeled} cases. For an $n\times n$ matrix
$A$, let $A(i|j)$ denote the $(n-1)\times (n-1)$ matrix obtained
from $A$ by deleting $i$th row and $j$th column.

\begin{definition}\label{main-labelled} {\rm A simple {\em labeled} graph $G$ with $n$ vertices is a {\em Riordan
graph} of order $n$ if its adjacency matrix ${\cal A}(G)$ is an
$n\times n$ symmetric $(0,1)$-matrix expressed by ${\cal
A}(G)=B+B^T$ for some binary Riordan matrix $B={\cal B}(tg,f)_n$,
i.e. $B$ is an $n\times n$ $(0,1)$-lower triangular matrix whose
main diagonal entries are all zeros such that
\begin{eqnarray*}\label{e:adj0}
B(1|n)={\cal B}(g,f)_{n-1}.
\end{eqnarray*}
We call the graph $G$ a Riordan graph corresponding to the Riordan matrix
$(g,f)$ and we denote such a Riordan graph $G$
  by $G_n(g,f)$, or simply by $G_n$ when the matrix $(g,f)$ is understood from the context, or it is not important.

A simple {\em unlabeled} graph is a {\em Riordan graph} if at least
one of its labeled copies is a Riordan graph.}
\end{definition}

Since the Riordan matrix $(g,f)=[r_{ij}]_{i,j\ge0}$ is indexed
from $i=0$ and $j=0$, if ${\cal A}(G_n)=[a_{ij}]_{1\le i,j\le n}$
where $a_{ii}=0$, then for $i>j\geq 1$,
\begin{eqnarray}\label{e:def2}
a_{ij}\equiv[t^{i-1}]tgf^{j-1}\equiv[t^{i-2}]gf^{j-1}=r_{i-2,j-1}.
 \end{eqnarray}


 Throughout this paper, we assume that every graph is a simple graph with vertex set $V=[n]$ and edge set $E$.

\begin{example} {\rm(Labeled Riordan graphs)}\end{example}
\begin{itemize}
\item [(i)] The Pascal graph $PG_6$ in Fig.\ 1 is the Riordan graph with 6 vertices given by the matrix $\left({1\over 1-t},{t\over 1-t}\right)$.
\item [(ii)] The Toeplitz graph $TG_6$ in Fig.\ 2 is the Riordan graph with 6 vertices given by the matrix $(1+t^2+t^4,t)$.
\item [(iii)] The Riordan graph with $n$ vertices given by the Catalan matrix $(C,tC)$ in Section~\ref{Riordan-matrix-sec} is called the
{\it Catalan graph} $CG_n$. See Fig.\ 3 for $CG_6$.
\end{itemize}
  \begin{center}
  \begin{tabular}{cc}
 $\mathcal{A}(CG_6)=\left(\begin{array}{cccccc}
        0 & 1 & 1 & 0 & 1 & 0 \\
        1 & 0 & 1 & 0 & 1 & 0 \\
        1 & 1 & 0 & 1 & 1 & 1 \\
        0 & 0 & 1 & 0 & 1 & 0 \\
        1 & 1 & 1 & 1 & 0 & 1 \\
        0 & 0 & 1 & 0 & 1 & 0
      \end{array}\right),$& \scalebox{1} 
 {
 \begin{pspicture}(0,0)(3.02625,0.0939063)
 \psdots[dotsize=0.12](0.08,1.2104688)
 \psdots[dotsize=0.12](1.48,1.2104688)
 \psdots[dotsize=0.12](2.88,1.2104688)
 \psdots[dotsize=0.12](0.38,-0.78953123)
 \psdots[dotsize=0.12](1.48,-1.3895313)
 \psdots[dotsize=0.12](2.58,-0.78953123)
 \usefont{T1}{ptm}{m}{n}
 \rput(0.07859375,1.5204687){2}
 \usefont{T1}{ptm}{m}{n}
 \rput(1.5009375,1.5204687){4}
 \usefont{T1}{ptm}{m}{n}
 \rput(2.8953125,1.5204687){6}
 \usefont{T1}{ptm}{m}{n}
 \rput(0.166875,-0.93953127){1}
 \usefont{T1}{ptm}{m}{n}
 \rput(1.3076563,-1.5395312){3}
 \usefont{T1}{ptm}{m}{n}
 \rput(2.6695313,-0.95953125){5}
 \psline[linewidth=0.04cm](0.08,1.2104688)(0.38,-0.78953123)
 \psline[linewidth=0.04cm](0.08,1.2304688)(1.48,-1.3695313)
 \psline[linewidth=0.04cm](0.08,1.2104688)(2.58,-0.78953123)
 \psline[linewidth=0.04cm](1.48,1.2104688)(1.48,-1.3495313)
 \psline[linewidth=0.04cm](1.48,1.2304688)(2.58,-0.74953127)
 \psline[linewidth=0.04cm](2.86,1.2104688)(1.48,-1.4095312)
 \psline[linewidth=0.04cm](2.88,1.2104688)(2.6,-0.76953125)
 \psline[linewidth=0.04cm](0.38,-0.78953123)(2.58,-0.78953123)
 \psline[linewidth=0.04cm](0.38,-0.78953123)(1.48,-1.3895313)
 \psline[linewidth=0.04cm](1.48,-1.3895313)(2.58,-0.78953123)
 \end{pspicture}
 }
  \end{tabular}
   \vskip.6pc
  {\centerline{Fig.\ 3: The Catalan graph $CG_6$ with 6 vertices}}
  \end{center}

\begin{example} {\rm(Unlabeled graphs)}\end{example}
\begin{itemize}
\item [(i)] Let us determine whether the prism graph $Pr_6$ in Fig.\ 4 is a Riordan graph.
\begin{center}
\psscalebox{1.0 1.0} 
{
\begin{pspicture}(0,-1.0985577)(3.3971155,1.0985577)
\psdots[linecolor=black, dotsize=0.2](1.6985576,1.0)
\psdots[linecolor=black, dotsize=0.2](0.098557666,-0.99999994)
\psdots[linecolor=black, dotsize=0.2](3.2985578,-0.99999994)
\psdots[linecolor=black, dotsize=0.2](1.6985576,0.20000003)
\psdots[linecolor=black, dotsize=0.2](1.0803758,-0.59999996)
\psdots[linecolor=black, dotsize=0.2](2.2985578,-0.59999996)
\psline[linecolor=black,
linewidth=0.04](1.6803758,1.0181818)(0.1349213,-0.9272727)
\psline[linecolor=black,
linewidth=0.04](0.11673948,-0.96363634)(3.280376,-0.98181814)
\psline[linecolor=black,
linewidth=0.04](1.7167395,1.0)(3.280376,-0.9454545)
\psline[linecolor=black,
linewidth=0.04](1.6803758,0.21818186)(1.062194,-0.61818177)
\psline[linecolor=black,
linewidth=0.04](1.6985576,0.20000003)(2.262194,-0.5636363)
\psline[linecolor=black,
linewidth=0.04](1.0803758,-0.58181816)(2.244012,-0.58181816)
\psline[linecolor=black,
linewidth=0.04](0.11673948,-0.96363634)(1.0803758,-0.58181816)
\psline[linecolor=black,
linewidth=0.04](1.6985576,1.0)(1.6985576,0.21818186)
\psline[linecolor=black,
linewidth=0.04](3.2985578,-0.96363634)(2.280376,-0.59999996)
\end{pspicture}
}
 \vskip.3pc {\centerline{ Fig.\ 4: The prism graph $Pr_6$ with 6 vertices}}
\end{center}

There is a vertex labeling whose adjacency matrix is ${\cal
A}(Pr_6)=B+B^T$ where $B={\mathcal B}(t^2+t^3+t^4,t)_6$ is a Toeplitz
matrix of order 6 as shown in Fig.\ 5. Thus, $Pr_6$ is
a 3-regular Riordan graph given by the Riordan matrix $(t+t^2+t^3,t)$; see Fig.\ 5.

\begin{center}
\begin{tabular}{cc}
$\mathcal{A}(Pr_6)=\left(\begin{array}{cccccc}
        0 & 0 & 1 & 1 & 1 & 0 \\
        0 & 0 & 0 & 1 & 1 & 1 \\
        1 & 0 & 0 & 0 & 1 & 1 \\
        1 & 1 & 0 & 0 & 0 & 1 \\
        1 & 1 & 1 & 0 & 0 & 0 \\
        0 & 1 & 1 & 1 & 0 & 0
      \end{array}\right),$& \psscalebox{1.0 1.0} 
{
\begin{pspicture}(0,-0.1)(3.4044867,1.402424)
\psdots[linecolor=black, dotsize=0.2](1.6985577,1.1075076)
\psdots[linecolor=black, dotsize=0.2](0.09855771,-0.8924924)
\psdots[linecolor=black, dotsize=0.2](3.2985578,-0.8924924)
\psdots[linecolor=black, dotsize=0.2](1.6985577,0.30750757)
\psdots[linecolor=black, dotsize=0.2](1.0803759,-0.49249244)
\psdots[linecolor=black, dotsize=0.2](2.2985578,-0.49249244)
\psline[linecolor=black,
linewidth=0.04](1.6803759,1.1256894)(0.13492134,-0.81976515)
\psline[linecolor=black,
linewidth=0.04](0.116739534,-0.8561288)(3.280376,-0.8743106)
\psline[linecolor=black,
linewidth=0.04](1.7167395,1.1075076)(3.280376,-0.83794695)
\psline[linecolor=black,
linewidth=0.04](1.6803759,0.32568938)(1.0621941,-0.51067424)
\psline[linecolor=black,
linewidth=0.04](1.6985577,0.30750757)(2.2621942,-0.4561288)
\psline[linecolor=black,
linewidth=0.04](1.0803759,-0.4743106)(2.2440124,-0.4743106)
\psline[linecolor=black,
linewidth=0.04](0.116739534,-0.8561288)(1.0803759,-0.4743106)
\psline[linecolor=black,
linewidth=0.04](1.6985577,1.1075076)(1.6985577,0.32568938)
\psline[linecolor=black,
linewidth=0.04](3.2985578,-0.8561288)(2.280376,-0.49249244)
\rput[bl](1.429783,0.2861637){\scriptsize1}
\rput[bl](0.8368976,-0.44348058){\scriptsize3}
\rput[bl](2.376028,-0.44031852){\scriptsize5}
\rput[bl](3.2044866,-1.25){\scriptsize2}
\rput[bl](1.6013246,1.3){\scriptsize4}
\rput[bl](0.0029055388,-1.25){\scriptsize6}
\end{pspicture}
}
\end{tabular}
 \vskip.5pc
{\centerline{ Fig.\ 5: A labeled prism graph $Pr_6$ with 6 vertices}}
\end{center}

\item [(ii)] The graph in Fig.\ 6 with 5 vertices is not a Riordan graph because there
is no possible vertex labeling of the graph whose adjacency matrix
can be expressed as a binary Riordan matrix (see Proposition
\ref{non-Riordan-graph}).
\end{itemize}

\begin{center}
\psscalebox{1. 1.} 
{
\begin{pspicture}(0,-1.4985577)(2.1971154,1.4985577)
\psdots[linecolor=black, dotsize=0.2](0.098557666,0.6)
\psdots[linecolor=black, dotsize=0.2](2.0985577,0.6)
\psdots[linecolor=black, dotsize=0.2](0.098557666,-1.4)
\psdots[linecolor=black, dotsize=0.2](2.0985577,-1.4)
\psline[linecolor=black,
linewidth=0.04](0.098557666,0.6)(2.0985577,0.6)
\psline[linecolor=black,
linewidth=0.04](2.0985577,0.6)(2.0985577,-1.4)
\psline[linecolor=black,
linewidth=0.04](0.098557666,0.6)(0.098557666,-1.4)
\psline[linecolor=black,
linewidth=0.04](0.098557666,0.6)(2.0985577,-1.4)
\psline[linecolor=black,
linewidth=0.04](2.0985577,0.6)(0.098557666,-1.4)
\psline[linecolor=black,
linewidth=0.04](0.098557666,-1.4)(2.0985577,-1.4)(2.0985577,-1.4)
\psdots[linecolor=black, dotsize=0.2](1.0785576,1.4000001)
\end{pspicture}
}  \vskip.3pc {\centerline{ Fig.\ 6: A non-Riordan graph with 5
vertices}}
\end{center}

\begin{definition} {\rm A Riordan graph $G_n(g,f)$ is {\em proper}
if the binary Riordan matrix ${\cal B}(g,f)_{n-1}$ is
proper.}\end{definition}

If a Riordan graph $G=G_n(g,f)$ is proper then obviously $G$ is a
connected graph, and the Riordan matrix $(g,f)$ is also proper
because $g(0)\equiv f'(0)\equiv 1$. The converse to this statement
is not true. For instance, $(1,2t+t^2)$ is a proper Riordan matrix
but $G_n(1,2t+t^2)$ is not a proper Riordan graph.

For any Riordan graph $G_n=G_n(g,f)$, we can think of the sequence
of induced subgraphs
$$G_1=\left<\{1\}\right>,\;G_2=\left<\{1,2\}\right>,\ldots,
  G_{n-1}=\left<\{1,2,\ldots,n-1\}\right>,
  $$
each one defined by the same pair of functions, showing the recursive nature of Riordan graphs. From applications point of view, this property implies that when a new node is added to a network, the entire network does not have to be reconfigured.

\begin{proposition}\label{num-graphs}
The number of Riordan graphs of order $n\ge1$ is $(4^{n-1}+2)/3$.
\end{proposition}
\noindent\textbf{Proof.} Let $G=G_n(g,f)$ be a labeled Riordan graph
and $i$ be the smallest index such that $g_i=[t^i]g\equiv 1$.

\begin{itemize}
\item If $i\ge n-1$ then $G$ is the null graph $N_n$.
\item If $0\le i\le n-2$ then we may assume that $g=t^i+g_{i+1}t^{i+1}+\cdots +g_{n-2}t^{n-2}\in{\mathbb Z}_2[t]$ and $f=f_1t+f_2t^2+\cdots + f_{n-i-2}t^{n-i-2}\in{\mathbb Z}_2[t]$.
\end{itemize}
It follows that the number of possibilities of the coefficients in
$\{0,1\}$ to create $G$ is
\begin{align*}
1+\sum_{i=0}^{n-2}2^{2(n-i-2)}={4^{n-1}+2\over 3}
\end{align*}
where the 1 corresponds to the null graph. \hfill{\rule{2mm}{2mm}}

\begin{remark} {\rm The sequence $(a_n)_{n\ge0}$, $a_n={4^n+2\over3}$ for the number of Riordan graphs on $n$
vertices has also various combinatorial interpretations as shown in
A047849 in the On-Line Encyclopedia of Integer Sequences
 (OEIS) \cite{oeis}:
$$1,2,6,22,86,342,1366,\ldots.
$$
This sequence has the generating function $(1-3t)/((1-t)(1-4t))$.}
\end{remark}

\begin{definition} {\rm Any Riordan matrix $(g,f)$ over $\mathbb{Z}$ naturally
defines the {\em infinite} graph
$$G:=G(g,f)=\lim_{n\to\infty}G_n(g,f),$$ which we call the {\em
infinite Riordan graph} corresponding to the Riordan matrix
$(g,f)$.}
\end{definition}

We note that even if an {\em unlabeled} graph is Riordan, its random
labeling is likely to result in a non-Riordan graph. The following
theorem gives necessary conditions for a graph to be Riordan. These
conditions are formulated in terms of the subdiagonal elements in
the adjacency matrix of a Riordan graph.

\begin{theorem}[Necessary conditions for Riordan
graphs]\label{thm1} Let $G=G_n(g,f)$ be a Riordan graph with vertex
set $V=[n]$ and edge set $E$. Then, one of the following
holds:
\begin{itemize}
\item[{\rm(i)}] $i(i+1)\notin E$ for $1\le i \le n-1$.
\item[{\rm(ii)}] $12\in E$ and $i(i+1)\notin E$ for $2\le i \le
n-1$.
\item[{\rm(iii)}] $i(i+1)\in E$ for $1\le i \le n-1$, i.e.\ $G$ has the Hamiltonian path $1\rightarrow 2\rightarrow\cdots \rightarrow~n$.
\end{itemize}
\end{theorem}
\noindent\textbf{Proof.} Let ${\cal A}(G)=[a_{i,j}]_{1\leq i,j\leq
n}$. From the definition of the Riordan matrix $(g,f)$, we have
$a_{i+1,i}\equiv g_0 f_1^{i-1}$ for $1\le i \le n-1$ where
$g_0=[t^0]g$ and $f_1=[t^1]f$. Going through the four possibilities
of choosing $g_0$ and $f_1$ in $\{0,1\}$, we obtain the required
result. \hfill{\rule{2mm}{2mm}}\\

It is known  that the number of {\it unlabeled} graphs on
$n$ vertices is given by the sequence $1,2,4,11,34,156,\ldots$, which
is A000088 in the OEIS \cite{oeis}.
\begin{theorem}\label{unlabeled}
All unlabeled graphs on at most four vertices are Riordan graphs.
\end{theorem}
\noindent\textbf{Proof.} As we observed above, the total number of
unlabeled graphs on at most four vertices is 18. The Fig.\ 7
justifies that such graphs are all Riordan graphs. Using the
matrices in Example~\ref{3by3example}, we provide the proper
labeling and the corresponding Riordan matrices $(g,f)$ in the
figure.\hfill{\rule{2mm}{2mm}}\\
\begin{center}
\scalebox{0.7} 
{
\begin{pspicture}(0,-4.52)(21.02,4.52)
\psdots[dotsize=0.2](1.5,3.08)
\usefont{T1}{ptm}{m}{n}
\rput(1.466875,3.37){1} \psdots[dotsize=0.2](4.0,3.04)
\psdots[dotsize=0.2](5.0,3.04)
\usefont{T1}{ptm}{m}{n}
\rput(3.966875,3.37){1}
\usefont{T1}{ptm}{m}{n}
\rput(4.998594,3.37){2} \psdots[dotsize=0.2](7.02,3.02)
\psdots[dotsize=0.2](8.02,3.02)
\usefont{T1}{ptm}{m}{n}
\rput(6.986875,3.35){1}
\usefont{T1}{ptm}{m}{n}
\rput(8.018594,3.35){2}
\psline[linewidth=0.04cm](7.04,3.04)(8.04,3.04)
\psdots[dotsize=0.2](10.5,3.8) \psdots[dotsize=0.2](9.8,2.9)
\psdots[dotsize=0.2](11.3,2.9)
\usefont{T1}{ptm}{m}{n}
\rput(10.466875,4.09){1}
\usefont{T1}{ptm}{m}{n}
\rput(9.558594,2.95){2}
\usefont{T1}{ptm}{m}{n}
\rput(11.547656,2.95){3} \psdots[dotsize=0.2](13.46,3.8)
\psdots[dotsize=0.2](12.76,2.9) \psdots[dotsize=0.2](14.26,2.9)
\usefont{T1}{ptm}{m}{n}
\rput(13.426875,4.07){1}
\usefont{T1}{ptm}{m}{n}
\rput(12.518594,2.93){2}
\usefont{T1}{ptm}{m}{n}
\rput(14.507656,2.93){3}
\psline[linewidth=0.04cm](13.44,3.8)(12.76,2.9)
\psdots[dotsize=0.2](16.48,3.8) \psdots[dotsize=0.2](15.78,2.9)
\psdots[dotsize=0.2](17.28,2.9)
\usefont{T1}{ptm}{m}{n}
\rput(16.446875,4.07){1}
\usefont{T1}{ptm}{m}{n}
\rput(15.538593,2.91){2}
\usefont{T1}{ptm}{m}{n}
\rput(17.527657,2.91){3}
\psline[linewidth=0.04cm](16.5,3.82)(15.82,2.94)
\psline[linewidth=0.04cm](15.82,2.92)(17.26,2.92)
\psdots[dotsize=0.2](19.5,3.8) \psdots[dotsize=0.2](18.8,2.9)
\psdots[dotsize=0.2](20.3,2.9)
\usefont{T1}{ptm}{m}{n}
\rput(19.446875,4.07){1}
\usefont{T1}{ptm}{m}{n}
\rput(18.558594,2.91){2}
\usefont{T1}{ptm}{m}{n}
\rput(20.547657,2.91){3}
\psline[linewidth=0.04cm](19.46,3.8)(18.78,2.92)
\psline[linewidth=0.04cm](18.84,2.92)(20.28,2.92)
\psline[linewidth=0.04cm](19.52,3.8)(20.26,2.94)
\psline[linewidth=0.04cm](0.0,1.5)(21.0,1.5)
\psline[linewidth=0.04cm](0.0,-1.5)(21.0,-1.5)
\psline[linewidth=0.04cm](3.0,4.5)(3.0,-4.5)
\psline[linewidth=0.04cm](6.0,4.5)(6.0,-4.5)
\psline[linewidth=0.04cm](9.0,4.5)(9.0,-4.5)
\psline[linewidth=0.04cm](12.0,4.5)(12.0,-4.5)
\psline[linewidth=0.04cm](15.0,4.5)(15.0,-4.5)
\psline[linewidth=0.04cm](18.0,4.5)(18.0,-4.5)
\psdots[dotsize=0.2](1.0,1.0) \psdots[dotsize=0.2](1.0,0.0)
\psdots[dotsize=0.2](2.0,0.0) \psdots[dotsize=0.2](2.0,1.0)
\usefont{T1}{ptm}{m}{n}
\rput(0.726875,1.01){1}
\usefont{T1}{ptm}{m}{n}
\rput(2.2585938,0.99){2}
\usefont{T1}{ptm}{m}{n}
\rput(0.7409375,0.01){4}
\usefont{T1}{ptm}{m}{n}
\rput(2.2276564,-0.03){3} \psdots[dotsize=0.2](4.02,0.98)
\psdots[dotsize=0.2](4.02,-0.02) \psdots[dotsize=0.2](5.02,-0.02)
\psdots[dotsize=0.2](5.02,0.98)
\usefont{T1}{ptm}{m}{n}
\rput(3.746875,0.99){1}
\usefont{T1}{ptm}{m}{n}
\rput(5.2785935,0.97){2}
\usefont{T1}{ptm}{m}{n}
\rput(3.7609375,-0.01){4}
\usefont{T1}{ptm}{m}{n}
\rput(5.2476563,-0.05){3}
\psline[linewidth=0.04cm](4.0,0.98)(4.98,0.98)
\psdots[dotsize=0.2](7.02,0.98) \psdots[dotsize=0.2](7.02,-0.02)
\psdots[dotsize=0.2](8.02,-0.02) \psdots[dotsize=0.2](8.02,0.98)
\usefont{T1}{ptm}{m}{n}
\rput(6.746875,0.99){1}
\usefont{T1}{ptm}{m}{n}
\rput(8.278594,0.97){2}
\usefont{T1}{ptm}{m}{n}
\rput(6.7609377,-0.01){4}
\usefont{T1}{ptm}{m}{n}
\rput(8.247656,-0.05){3}
\psline[linewidth=0.04cm](7.0,0.96)(8.0,0.02)
\psline[linewidth=0.04cm](8.0,0.96)(7.0,0.02)
\psdots[dotsize=0.2](10.02,0.98) \psdots[dotsize=0.2](10.02,-0.02)
\psdots[dotsize=0.2](11.02,-0.02) \psdots[dotsize=0.2](11.02,0.98)
\usefont{T1}{ptm}{m}{n}
\rput(9.746875,0.99){1}
\usefont{T1}{ptm}{m}{n}
\rput(11.278594,0.97){2}
\usefont{T1}{ptm}{m}{n}
\rput(9.760938,-0.01){4}
\usefont{T1}{ptm}{m}{n}
\rput(11.247656,-0.05){3}
\psline[linewidth=0.04cm](10.02,0.98)(10.98,0.98)
\psline[linewidth=0.04cm](10.02,0.94)(11.0,0.04)
\psdots[dotsize=0.2](13.04,0.98) \psdots[dotsize=0.2](13.04,-0.02)
\psdots[dotsize=0.2](14.04,-0.02) \psdots[dotsize=0.2](14.04,0.98)
\usefont{T1}{ptm}{m}{n}
\rput(12.766875,0.99){1}
\usefont{T1}{ptm}{m}{n}
\rput(14.2985935,0.97){2}
\usefont{T1}{ptm}{m}{n}
\rput(12.780937,-0.01){4}
\usefont{T1}{ptm}{m}{n}
\rput(14.267656,-0.05){3}
\psline[linewidth=0.04cm](13.04,0.98)(14.0,0.98)
\psline[linewidth=0.04cm](14.02,1.0)(14.02,0.02)
\psline[linewidth=0.04cm](13.04,0.0)(14.06,0.0)
\psdots[dotsize=0.2](16.0,0.98) \psdots[dotsize=0.2](16.0,-0.02)
\psdots[dotsize=0.2](17.0,-0.02) \psdots[dotsize=0.2](17.0,0.98)
\usefont{T1}{ptm}{m}{n}
\rput(15.726875,0.99){1}
\usefont{T1}{ptm}{m}{n}
\rput(17.258595,0.97){2}
\usefont{T1}{ptm}{m}{n}
\rput(15.740937,-0.01){4}
\usefont{T1}{ptm}{m}{n}
\rput(17.227655,-0.05){3}
\psline[linewidth=0.04cm](16.0,0.98)(16.96,0.98)
\psline[linewidth=0.04cm](16.0,0.94)(16.0,0.04)
\psline[linewidth=0.04cm](16.94,0.94)(16.02,0.0)
\psdots[dotsize=0.2](19.02,1.0) \psdots[dotsize=0.2](19.02,0.0)
\psdots[dotsize=0.2](20.02,0.0) \psdots[dotsize=0.2](20.02,1.0)
\usefont{T1}{ptm}{m}{n}
\rput(18.746876,1.01){1}
\usefont{T1}{ptm}{m}{n}
\rput(20.278593,0.99){2}
\usefont{T1}{ptm}{m}{n}
\rput(18.760937,0.01){4}
\usefont{T1}{ptm}{m}{n}
\rput(20.247656,-0.03){3}
\psline[linewidth=0.04cm](19.02,1.0)(19.98,1.0)
\psline[linewidth=0.04cm](19.0,0.94)(19.0,0.04)
\psline[linewidth=0.04cm](19.02,1.0)(20.0,0.04)
\psdots[dotsize=0.2](1.02,-2.0) \psdots[dotsize=0.2](1.02,-3.0)
\psdots[dotsize=0.2](2.02,-3.0) \psdots[dotsize=0.2](2.02,-2.0)
\usefont{T1}{ptm}{m}{n}
\rput(0.746875,-1.99){1}
\usefont{T1}{ptm}{m}{n}
\rput(2.2785938,-2.01){2}
\usefont{T1}{ptm}{m}{n}
\rput(0.7609375,-2.99){4}
\usefont{T1}{ptm}{m}{n}
\rput(2.2476563,-3.03){3}
\psline[linewidth=0.04cm](1.02,-1.98)(1.98,-1.98)
\psline[linewidth=0.04cm](2.0,-1.98)(2.0,-2.96)
\psline[linewidth=0.04cm](1.02,-2.98)(2.04,-2.98)
\psline[linewidth=0.04cm](1.0,-2.02)(1.0,-2.96)
\psdots[dotsize=0.2](4.0,-2.0) \psdots[dotsize=0.2](4.0,-3.0)
\psdots[dotsize=0.2](5.0,-3.0) \psdots[dotsize=0.2](5.0,-2.0)
\usefont{T1}{ptm}{m}{n}
\rput(3.726875,-1.99){1}
\usefont{T1}{ptm}{m}{n}
\rput(5.2585936,-2.01){2}
\usefont{T1}{ptm}{m}{n}
\rput(3.7409375,-2.99){4}
\usefont{T1}{ptm}{m}{n}
\rput(5.2276564,-3.03){3}
\psline[linewidth=0.04cm](4.0,-1.98)(4.96,-1.98)
\psline[linewidth=0.04cm](4.98,-1.98)(4.98,-2.96)
\psline[linewidth=0.04cm](4.0,-2.98)(5.02,-2.98)
\psline[linewidth=0.04cm](4.98,-2.02)(3.98,-2.96)
\psdots[dotsize=0.2](7.02,-2.0) \psdots[dotsize=0.2](7.02,-3.0)
\psdots[dotsize=0.2](8.02,-3.0) \psdots[dotsize=0.2](8.02,-2.0)
\usefont{T1}{ptm}{m}{n}
\rput(6.746875,-1.99){1}
\usefont{T1}{ptm}{m}{n}
\rput(8.278594,-2.01){2}
\usefont{T1}{ptm}{m}{n}
\rput(6.7609377,-2.99){4}
\usefont{T1}{ptm}{m}{n}
\rput(8.247656,-3.03){3}
\psline[linewidth=0.04cm](7.02,-1.98)(7.98,-1.98)
\psline[linewidth=0.04cm](8.0,-1.98)(8.0,-2.96)
\psline[linewidth=0.04cm](7.02,-2.98)(8.04,-2.98)
\psline[linewidth=0.04cm](7.0,-2.02)(7.0,-2.96)
\psline[linewidth=0.04cm](7.98,-2.0)(7.06,-2.96)
\psdots[dotsize=0.2](10.02,-2.0) \psdots[dotsize=0.2](10.02,-3.0)
\psdots[dotsize=0.2](11.02,-3.0) \psdots[dotsize=0.2](11.02,-2.0)
\usefont{T1}{ptm}{m}{n}
\rput(9.746875,-1.99){1}
\usefont{T1}{ptm}{m}{n}
\rput(11.278594,-2.01){2}
\usefont{T1}{ptm}{m}{n}
\rput(9.760938,-2.99){4}
\usefont{T1}{ptm}{m}{n}
\rput(11.247656,-3.03){3}
\psline[linewidth=0.04cm](10.02,-1.98)(10.98,-1.98)
\psline[linewidth=0.04cm](11.0,-1.98)(11.0,-2.96)
\psline[linewidth=0.04cm](10.02,-2.98)(11.04,-2.98)
\psline[linewidth=0.04cm](10.0,-2.02)(10.0,-2.96)
\psline[linewidth=0.04cm](10.98,-2.0)(10.06,-2.96)
\psline[linewidth=0.04cm](10.02,-2.0)(10.98,-2.98)
\usefont{T1}{ptm}{m}{n}
\rput(1.4915625,2.15){$G_1(0,t)$}
\usefont{T1}{ptm}{m}{n}
\rput(4.506094,2.15){$G_2(0,t)$}
\usefont{T1}{ptm}{m}{n}
\rput(7.458594,2.15){$G_2(1,t)$}
\usefont{T1}{ptm}{m}{n}
\rput(10.527344,2.15){$G_3(0,t)$}
\usefont{T1}{ptm}{m}{n}
\rput(13.581718,2.15){$G_3(1,0)$}
\usefont{T1}{ptm}{m}{n}
\rput(16.50672,2.15){$G_3(1,t)$}
\usefont{T1}{ptm}{m}{n}
\rput(19.565624,2.15){$G_3(1+t,t)$}
\usefont{T1}{ptm}{m}{n}
\rput(1.5232812,-0.9){$G_4(0,t)$}
\usefont{T1}{ptm}{m}{n}
\rput(4.4659376,-0.9){$G_4(1,0)$}
\usefont{T1}{ptm}{m}{n}
\rput(7.5171876,-0.9){$G_4(t,t)$}
\usefont{T1}{ptm}{m}{n}
\rput(10.501562,-0.9){$G_4(1+t,0)$}
\usefont{T1}{ptm}{m}{n}
\rput(13.516094,-0.9){$G_4(1,t)$}
\usefont{T1}{ptm}{m}{n}
\rput(16.5,-0.9){$G_4(1+t^2,t^2)$}
\usefont{T1}{ptm}{m}{n}
\rput(19.6,-0.9){$G_4(1+t+t^2,0)$}
\usefont{T1}{ptm}{m}{n}
\rput(1.4917188,-3.8){$G_4(1+t^2,t)$}
\usefont{T1}{ptm}{m}{n}
\rput(4.516719,-3.8){$G_4(1,t+t^2)$}
\usefont{T1}{ptm}{m}{n}
\rput(7.52,-3.8){$G_4(1+t^2,t+t^2)$}
\usefont{T1}{ptm}{m}{n}
\rput(10.52,-3.8){$G_4(1+t+t^2,t)$}
\end{pspicture}
}  \vskip.3pc {\centerline{ Fig.\ 7: All unlabeled graphs on at most
4 vertices which are Riordan graphs}}
\end{center}

We note that the Riordan matrices ${\cal B}(1,t+t^2)_3$ and ${\cal B}(1+t,t+t^2)_3$ correspond to the same unlabeled graph, and  ${\cal B}(1+t^2,t+t^2)_3$, ${\cal B}(1+t,t)_3$ and ${\cal B}(1+t+t^2,t+t^2)_3$ correspond to the same unlabeled graph. On the other hand, the following proposition shows that not all
unlabeled graphs on $n$ vertices are Riordan for $n\geq 5$.

\begin{proposition}\label{non-Riordan-graph} The unlabelled graph $H_{n+1}\cong K_n\cup K_1$ obtained from a complete graph $K_n$ by adding an isolated vertex is not Riordan for $n\geq 4$.\end{proposition}

\noindent{\bf Proof.} Let $H_{n+1}\cong K_n\cup K_1$. Suppose that
there exist $g$ and $f$ such that a labeled copy of $H_{n+1}$ is the
Riordan graph $G_n(g,f)$. We consider two cases depending on whether
the isolated vertex in $H_{n+1}$ is labeled by $1$ or not. Assume
that the isolated vertex is labeled by 1. Since there are no edges
$1i\in E$ for $i=2,\ldots,n+1$, we have $g=0$ so that $G_n(g,f)$ is
the null graph $N_n$. This leads to a contradiction. Now, let $i\neq1$
be the label of the isolated vertex and
$\mathcal{A}(H_{n+1})=[a_{i,j}]_{1\leq i,j\leq n+1}$. Since $n$
sub-diagonal entries of $\mathcal{A}(H_{n+1})$ by Theorem~\ref{thm1}
are of the form
\begin{align*}
(a_{2,1},\ldots,a_{n+1,n})=\left\{
\begin{array}{ll}
(\underbrace{1,\ldots,1}_{(i-2){\rm
times}},0,0,\underbrace{1,\ldots,1}_{(n-i){\rm\ times}})
 & \text{if $i\neq n$} \\
(\underbrace{1,\ldots,1}_{(n-1){\rm\ times}},0) & \text{if $i=n,$}%
\end{array}
\right.
\end{align*}
 this also leads to a contradiction. Hence the
proof follows. \hfill{\rule{2mm}{2mm}}

\subsection{Operations on Riordan graphs}\label{product-sec}

 There are many graph operations studied in the literature (see \cite{BM,IK} and
references there in). Basically, graph operations produce new graphs
from initial ones. However, in general we cannot guarantee that a
particular operation applied to Riordan graphs results in a Riordan
graph. In this section, we consider graph operations that preserve
Riordan graphs. An example of such an operation is the {\em ring sum} of
two graphs defined as follows.

\begin{definition} {\rm Given two graphs $G_1=(V_1,E_1)$ and $G_2=(V_2,E_2)$ we define the {\em ring sum} $G_1\oplus G_2=(V_1\cup V_2,(E_1\cup E_2) - (E_1\cap E_2))$.
Thus, an edge is in $G_1\oplus G_2$ if and only if it is an edge in
$G_1$, or and edge in $G_2$, but not both.} \end{definition}

The ring sum is well defined on Riordan graphs $G_n(g,f)$ and
$G_n(h,f)$ with a fixed $f$ and a fixed vertex set (e.g. $[n]$) due
to the fact that $(g,f)+(h,f)=(g+h,f)$,
 so
\begin{eqnarray*} G_n(g,f)\oplus G_n(h,f)=G_n(g+h,f).
\end{eqnarray*}

 Every simple graph $G$ with $n$ vertices has an $n\times n$ adjacency matrix  ${\cal A}(G) =B+B^T$, where $B$ is an $n\times n$ $(0,1)$-lower triangular matrix with
 zero diagonal. In what follows we give a new graph operation defined for any simple graphs, not necessary Riordan graphs.

 \begin{definition} \label{def1} {\rm Let $G$ and $G'$ be simple labeled graphs with the same vertex set
 $V=[n]$ whose adjacency matrices of order $n$ are of the form ${\cal A}(G)=B+B^T$ and ${\cal
 A}(G')=B'+B'^T$, respectivelty. We define the {\it R-product} $\otimes_R$ of $G$ and $G'$ to be the graph $G\otimes_R G'$ with the vertex set $V$
 whose adjacency matrix of order $n$ is of the form:
 $$
 {\cal A}\left(G\otimes_R G'\right)=B_R+B_R^T
 $$
 where $B_R(1|n)=B(1|n)B'(1|n)$.}
 \end{definition}

 In particular, if $G=G_{n}(g,f)$ and
 $G'=G_{n}(h,\ell)$ then the {\it R-product} $\otimes_R$ of $G$ and $G'$ is
the Riordan graph (with the vertex set $V$) given by
\begin{eqnarray}\label{e:de}
G_{n}(g,f) \otimes_R G_{n}(h,\ell)=G_{n}(g h(f),\ell(f)).
\end{eqnarray}

\begin{example}
 {\rm Let $G$ be the Fibonacci graph $G_n(1,t+t^2)$ and $G'=G_n(\frac{1}{1-t},\frac{t}{1-t})$ be the Pascal graph.
 Then, the R-product of these  graphs is
 $$G\otimes_R G'=G_n\left(\frac{1}{1-t-t^2},\frac{t+t^2}{1-t-t^2}\right).
 $$
 }
\end{example}

Using \eref{e:de} we immediately obtain the following theorem.

\begin{theorem} Every Riordan graph $G_n(g,f)$ can be expressed as
the R-product of $G_n(g,t)$ and $G_n(1,f)$. Furthermore, the set of
all proper Riordan graphs forms a group under the R-product
$\otimes_R$ given by \eref{e:de}. The identity of the group is the
path graph $G_n(1,t)$.
\end{theorem}

To give a combinatorial interpretation for the R-product on Riordan
graphs, let $B=[b_{ij}]$ be an $n\times n$ (0,1)-lower triangular
matrix with zero diagonal. We associate with $B$ a {\it digraph}
$D_c(B)$ such that there is a colored arc from $i$ to $j$ using
color $c$ if and only if $b_{ij}=1$.

Consider Riordan graphs $G_{n}(g,f)$ and $G_{n}(h,\ell)$. Let
$B_1={\cal B}(tg,f)_n$, $B_2={\cal B}(t,t)_n^T$ and $B_3={\cal
B}(th,\ell)_n$ be binary Riordan matrices. Assume that
$D_r(B_1)=(V,E_r)$, $D_g(B_2)=(V,E_g)$, and $D_b(B_3)=(V,E_b)$ are
the digraphs with the same vertex set $V=[n]$ and the arc sets
respectively $E_r$, $E_g$ and $E_b$, where $r,g$ and $b$ denote
colors red, green, and blue.

 \begin{definition}\label{rgb} {\rm Given Riordan graphs $G=G_{n}(g,f)$ and $H=G_{n}(h,\ell)$ with the vertex set $V=[n]$, the {\it RGB-graph}
$\mathcal{D}_n=\mathcal{D}(G,H)$ is defined to be the digraph $\left(V,E_r\cup E_g\cup E_b\right)$ where $V=[n]$. In particular, a directed walk
of length $3$ in an RGB-graph is called an {\em RGB-walk} if its
first arc is red, the second arc is green, and the third arc
is blue.}
\end{definition}

\begin{example}\label{exam}
{\rm Let $G=G_{6}\left({1\over 1-t},{t\over 1-t}\right)$ and
$H=G_{6}\left({1\over 1-t^2},t^2\right)$. Then, the RGB-graph
$\mathcal{D}_6=\mathcal{D}(G,H)$ is shown in Fig.\ 8, where solid (resp., dotted, dashed) arrows represent red (resp., green,  blue) arrows.}
\end{example}

\begin{center}
\scalebox{1.3} 
{
\begin{pspicture}(0,-1.74)(10.22,1.74)
\psline[linecolor=black, linewidth=0.03, linestyle=dotted,
dotsep=0.10583334cm](0.1,-0.11)(10.1,-0.11)
\psbezier[linewidth=0.01,linecolor=black](2.1,-0.09)(2.1,0.22236164)(0.14,0.21)(0.14,-0.095532425)
\psbezier[linewidth=0.01,linecolor=black](4.08,-0.07)(4.08,0.50142854)(0.12,0.53)(0.12,-0.04142857)
\psbezier[linewidth=0.01,linecolor=black](6.06,-0.07)(6.06,0.9056097)(0.1,0.93)(0.1,-0.045609757)
\psbezier[linewidth=0.01,linecolor=black](8.06,-0.05)(8.06,1.35)(0.15975063,1.1850997)(0.06,-0.05)
\psbezier[linewidth=0.01,linecolor=black](10.1,-0.07)(10.08,1.6942942)(0.02,1.73)(0.02,-0.07941049)
\psbezier[linewidth=0.01,linecolor=black](8.06,-0.11)(8.06,0.83657944)(2.1,0.91)(2.12,-0.010487805)
\psbezier[linewidth=0.01,linecolor=black](6.06,-0.11)(6.1,0.29)(4.12,0.29)(4.1,-0.09)
\psbezier[linewidth=0.01,linecolor=black](8.04,-0.11)(8.04,0.46142858)(4.08,0.49)(4.08,-0.08142857)
\psbezier[linewidth=0.01,linecolor=black](4.02,-0.07)(4.02,0.24236165)(2.08,0.23)(2.08,-0.07553242)
\psbezier[linewidth=0.01,linecolor=black](10.04,-0.07)(10.06,0.20236164)(8.1,0.23)(8.1,-0.07553242)
\psbezier[linecolor=black, linewidth=0.015, linestyle=dashed,
dash=0.17638889cm
0.10583334cm](0.08,-0.15)(0.08,-0.51)(2.1,-0.51)(2.1,-0.15)
\psbezier[linecolor=black, linewidth=0.015, linestyle=dashed,
dash=0.17638889cm
0.10583334cm](0.06,-0.15363637)(0.06,-1.07)(6.12,-1.07)(6.12,-0.11)
\psbezier[linecolor=black, linewidth=0.015, linestyle=dashed,
dash=0.17638889cm
0.10583334cm](0.04,-0.19)(0.04,-1.6924391)(10.08,-1.73)(10.1,-0.19)
\psbezier[linecolor=black, linewidth=0.015, linestyle=dashed,
dash=0.17638889cm
0.10583334cm](2.12,-0.15)(2.12,-0.71585363)(6.04,-0.73)(6.04,-0.16414635)
\psbezier[linecolor=black, linewidth=0.015, linestyle=dashed,
dash=0.17638889cm
0.10583334cm](2.08,-0.11)(2.08,-1.31)(10.08,-1.41)(10.08,-0.17)
\psbezier[linecolor=black, linewidth=0.015, linestyle=dashed,
dash=0.17638889cm
0.10583334cm](4.08,-0.15)(4.08,-1.0475609)(10.06,-1.07)(10.06,-0.17243902)
\psbezier[linewidth=0.01,linecolor=black](8.06,-0.09)(8.06,0.22236164)(6.1,0.21)(6.1,-0.095532425)
\rput(-0.016279142,
-0.0059292787){\psdots[dotsize=0.1632457,linecolor=black,dotangle=-26.0,dotstyle=triangle*](5.14,1.29)}
\rput(-0.01627916,
-0.0059293453){\psdots[dotsize=0.1632457,linecolor=black,dotangle=-26.0,dotstyle=triangle*](4.12,0.96)}
\rput(-0.01627916,
-0.0059293453){\psdots[dotsize=0.1632457,linecolor=black,dotangle=-26.0,dotstyle=triangle*](3.14,0.699)}
\rput(-0.01627916,
-0.0059293453){\psdots[dotsize=0.1632457,linecolor=black,dotangle=-26.0,dotstyle=triangle*](2.14,0.4)}
\rput(-0.01627916,
-0.0059293453){\psdots[dotsize=0.1632457,linecolor=black,dotangle=-26.0,dotstyle=triangle*](1.12,0.15)}
\rput(-0.01627916,
-0.0059293453){\psdots[dotsize=0.1632457,linecolor=black,dotangle=-26.0,dotstyle=triangle*](3.14,0.18)}
\rput(-0.016279142,
-0.0059292787){\psdots[dotsize=0.1632457,linecolor=black,dotangle=-26.0,dotstyle=triangle*](5.16,0.199)}
\rput(-0.016279142,
-0.0059292787){\psdots[dotsize=0.1632457,linecolor=black,dotangle=-26.0,dotstyle=triangle*](6.16,0.345)}
\rput(-0.016279142,
-0.0059292787){\psdots[dotsize=0.1632457,linecolor=black,dotangle=-26.0,dotstyle=triangle*](5.12,0.66)}
\rput(-0.016279142,
-0.0059292787){\psdots[dotsize=0.1632457,linecolor=black,dotangle=-26.0,dotstyle=triangle*](7.16,0.149)}
\rput(-0.016279142,
-0.0059292787){\psdots[dotsize=0.1632457,linecolor=black,dotangle=-26.0,dotstyle=triangle*](9.14,0.155)}
\rput(-0.016279226,
-0.00592933){\psdots[dotsize=0.1632457,linecolor=black,dotangle=-26.0,dotstyle=triangle*](1.14,-0.4)}
\rput(-0.016279226,
-0.00592933){\psdots[dotsize=0.1632457,linecolor=black,dotangle=-26.0,dotstyle=triangle*](2.92,-0.817)}
\rput(-0.016279226,
-0.00592933){\psdots[dotsize=0.1632457,linecolor=black,dotangle=-26.0,dotstyle=triangle*](4.12,-0.565)}
\rput(-0.01627921,
-0.0059292633){\psdots[dotsize=0.1632457,linecolor=black,dotangle=-26.0,dotstyle=triangle*](5.14,-1.32)}
\rput(-0.01627921,
-0.0059292633){\psdots[dotsize=0.1632457,linecolor=black,dotangle=-26.0,dotstyle=triangle*](7.12,-0.82)}
\rput(-0.01627921,
-0.0059292633){\psdots[dotsize=0.1632457,linecolor=black,dotangle=-26.0,dotstyle=triangle*](6.14,-1.035)}
\rput(-0.02501505,
0.07058902){\psdots[dotsize=0.1632457,linecolor=black,dotangle=-207.0,dotstyle=triangle*](1.1,-0.18)}
\rput(-0.02501505,
0.07058902){\psdots[dotsize=0.1632457,linecolor=black,dotangle=-207.0,dotstyle=triangle*](3.08,-0.18)}
\rput(-0.025014764,
0.07058895){\psdots[dotsize=0.1632457,linecolor=black,dotangle=-207.0,dotstyle=triangle*](5.12,-0.18)}
\rput(-0.025014764,
0.07058895){\psdots[dotsize=0.1632457,linecolor=black,dotangle=-207.0,dotstyle=triangle*](7.12,-0.18)}
\rput(-0.020474859,
0.08949902){\psdots[dotsize=0.1632457,linecolor=black,dotangle=-207.0,dotstyle=triangle*](9.1,-0.19)}
\usefont{T1}{ptm}{m}{n}
\rput(-0.1,-0.225){\scriptsize $1$}
\usefont{T1}{ptm}{m}{n}
\rput(1.85,-0.225){\scriptsize $2$}
\usefont{T1}{ptm}{m}{n}
\rput(3.85,-0.225){\scriptsize $3$}
\usefont{T1}{ptm}{m}{n}
\rput(5.85,-0.225){\scriptsize $4$}
\usefont{T1}{ptm}{m}{n}
\rput(7.85,-0.225){\scriptsize $5$}
\usefont{T1}{ptm}{m}{n}
\rput(10.4,-0.225){\scriptsize $6$} \psdots[dotsize=0.2](0.1,-0.11)
\psdots[dotsize=0.2](2.1,-0.11) \psdots[dotsize=0.2](4.1,-0.11)
\psdots[dotsize=0.2](6.1,-0.11) \psdots[dotsize=0.2](8.1,-0.11)
\psdots[dotsize=0.2](10.1,-0.11)
\end{pspicture}
}\vskip.2pc\centerline{ Fig.\ 8: The RGB-graph
$\mathcal{D}_6=\mathcal{D}(G,H)$}
\end{center}

\begin{theorem}\label{e:th8}
For Riordan graphs $G=G_{n}(g,f)$ and $H=G_{n}(h,\ell)$, let
$B_1=[r_{ij}]$, $B_2=[g_{ij}]$ and $B_3=[b_{ij}]$ be the same as in
Definition~\ref{rgb}. Let $\omega_{ij}$, for $i>j$, be the number
of RGB-walks from $i$ to $j$ in the RGB-graph $\mathcal{D}(G,H)$.
Then,
$$
\omega_{ij}=\sum_{k=j}^{i-1}r_{i,k}g_{k,k+1}b_{k+1,j}.
$$
In particular, two vertices $i$ and $j$ in $G \otimes_R H$ are
adjacent if and only if $\omega_{ij}\equiv 1$.
\end{theorem}
\textbf{Proof.} Since $B_1={\cal B}(tg,f)_{n}$, $B_2={\cal
B}(t,t)_{n}^T$ and $B_3={\cal B}(th,\ell)_{n}$, $B_1$ (resp., $B_2$;
$B_3$) is the adjacency matrix of the directed subgraph of
$\mathcal{D}(G,H)$ formed by the red (resp., green; blue) edges. All
entries in $B_2$ are 0 except $g_{k,k+1}=1$ for $1\leq k\leq n-1$.
Thus, if $B_1B_2B_3=[d_{ij}]_{1\leq i,j\leq n}$ then
\begin{align}\label{e:entries1}
d_{ij}=\left\{
\begin{array}{ll}
\sum_{k=j}^{i-1}r_{i,k}g_{k,k+1}b_{k+1,j}&
\text{if $i>j$} \\
0 & \text{otherwise.}
\end{array}
\right.
\end{align}
It implies that if $i>j$ then $d_{i,j}$ counts the number of
RGB-walks from $i$ to~$j$ in the digraph $\mathcal{D}(G,H)$. Now, let
${\cal A}=[a_{ij}]_{1\leq i,j\leq n}$ be the adjacency matrix of
$G\otimes_R H$. Since
\begin{eqnarray*}
{\cal B}(gh(f),\ell(f))_{n-1}\equiv{\cal B}(g,f)_{n-1}{\cal
B}(h,\ell)_{n-1},
\end{eqnarray*}
we have
\begin{align}\label{e:adj}
a_{ij}\equiv\sum_{k=j}^{i-1}r_{i,k}b_{k+1,j}\;\;\text{if $i>j$}.
\end{align}
Since $g_{k,k+1}=1$ for $1\leq k\leq n-1$ it follows from
\eref{e:entries1} and \eref{e:adj} that $d_{i,j}\equiv a_{i,j}$ if
$i>j$. It means that, for $i>j$, $i$ and $j$  are adjacent, i.e.\
$a_{ij}=1$ if and only if $d_{ij}$ is odd. Hence the proof follows.
\hfill{\rule{2mm}{2mm}}

\begin{example}\label{ex5} {\rm Let $G$ and $H$ be the same Riordan graphs as in Example~\ref{exam}.
 Since $B_1={\cal
B}\left({t\over 1-t},{t\over 1-t}\right)_6$, $B_2={\cal
B}(t,t)_{6}^T$ and $B_3={\cal B}\left({t\over
1-t^2},t^2\right)_{6}$, we obtain
\begin{eqnarray*}
D:=B_1B_2B_3=[d_{i,j}]=\left(\begin{array}{cccccc}
0 & 0 & 0 & 0 & 0 & 0 \\
1 & 0 & 0 & 0 & 0 & 0 \\
1 & 0 & 0 & 0 & 0 & 0 \\
2 & 1 & 0 & 0 & 0 & 0 \\
2 & 1 & 0 & 0 & 0 & 0 \\
2 & 1 & 1 & 0 & 0 & 0
\end{array}\right).
\end{eqnarray*}
 By counting RGB-walks in Fig.\ 8,
we see that if $i>j$ then $d_{i,j}$ counts the number of RGB-walks
from a vertex $i$ to a vertex $j$. For instance, $d_{5,1}=2$
because there are two RGB-walks, $5{\rightarrow} 1{ \rightarrow} 2 {
\rightarrow} 1$ and $5{\rightarrow}3{\rightarrow} 4 {\rightarrow}1$
from the vertex 5 to the vertex 1.

In addition, we obtain the adjacency matrix of $G\otimes_R
H=G_n\left(1+t,{t^2\over1-t^2}\right)$ as follows:
\begin{eqnarray*}
D+D^T\equiv\left(\begin{array}{cccccc}
0 & 1 & 1 & 0 & 0 & 0 \\
1 & 0 & 0 & 1 & 1 & 1 \\
1 & 0 & 0 & 0 & 0 & 1 \\
0 & 1 & 0 & 0 & 0 & 0 \\
0 & 1 & 0 & 0 & 0 & 0 \\
0 & 1 & 1 & 0 & 0 & 0
\end{array}\right),
\end{eqnarray*}
which should be the case by Theorem~\ref{e:th8}.}
\end{example}

\subsection{Families of Riordan graphs}\label{families-sec}

Recall that every Riordan graph $G$ with $n$ vertices is
determined by some Riordan matrix $(g,f)$ over ${\mathbb Z}$ such
that $(g,f)_{n-1}\equiv {\cal B}(g,f)_{n-1}$. Such a graph $G$ is a Riordan graph corresponding to the Riordan matrix $(g,f)$
and is denoted by $G_n(g,f)$. Thus the families of Riordan graphs are
naturally defined by the corresponding classes of Riordan matrices.

In this section, we introduce several classes of Riordan graphs
and give examples of graphs in these classes. The names of the
classes come from the widely used names of Riordan matrices
defining the respective Riordan graphs; such matrices are obtained
by imposing various restrictions on the pairs of functions $(g,f)$.
Also, in Section~\ref{Four-families} we introduce {\em
o-decomposable}, {\em e-decomposable}, {\em io-decomposable} and
{\em ie-decomposable} Riordan graphs.

Note that the most general definition of the {\em null graphs} $N_n$
(also known as the {\em empty graphs}) in our terms is $G_n(0,f)$
for any $f$ where $f(0)=0$.
Also note that the empty graphs, the star graphs $G_n(\frac{1}{1-t},0)$, and the {\em complete $k$-ary trees} for $k\ge2$ defined by $G_n(1+t+\cdots+t^{k-1},t^k)$
are examples of non-proper Riordan graphs; other examples of non-proper Riordan graphs can be obtained from (v) in Theorem~\ref{e:th13},
and even more such examples are discussed at the end of this subsection. However, most of Riordan graphs considered in this paper are proper.
Inspired by Riordan subgroups given in Section~\ref{Riordan-matrix-sec}, we have the following classes of Riordan graphs.  \\[-2mm]

\noindent {\bf Riordan graphs of the {\em Appell type}.} This class
of graphs is defined by an {\em Appell matrix} $(g,t)$. Let
$G=G_n(g,t)$ be a Riordan graph of the Appell type and
$g=\sum_{i=0}^{n-2}g_{i}t^i\in{\mathbb Z}_2[t]$. Then the adjacency
matrix $\mathcal{A}(G)$ is given by the symmetric
$(0,1)$-Toeplitz matrix
\begin{align*}
\mathcal{A}(G_n)=\left(\begin{array}{cccccc}
  0      & g_0    & g_1    & \cdots & \cdots & g_{n-2}\\
  g_0    & 0      & g_0    & \ddots &        & \vdots\\
  g_1    & g_0    & \ddots & \ddots & \ddots & \vdots\\
  \vdots &\ddots  & \ddots & \ddots &  g_0   & g_1\\
  \vdots &        & \ddots & g_0    & 0      & g_0\\
  g_{n-2}& \cdots & \cdots & g_1    & g_0    &0\\
\end{array}\right).
\end{align*}
Thus, the class of Riordan graphs of the Appell type is precisely the
class of Toeplitz graphs. Examples of graphs in this class are
\begin{itemize}
\item the null graphs $N_n$ defined by $G_n(0,t)$;
\item the path graphs $P_n$ defined by $G_n(1,t)$;
\item the complete graphs $K_n$ defined by
$G_n\left(\frac{1}{1-t},t\right)$; and
\item the complete bipartite graphs $K_{\lfloor{n\over2}\rfloor, \lceil{n\over2}\rceil}$ defined by $G_n({1\over 1-t^2},t).$
\end{itemize}

\noindent {\bf Riordan graphs of the {\em Bell type}.} This class of
graphs is defined by a {\em Bell matrix} $(g,tg)$. Examples of
graphs in this class are
\begin{itemize}
\item the null graphs $N_n$ defined by $G_n(0,0)$;
\item the path graphs $P_n$ defined by $G_n(1,t)$;
\item the Pascal graphs $PG_n$ defined by $G_n({1\over 1-t},{t\over 1-t})$;
\item the {\em Catalan graphs} $CG_n$ defined by $G_n\left({1-\sqrt{1-4t}\over 2t},{1-\sqrt{1-4t}\over 2}\right)$; and
\item the {\em Motzkin graphs} $MG_n$ defined by $G_{n}\left( \frac{1-t-\sqrt{1-2t-3t^{2}}}{2t^{2}},\frac{1-t-
\sqrt{1-2t-3t^{2}}}{2t}\right)$.
\end{itemize}

\noindent {\bf Riordan graphs of the {\em Lagrange type}.} This
class of graphs is defined by a {\em Lagrange matrix} $(1,f)$, and
it is trivially related to Riordan graphs of the Bell type. Indeed,
letting $f=tg$, we see that removing the vertex 1 in $G_n(1,tg)$
gives the graph $G_{n-1}(g,tg)$ of the Bell type. Conversely, given
a graph $G_{n-1}(g,tg)$, we can always relabel each vertex $i$ by
$i+1$, and add a new vertex labelled by 1 and connected to the
vertex 2, to obtain the graph $G_n(1,tg)$ of the Lagrange type.
Examples of graphs in this class are
\begin{itemize}
\item the path graphs $P_n$ defined by $G_n(1,t)$; and
\item the {\em Fibonacci graph} $FG_n$ defined by $G_n\left(1,t+t^2\right)$.
\end{itemize}

\noindent {\bf Riordan graphs of the {\em checkerboard
type}.} This class of graphs is defined by a {\em checkerboard
matrix} $(g,f)$ such that $g$ is an even function and $f$ is an odd
function. Examples of graphs in this class are
\begin{itemize}
\item the path graphs $P_n$ defined by $G_n(1,t)$; and
\item the complete bipartite graphs $K_{\lfloor{n\over2}\rfloor, \lceil{n\over2}\rceil}$ defined by $G_n({1\over 1-t^2},t)$.
\end{itemize}

\noindent {\bf Riordan graphs of the {\em derivative type}.} This
class of graphs is defined by a {\it derivative matrix}
$(f',f)$. Examples of graphs in this class are
\begin{itemize}
\item the null graphs $N_n$ defined by $G_n(0,0)$; and
\item the path graphs $P_n$ defined by $G_n(1,t)$.
\end{itemize}

\noindent {\bf Riordan graphs of the {\em hitting time type}.} This
class of graphs is defined by a {\em hitting time matrix}
$(tf'/f,f)$, $f'(0)=1$, and it is trivially related to Riordan
graphs of the derivative type. Indeed, removing the first row and
the first column of the adjacency matrix ${\cal A}(G)$ of a graph
$G=G_n(tf'/f,f)$ of a hitting time type, we obtain the graph
$G_{n-1}(f',f)$ of the derivative type which corresponds to
removing the vertex 1 of $G$. Conversely, given a Riordan graph
$G_{n-1}(f',f)$ of the derivative type, one can relabel each vertex
$i$ by $i+1$, and add a new vertex, labeled by 1, that is connected
to the vertices defined by the coefficients of the function $tf'/f$
to obtain $G_n(tf'/f,f)$.

Thus, Riordan graphs of the Bell type (Section~\ref{Bell-type-sec}) and the derivative type (Section~\ref{derivative-sec}) are of the main interest in Section~\ref{Four-families}.

 As is mentioned above, more classes of Riordan graphs can be
introduced using the operations, ring sum $\oplus$ and R-product
$\otimes_R$ defined in Section~\ref{product-sec}. Indeed, to
illustrate this idea, note that the ring sum $\oplus$  of a Riordan
graph $G_n(g,tg)$ of the Bell type and the Riordan graph
$G_n((tg)',tg)$ of the derivative type is well defined. Such a sum
results in a new class of graphs defined by Riordan matrices of
the form $(tg',tg)$. Indeed, since $2g\equiv 0$ (mod 2), we have
$$G_n(g,tg)\oplus G_n((tg)',tg)=G_n(g+g+tg',tg)=G_n(tg',tg).$$
Note that $G_n(tg',tg)$ is not proper, as is the ring sum of any two
proper Riordan graphs.

We end the subsection by noticing that the class of Riordan graphs of the Appell type (i.e.\ Toeplitz graphs) are closed under the operations $\oplus$ and $\otimes_R$.
It is not difficult to show that this class of graphs on $n$ vertices forms a {\em commutative ring} with the identity element $G_n(1,z)$ and the zero element $G_n(0,z)=N_n$.

\subsection{The complement of a Riordan graph}\label{complement-sec}

For any Riordan graph $G_n=G_n(g,t)$ of the Appell type, the ring
sum $G_n(g,t)\oplus G_n\left(\frac{1}{1-t},t\right)$ gives the {\em
complement} of $G_n$, i.e.\ the graph in which edges of $G_n$ become
non-edges, and vice versa.

In general, it is not true that the complement of a Riordan (labeled
or unlabeled) graph is Riordan. Indeed, the complement of the star
graph $G_n\left(\frac{1}{1-t},0\right)$, $n\ge5$, is a labeled copy
of the graph $H_{n}$  in Proposition~\ref{non-Riordan-graph}, which
is non-Riordan. Thus, Riordan graphs can become non-Riordan, and
thus versa, under taking the complement.

Note that the complement of the Riordan graph $G_4(1+t^2,t^2)$ in
Fig.\ 7 is the Riordan graph $G_4(1+t+t^2,t^3)$ showing that the
operation of the complement preserves the property of being Riordan
for some graphs of non-Appell type.

\section{Structural properties of Riordan graphs}\label{struct-sec}

We begin with basic properties of a Riordan graph $G=G_n(g,f)$,
which can be directly determined in terms of column generating
functions of the binary Riordan matrix ${\cal
B}(g,f)=[b_{ij}]_{i,j\ge0}$ where $b_{ij}\equiv[t^i]gf^j\in
\{0,1\}$. Let
$$
{\cal B}_n^{(j)}={\cal B}_n^{(j)}(t):=\sum_{i=j}^n
b_{ij}t^i\in{\mathbb Z}_2[t],\;\;j=0,1,\ldots
$$
be a polynomial in $t$ of degree $n$ which generates the
(0,1)-elements $b_{ij}$ in the $j$th column of ${\cal B}(g,f)$. Thus
if ${\cal A}(G)=B+B^T$ is the $n\times n$ adjacency matrix of $G$
then $B(1|n)$ can be expressed in terms of polynomials in $t$ of
degree $n-2$:
$$
B(1|n)={\cal B}(g,f)_{n-1}:=\left[{\cal B}_{n-2}^{(0)},{\cal
B}_{n-2}^{(1)},\ldots,{\cal B}_{n-2}^{(n-2)}\right].
$$

In this paper, $d_G(i)$ denotes the {\em degree} of a vertex $i$ in
a graph $G$. If $G$ is understood from the context, we simply write
$d(i)$.

 \begin{theorem}[Basic Properties]\label{lem1}
 Let $G=G_n(g,f)$ be a Riordan graph and let ${\cal B}(g,f)=[b_{ij}]_{i,j\ge0}$. Then,
 \begin{itemize}
\item[{\rm(i)}] For $i>j\ge1$, $ i j\in E$ if and only if
$[t^{i-2}]gf^{j-1}\equiv b_{i-2,j-1}=1$.
\item[{\rm(ii)}] $ d(1)={\cal B}_{n-2}^{(0)}(1)$.
\item[{\rm(iii)}] $d(n)=\sum_{j=0}^{n-2}[t^{n-2}]{\cal B}_{n-2}^{(j)}(t)$.
\item[{\rm(iv)}] For $k\not\in\{ 1, n\}$, $d(k)={\cal B}_{n-2}^{(k-1)}(1)
+\sum_{j=0}^{k-2}[t^{k-2}]{\cal B}_{n-2}^{(j)}(t)$.
\item[{\rm(v)}] $|E|=\sum_{j=0}^{n-2}{\cal B}_{n-2}^{(j)}(1)$.
\item[{\rm(vi)}] If $G$ is proper then it has the Hamiltonian
path $1\rightarrow 2\rightarrow \cdots \rightarrow n$. In addition,
if $b_{n-2,0}=[t^{n-2}]{\cal B}_{n-2}^{(0)}(t)=1$ then $G$ has the
Hamiltonian cycle $1\rightarrow 2\rightarrow \cdots \rightarrow
n\rightarrow 1$.
\end{itemize}
\end{theorem}
\noindent{\bf Proof.} The (i)--(iii) are straightforward from
\eref{e:def2}. The (iv) follows from the fact that the degree of a
vertex $k$ is the summation of all entries located in both the
$(k-1)$th column and the $(k-2)$th row of ${\cal B}(g,f)_{n-1}$. The
(v) follows from the fact that the number of edges in $G$ is equal
to the number of 1s in ${\cal B}(g,f)_{n-1}$. The (vi) follows from
the fact that if $G$ is proper then all entries of the subdiagonal
in its adjacency matrix are 1s, i.e.\ $i$ is adjacent to $i+1$ for
$i=1,\ldots,n-1$. \hfill{\rule{2mm}{2mm}} \vskip.05pc
The {\em matching number} $\beta(G)$ is the size of a {\em
maximal matching} in a graph $G$.

  \begin{theorem}
 Let $G=G_n(g,f)$ be a proper Riordan graph. Then $\beta(G)=\lfloor\frac{n}{2}\rfloor$.
 \end{theorem}
\noindent{\bf Proof.} By (vi) in Theorem~\ref{lem1}, $G$ has the
Hamiltonian path $1\rightarrow 2\rightarrow \cdots \rightarrow n$.
Therefore, for even and odd $n$, respectively, maximal matchings are
$\{12,\, 34,\, \ldots,\, (n-1)n\}$ and $\{12,\, 34,\, \ldots,\,
(n-2)(n-1)\}$.
This completes the proof.\hfill{\rule{2mm}{2mm}}\\

From now on, the following theorems will be useful for studying the
structural properties of Riordan graphs in next sections.

\begin{lemma}\label{e:lem2}
Let $g,f\in\mathbb{Z}_2[[t]]$ where $f(0)=0$. Then, for an integer
$s\ge1$,
\begin{eqnarray}\label{e:pp}
\left(g\circ f\right)^{2^s}\equiv g\circ f^{2^s}\;\;{(\rm{mod}}\;2).
\end{eqnarray}
\end{lemma}
\textbf{Proof.} We proceed by induction on $s\ge1$. Let $s=1$ and
let $g=g(t):=\sum_{n\ge0}g_nt^n$. Since $g_n\in\{0,1\}$, $g_n^2=g_n$
for all $n\ge0$. Thus, we obtain
\begin{align*}
(g\circ f)^2&=g(f)^2=\left(\sum_{n\ge0}g_nf^n\right)^2=\sum_{n\ge0}\left(\sum_{i=0}^ng_ig_{n-i}\right)f^n\\
&=g_0^2+2g_0g_1f+(2g_0f_1+g_1^2)f^2+(2g_0g_3+2g_1g_2)f^3+\cdots\\
&\equiv\sum_{n\ge0}g_{n}^2f^{2n}\equiv\sum_{n\ge0}g_{n}(f^{2})^n=g(f^2)=g\circ
f^2,
\end{align*}
which proves \eref{e:pp} when $s=1$. Next, by inductive hypothesis
we obtain
\begin{eqnarray*}
\left(g\circ f\right)^{2^s}\equiv\left(g\circ f^2\right)^{2^{s-1}}\equiv
g\circ f^{2^s}=g\left(f^{2^s}\right).
\end{eqnarray*}
Hence the proof is complete.\hfill{\rule{2mm}{2mm}}\\

\subsection{Decomposition of Riordan graphs}\label{decomp-sec}
For any simple graph $G=(V,E)$ of order $n$ with adjacency matrix
${\cal A}(G)$ and $V=[n]$, there exists an $n\times n$ permutation
matrix $P$ such that
\begin{align}\label{permutation-product}
P\mathcal{A}(G)P^T=\left(\begin{array}{cc}
                           \mathcal{A}(\left<V_1\right>) & \tilde{B} \\
                           \tilde{B}^T & \mathcal{A}(\left<V_2\right>)
                         \end{array}
\right)
\end{align}
where $\left<V_1\right>$ and $\left<V_2\right>$ are respectively
induced subgraphs of disjoint vertex sets $V_1$ and $V_2$ in which
$V_1\bigcup V_2=V$. Theorem \ref{e:th} gives a description for the
adjacency matrix in the case of a Riordan graph where $V_1=V_o$ and
$V_2=V_e$ are assumed to be the sets of odd and even labeled
vertices, respectively.

\begin{theorem}\label{e:th}{\rm (Riordan Graph Decomposition)}  Let $G=G_n(g,f)$ be a Riordan graph of order $n$ and let $n=n_1+n_2$ where $n_1=\lceil
n/2\rceil$ and $n_2=\lfloor n/2\rfloor$. If $P$ is the $n\times n$
permutation matrix of the form
$P=\left(e_{1},e_3,\ldots,e_{2n_1-1},e_{2},e_4,\ldots,2n_2\right)^{T}$
where $e_{i}$ denotes the unit column vector with $1$ in the $i$th
position, then
\begin{eqnarray}\label{e:bm}
P\mathcal{A}(G)P^T=\left(
\begin{array}{cc}
X & \tilde{B} \\
\tilde{B}^{T} & Y
\end{array}
\right)
\end{eqnarray}
where $X={\cal A}\left(\left<V_o\right>\right)$, $Y={\cal
A}\left(\left<V_e\right>\right)$, and
\begin{itemize}
\item[{\rm(i)}] $\left<V_o\right>$ is the Riordan graph of order $n_1$ corresponding to the Riordan matrix $\left(g^{\prime}(\sqrt{t}) ,f(t)\right)$;
\item[{\rm(ii)}] $\left<V_e\right>$ is the Riordan graph of order $n_2$ corresponding to the Riordan matrix
$\left(\left(\frac{gf}{t}\right)^{\prime}(\sqrt{t}),f(t)\right)$;
\item[{\rm(iii)}] $\tilde{B}={\cal B}(t(gf)^{\prime }(\sqrt{t}),f(t))_{n_1\times n_2}+{\cal B}((tg)^{\prime}(\sqrt{t}),f(t))
_{n_2\times n_1}^{T}.$
\end{itemize}
\end{theorem}
\noindent\textbf{Proof.}  Regarding to the given permutation matrix
$P$, \eref{e:bm} is a direct consequence of
\eref{permutation-product}. If we let ${\cal A}(G)=(b_{ij})_{1\le
i,j\le n}=B+B^T$ where $B={\cal B}(tg,f)_n$ then
$$
P\mathcal{A}(G)P^T=PBP^T+(PBP^T)^T
$$ where for $m_1=2n_1-1$ and $m_2=2n_2$,
{\small\begin{eqnarray}\label{e:deco} PBP^T&=&\left(
 \begin{array}{ccccccccccc}
 0 &  & && & \vline&0&&&& \\
 b_{31} & 0 &&O& & \vline&b_{32}&0&&O&  \\
 b_{51} & b_{53} & 0 & & & \vline &b_{52}&b_{54}&0&& \\
 \vdots  & \vdots & \ddots & \ddots& & \vline &\vdots&\vdots&\ddots&\ddots&  \\
 b_{m_1,1}  & b_{m_1,3}& \cdots  &  b_{m_1,m_1-2}  &0&
 \vline&b_{m_1,2}&b_{m_1,4}&\cdots&b_{m_1,m_1-1}&0\\
 \hline
 b_{21} &  & && & \vline&0&&&& \\
 b_{41} & b_{43} &&O& & \vline&b_{42}&0&&O&  \\
 b_{61} & b_{63} &b_{65} & & & \vline &b_{62}&b_{64}&0&& \\
 \vdots  & \vdots & \ddots & \ddots& & \vline &\vdots&\vdots&\ddots&\ddots&  \\
 b_{m_2,1}  & b_{m_2,3}& \cdots  &  b_{m_2,m_2-3}  &b_{m_2,m_2-1}&
 \vline&b_{m_2,2}&b_{m_2,4}&\cdots&b_{m_2,m_2-2}&0
 \end{array} \right)\nonumber\\
 &:=&\left(\begin{array}{cc}{\hat X}_{n_1\times n_1}&Z_{n_1\times n_2}\\W_{n_2\times n_1}&{\hat Y}_{n_2\times
 n_2}\end{array}\right).
\end{eqnarray}}
In particular, if $n$ is odd then $PBP^T$ can be obtained from
\eref{e:deco} by deleting the last row and last column.

 Now we show that
${\hat X}$, ${\hat Y}$, $Z$, $W$ are Riordan matrices. First recall
that if $(g,f)=(r_{ij})_{i,j\ge0}$ where $r_{ij}=[t^i]gf^j$ then
from \eref{e:def2} and $i\ge2$ and $j\ge1$, we have
\begin{eqnarray}\label{e:moj}
b_{ij}\equiv[t^{i-1}]tgf^{j-1}\equiv[t^{i-2}]gf^{j-1}=r_{i-2,j-1}.
\end{eqnarray}

(i) Let ${\hat X}=(x_{ij})_{1\le i,j\le n_1}$. Then,
\begin{align*}
x_{ij}=\left\{
\begin{array}{ll}
b_{2i-1,2j-1}&
\text{if $1\le j<i\le n_1$} \\
0 & \text{otherwise.}
\end{array}
\right.
\end{align*}
Using \eref{e:moj} and the property $[t^i]h\equiv [t^{i-1}]h^\prime$
for any $h\in{\mathbb Z}[[t]]$ together with $g(t)^2\equiv g(t^2)$
from Lemma~\ref{e:lem2}, we obtain, for $1\leq j< i\leq n_1$,
\begin{eqnarray*}
x_{ij}&=&b_{2i-1,2j-1}\equiv
r_{2i-3,2j-2}=[t^{2i-3}]gf^{2j-2}\equiv[t^{2i-4}]\left(gf^{2j-2}\right)^\prime\\
&\equiv& [t^{2i-4}]g^\prime f^{2j-2}\equiv [t^{2i-2}]t^2g^\prime(t)
f^{2(j-1)}(t)\equiv
[t^{2(i-1)}]t^2g^\prime(t) f^{j-1}(t^2)\\
&\equiv&[t^{i-1}]tg^\prime(\sqrt t) f^{j-1}(t).
\end{eqnarray*}
Thus ${\hat X}$ is an $n_1\times n_1$ binary Riordan matrix given by
$$
{\hat X}={\cal B}\left(tg^\prime(\sqrt t), f(t)\right)_{n_1}.
$$
Since ${\cal A}\left(\left<V_o\right>\right)={\hat X}+{\hat X}^T$ it
follows that $\left<V_o\right>$ is a Riordan graph on $n_1(=\lceil
n/2\rceil)$ vertices of $\left(g^\prime(\sqrt t), f(t)\right)$.
\vskip.3pc

(ii) Let ${\hat Y}=(y_{ij})_{1\le i,j\le n_2}$. Then,
\begin{align*}
y_{ij}=\left\{
\begin{array}{ll}
b_{2i,2j}&
\text{if $1\le j<i\le n_2$} \\
0 & \text{otherwise.}
\end{array}
\right.
\end{align*}
By a similar method as (i), we obtain, for $1\leq j< i\leq n_2$,
\begin{eqnarray*}
y_{ij}&=&b_{2i,2j}\equiv r_{2i-2,2j-1}=[t^{2i-2}]gf^{2j-1}=[t^{2i-3}]\frac{gf^{2j-1}}{t}\\
&\equiv& [t^{2i-4}]\left(\frac{gf}{t}f^{2j-2}\right)^{\prime}\equiv
[t^{2i-4}]\left(\frac{gf}{t}\right)^{\prime} f^{2j-2}\equiv
[t^{2(i-2)}]\left(\frac{gf}{t}\right)^{\prime}(t) f^{j-1}(t^2)\\
&\equiv&[t^{i-1}]t\left(\frac{gf}{t}\right)^{\prime}(\sqrt{t})
f^{j-1}(t).
\end{eqnarray*}
Thus, ${\hat Y}$ is an $n_2\times n_2$ binary Riordan matrix given by
$$
{\hat Y}={\cal
B}\left(t\Big(\frac{gf}{t}\Big)^{\prime}(\sqrt{t}),f(t)\right)_{n_2}.
$$
Since ${\cal A}\left(\left<V_e\right>\right)={\hat Y}+{\hat Y}^T$ it
follows that $\left<V_e\right>$ is a Riordan graph on $n_2(=\lfloor
n/2\rfloor)$ vertices of
$\Big(\Big(\frac{gf}{t}\Big)^{\prime}(\sqrt{t}),f(t)\Big)$.
\vskip.3pc

 (iii) Let $Z=(z_{ij})_{1\le i\le n_1, 1\le j\le n_2}$.
Then,
\begin{align*}
z_{ij}=\left\{
\begin{array}{ll}
b_{2i-1,2j}&
\text{if $1\le j<i\le n_1$} \\
0 & \text{otherwise.}
\end{array}
\right.
\end{align*}
By a similar method as (i), we obtain, for $1\leq j< i\leq n_1$,
\begin{eqnarray*}
z_{ij}&=&b_{2i-1,2j}\equiv r_{2i-3,2j-1}=[t^{2i-3}]gf^{2j-1} \equiv[t^{2i-4}]\left((gf)f^{2j-2}\right)^{\prime}\\
 &\equiv&[t^{2i-4}](gf)^{\prime}f^{2j-2}\equiv [t^{2(i-2)}](gf)^{\prime}(t)f^{j-1}(t^2)\\
  &\equiv&[t^{i-1}]t(gf)^{\prime}(\sqrt{t})
  f^{j-1}(t).
\end{eqnarray*}
Thus, $Z$ is an $n_1\times n_2$ binary Riordan matrix given by
$$Z={\cal B}\left(t(gf)^{\prime}(\sqrt{t}),f(t)\right)_{n_1\times n_2}.$$

 Let $W=(w_{ij})_{1\le i\le n_2, 1\le j\le n_1}$. Then,
 \begin{align*}
w_{ij}=\left\{
\begin{array}{ll}
b_{2i,2j-1}&
\text{if $1\le j\le i\le n_2$} \\
0 & \text{otherwise.}
\end{array}
\right.
\end{align*}
 By a similar method as above, one can easily show that
 \begin{eqnarray*}
w_{ij}&=&[t^{i-1}]\,(tg)^{\prime}(\sqrt{t}) f^{j-1}(t).
\end{eqnarray*}
Thus $W$ is an $n_2\times n_1$ Riodan matrix given by
$\left((tg)^{\prime}(\sqrt{t}),f(t)\right)$. Since $\tilde{B}=Z+W^T$, the
proof is complete. \hfill{\rule{2mm}{2mm}}

\begin{definition}
{\rm Let {\em o-decomposable Riordan graphs}, standing for {\em odd
decomposable Riordan graphs}, be the class of graphs defined by
requiring $Y=O$ in (\ref{e:bm}), where $O$ is the zero matrix of
$Y$'s size. Also, let {\em e-decomposable Riordan graphs}, standing
for {\em even decomposable Riordan graphs}, be the class of graphs
defined by requiring $X=O$ in (\ref{e:bm}).}
\end{definition}

The parts (ii) and (iii) in the following theorem justify our choice
for the names of o-decomposable and e-decomposable Riordan graphs.

\begin{theorem}\label{e:th13}
Let $g\not\equiv0$ and $V_o$ and $V_e$ be, respectively, the odd and
even labelled vertex sets of a Riordan graph $G=G_{n}\left(
g,f\right)$. Then,
\begin{itemize}
\item[{\rm(i)}] For even $n$, $\left\langle V_{o}\right\rangle \cong \left\langle
V_{e}\right\rangle$ if and only if $[t^{2m-1}]g\equiv [t^{2m}]gf$
for all $m\ge1$.
\item[{\rm(ii)}] The induced subgraph $\left\langle V_{o}\right\rangle$ is
a null graph, i.e.\ $G$ is e-decomposable if and only if
$[t^{2m-1}]g\equiv 0$ for all $m\ge1$, that is, $g$ is an even
function modulo $2$.
\item[{\rm(iii)}] The induced subgraph $\left\langle V_{e}\right\rangle$ is a null
graph, i.e.\ $G$ is o-decomposable if and only if
$[t^{2m}]gf\equiv0$ for all $m\ge1$, that is, $gf$ is an odd
function modulo $2$.
\item[{\rm(iv)}] $G$ is a bipartite graph with parts $V_{o}$ and $
V_{e}$ if and only if $[t^{2m+1}]g\equiv [t^{2m}]f\equiv 0$ for all
$m\ge0$, that is, $G$ is of the checkerboard type.
\item[{\rm(v)}] There is no edge between a vertex $i\in V_o$ and a vertex $j\in V_e$ if and only if $[t^{2m}]g\equiv [t^{2m}]f\equiv
0$ for all $m\ge0$, i.e.\ $G\cong \left<V_o\right> \bigcup
\left<V_e\right>$.
\end{itemize}
\end{theorem}
\textbf{Proof.} (i) Let $n$ be even. From Theorem~\ref{e:th},
$\left\langle V_{o}\right\rangle \cong \left\langle
V_{e}\right\rangle$ if and only if
\begin{align}\label{e:eq10}
g^{\prime}\equiv \left(\frac{gf}{t}\right)^{\prime}.
\end{align}
Since $[t^{2m-1}]h'=2m[t^{2m}]h\equiv0$ and
$[t^{2m-2}]h'\equiv[t^{2m-1}]h$ for all $h\in\mathbb{Z}[[t]]$, the
equation in \eref{e:eq10} is equivalent to
$$[t^{2m-1}]g\equiv [t^{2m-1}]{gf\over t}\quad\Leftrightarrow\quad
[t^{2m-1}]g\equiv [t^{2m}]gf
$$ which proves (i).

(ii) From Theorem \ref{e:th}, the induced subgraph $\left\langle
V_{o}\right\rangle$ is a null graph if and only if the matrix $X$ in
\eref{e:bm} is a zero matrix, i.e.\
$$
 g'\equiv0 \quad
\Leftrightarrow\quad [t^{2m-1}]g\equiv 0 \text{ for all $m\ge1$}.
$$
Hence the proof follows.

(iii) From Theorem \ref{e:th}, the induced subgraph $\left\langle
V_{e}\right\rangle$ is a null graph if and only if the matrix $Y$ in
\eref{e:bm} is a zero matrix, i.e.\
\begin{align*}
 \left(\frac{gf}{t}\right)'\equiv0 \;\;
\Leftrightarrow\;\;[t^{2m}]gf\equiv 0 \text{ for all $m\ge1$}.
\end{align*}
Hence the proof follows.

(iv) From Theorem \ref{e:th}, $G$ is a bipartite graph with parts
$V_{o}$ and $V_{e} $ if and only if the matrices $X$ and $Y$ in
\eref{e:bm} are zero matrices, i.e.\
\begin{align*}
g'\equiv0\;\;\textrm{and}\;\;(gf/t)'\equiv0\quad&\Leftrightarrow\quad
g'\equiv0\;\;{\rm and}\;\;(f/t)'\equiv0\\
&\Leftrightarrow\quad \textrm{$[t^{2m+1}]g\equiv [t^{2m}]f\equiv 0$
for all $m\ge0$.}
\end{align*}
Hence the proof follows.

(v) Form Theorem \ref{e:th}, there is no edge between a vertex $i\in
V_o$ and a vertex $j\in V_e$  if and only if the matrix $B$ in
\eref{e:bm} is a zero matrix, i.e.\
\begin{align*}
(gf)^{\prime
}(\sqrt{t})\equiv0\;\;\textrm{and}\;\;(tg)^{\prime}(\sqrt{t})\equiv0\quad&\Leftrightarrow\quad
(tg)'\equiv0\;\;and\;\;(f/t)'\equiv0\\
&\Leftrightarrow\quad \textrm{$[t^{2m}]g\equiv [t^{2m}]f\equiv 0$
for all $m\ge0$.}
\end{align*}
Hence the proof follows.
 \hfill{\rule{2mm}{2mm}}


\subsection{Fractal properties of Riordan graphs}\label{A-seq-sec}
A {\em fractal} is an object exhibiting similar patterns at
increasingly small scales. Thus, fractals use the idea of a detailed
pattern that repeats itself.

%
In this section, we show that every Riordan graph $G_n(g,f)$ with
$f'(0)=1$ has fractal properties by using the notion of the
$A$-sequence of a Riordan matrix.

\begin{definition}{\rm Let $G$ be a graph. A pair of vertices $(u,v)$ in $G$ is
a {\em cognate pair} of a pair of vertices $(i,j)$ in $G$ if
\begin{itemize}
\item $|i-j|=|u-v|$ and
\item $i$ is adjacent to $j$  $\Leftrightarrow$ $u$ is adjacent to
$v$.
\end{itemize}}
\end{definition}

We denote the cognate pairs as $(i,j)\sim(u,v)$. The set of all
cognate pairs of $(i,j)$ is denoted by cog$(i,j)$. Clearly,
the relation $\sim$ on the set cog$(i,j)$ is an equivalence
relation.

\begin{definition} {\rm The $A$-sequence of the binary Riordan matrix ${\cal B}(g,f)$ defining a proper Riordan graph $G_n(g,f)$ with $f'(0)=1$ is called the {\em binary $A$-sequence of the graph}.}
\end{definition}

Let $G=G_n(g,f)$ be a proper Riordan graph. If $G$ is not a graph of
Appell type, i.e. $f\ne t$ then the binary $A$-sequence of $G$ has
of the form
\begin{eqnarray}\label{e:aaaa}
A=(\tilde{a}_k)_{k\ge0}=(1,\underbrace{0,\ldots,0}_{\ell \;{\rm
times}},1,\tilde{a}_{\ell+2},\ldots),\;\ell\ge0,\;\tilde{a}_k\in\{0,1\}.
\end{eqnarray}

The following theorem gives a relationship between cognate pairs and
the binary $A$-sequence of a proper Riordan graph.

\begin{theorem}\label{e:new} Let $G=G_n(g,f)$ be a proper Riordan graph
with a binary $A$-sequence $A=({\tilde a}_k)_{k\ge0}$ of the form
\eref{e:aaaa}. Assume that $\ell\ge0$ is the maximum number of
consecutive zero terms starting from $\tilde{a}_1$ in the
$A$-sequence. For a given vertex pair $(i,j)$ with $i>j$, if $s\ge0$
is the least integer satisfying
\begin{align}\label{e:cog1}
s\ge{\rm log}_2{i-j\over \ell+1},
\end{align}
then $(i+m2^s,j+m2^s)$ is a cognate pair of $(i,j)$ for each integer
$m$ where
\begin{align}\label{cog2}
{1-j\over 2^s}\le m\le {n-i\over2^s}.
\end{align}
 \end{theorem}
\textbf{Proof.} Let $\mathcal{A}(G)=[b_{ij}]_{1\le i,j\le n}=B+B^T$
where $b_{ij}=b_{ji}$ and $B=(tg,f)_n$. Then we may assume that
$g,f\in{\mathbb Z}_2[[t]]$ and $B=[b_{ij}]_{1\le i,j\le n}$ where
$b_{ij}=0$ if $i\le j$.
 Since the $A$-sequence of $G$ has the form
of \eref{e:aaaa}, its generating function is the polynomial over
${\mathbb Z}_2$:
$$
A(t)=1+t^{\ell+1}+\tilde{a}_{\ell+2}t^{\ell+2}+\cdots+\tilde{a}_{n-3}t^{n-3}.
$$
Using $f\equiv tA(f)$ it follows from Lemma \ref{e:lem2} that if
$k=2^s$ is an integer for $s\ge0$ then
 \begin{eqnarray*}\label{cognate-equation}
b_{i+k,j+k}&\equiv&\left[t^{i+k-2}\right]
gf^{j+k-1}\equiv\left[t^{i+k-2}\right]
gf^{j-1}\left(tA(f)\right)^{k}\equiv\left[t^{i-2}\right]gf^{j-1}A(f^k)\\
&=&\left[t^{i-2}\right]gf^{j-1}\left(1+f^{k(\ell+1)}+\sum_{m={\ell+2}}^{n-3}\tilde{a}_mf^{km}\right)\\
&=&b_{ij}+b_{i,j+k(\ell+1)}+\sum_{m={\ell+2}}^{n-3}\tilde{a}_mb_{i,j+km}.\end{eqnarray*}
Since $B=[b_{ij}]_{1\le i,j\le n}$ is a lower triangular matrix with
zero diagonal, if $j+k(\ell+1)\ge i$ then $b_{i,j+k(\ell+1)}=0$ so
that $b_{i,j+km}=0$ for all $m\ge{\ell+2}$. Thus, if $s$ satisfies
\eref{e:cog1} then
$$
b_{i+k,j+k}\equiv b_{i,j},
$$
i.e. $(i,j)\sim(i+2^s,j+2^s)$. Moreover, by transitivity of cognate
pairs we have $(i,j)\sim(i+m2^s,j+m2^s)$ where $m$ is an integer
such that $1\le j+m2^s <i+m2^s\le n$, i.e. ${1-j\over 2^s}\le m\le
{n-i\over2^s}$ for the least integer satisfying \eref{e:cog1}. Hence
we complete the proof.
 \hfill{\rule{2mm}{2mm}}
\vskip.5pc

 In particular, let $\ell=0$. If $i$ is adjacent to $j$
then the pairs cognate with $(i,j)$ are those connected by edges in
the following figures:

  \begin{center}
  \begin{tabular}{cc}
\scalebox{0.7} 
{
\begin{pspicture}(0,0)(12.2,1.01)
\psdots[dotsize=0.2](4.08,-0.09) \psdots[dotsize=0.2](6.08,-0.09)
\usefont{T1}{ptm}{m}{n}
\rput(4.1,-0.43){$i$}
\usefont{T1}{ptm}{m}{n}
\rput(6.08875,-0.43){$j$}
\psarc[linewidth=0.04](5.08,0.0){0.99}{0.0}{180.0}
\psdots[dotsize=0.2](2.08,-0.09) \psdots[dotsize=0.2](8.08,-0.09)
\psarc[linewidth=0.04](3.08,0.0){0.99}{0.0}{180.0}
\psarc[linewidth=0.04](7.08,0.0){0.99}{0.0}{180.0}
\psdots[dotsize=0.2](10.08,-0.09)
\psarc[linewidth=0.04](9.1,0.0){0.99}{0.0}{180.0}
\psdots[dotsize=0.08](10.68,-0.09)
\psdots[dotsize=0.08](11.08,-0.09)
\psdots[dotsize=0.08](11.48,-0.09) \psdots[dotsize=0.08](0.68,-0.09)
\psdots[dotsize=0.08](1.08,-0.09) \psdots[dotsize=0.08](1.48,-0.09)
\usefont{T1}{ptm}{m}{n}
\rput(2.1,-0.43){$i-2^s$}
\usefont{T1}{ptm}{m}{n}
\rput(8.1,-0.43){$j+2^s$}
\usefont{T1}{ptm}{m}{n}
\rput(10.1,-0.43){$j+2\cdot2^s$} \psdots[dotsize=0.2](4.08,-0.09)
\psdots[dotsize=0.2](6.08,-0.09) \psdots[dotsize=0.2](4.08,-0.09)
\psdots[dotsize=0.2](6.08,-0.09) \psdots[dotsize=0.2](2.08,-0.09)
\psdots[dotsize=0.2](8.08,-0.09) \psdots[dotsize=0.2](10.08,-0.09)
\end{pspicture}
}& for $|i-j|=2^s$\\
\scalebox{0.6} 
{
\begin{pspicture}(0,0)(18.2,2)
\psdots[dotsize=0.2](6.08,-0.09) \psdots[dotsize=0.2](8.08,-0.09)
\usefont{T1}{ptm}{m}{n}
\rput(6.1,-0.43){$i$}
\usefont{T1}{ptm}{m}{n}
\rput(8.1,-0.43){$j$}
\psarc[linewidth=0.04](7.095,0.0){0.99}{0.0}{180.0}
\psdots[dotsize=0.2](4.08,-0.09) \psdots[dotsize=0.2](10.08,-0.09)
\psdots[dotsize=0.2](12.08,-0.09) \psdots[dotsize=0.08](16.68,-0.09)
\psdots[dotsize=0.08](17.08,-0.09)
\psdots[dotsize=0.08](17.48,-0.09) \psdots[dotsize=0.08](0.68,-0.09)
\psdots[dotsize=0.08](1.08,-0.09) \psdots[dotsize=0.08](1.48,-0.09)
\usefont{T1}{ptm}{m}{n}
\rput(4.1,-0.43){$j-2^s$}
\usefont{T1}{ptm}{m}{n}
\rput(10.1,-0.43){$i+2^s$}
\usefont{T1}{ptm}{m}{n}
\rput(12.1,-0.43){$j+2^s$} \psdots[dotsize=0.2](2.08,-0.09)
\usefont{T1}{ptm}{m}{n}
\rput(2.1,-0.43){$i-2^s$} \psdots[dotsize=0.2](14.08,-0.09)
\usefont{T1}{ptm}{m}{n}
\rput(14.1,-0.43){$i+2\cdot2^s$} \psdots[dotsize=0.2](16.08,-0.09)
\usefont{T1}{ptm}{m}{n}
\rput(16.046719,-0.43){$j+2\cdot2^s$}
\psarc[linewidth=0.04](3.075,0.0){0.99}{0.0}{180.0}
\psarc[linewidth=0.04](11.075,0.0){0.99}{0.0}{180.0}
\psarc[linewidth=0.04](15.075,0.0){0.99}{0.0}{180.0}
\psdots[dotsize=0.2](6.08,-0.09) \psdots[dotsize=0.2](8.08,-0.09)
\psdots[dotsize=0.2](4.08,-0.09) \psdots[dotsize=0.2](10.08,-0.09)
\psdots[dotsize=0.2](12.08,-0.09) \psdots[dotsize=0.2](2.08,-0.09)
\psdots[dotsize=0.2](14.08,-0.09) \psdots[dotsize=0.2](16.08,-0.09)
\end{pspicture}
}& for $|i-j|< 2^s$
  \end{tabular}
  \end{center}
\vskip.3pc
  \centerline{Fig.\ 9: Cognate pairs in a Riordan graph}
\vskip.5pc

\begin{remark} {\rm We note that not every cognate pairs of $(i,j)$
can be expressed as the form $(i+m2^s,j+m2^s)$ for some integers $m$
and $s\ge0$. It might be interesting to find all cognate pairs of a
vertex pair $(i,j)$ of a Riordan graph.}
\end{remark}

The following theorem shows that every proper Riordan graph
$G_n(g,f)$ with $f\ne t$ has a fractal property.

\begin{theorem}\label{e:coro2} Let $G=G_n(g,f)$ be the same Riordan
graph in Theorem~\ref{e:new}. Consider a cognate pair
$(i+m2^s,j+m2^s)$ of a vertex pair $(i,j)$ with $i>j$. Then, the
associated two induced subgraphs are isomorphic as follows:
\begin{eqnarray*}
\left<\{j,j+1,\ldots,i\}\right>\cong
\left<\{j+m2^s,j+1+m2^s,\ldots,i+m2^s\}\right>.
\end{eqnarray*}
\end{theorem}
 \textbf{Proof.} Let $X=\{j,j+1,\ldots,i\}$ and $Y=\{j+m2^s,j+1+m2^s,\ldots,i+m2^s\}$.
 Consider a bijection $\sigma:\,X\rightarrow Y$ with
 $\sigma(k)=k+m2^s$ for $k\in X$. It is enough to show that for any two vertices $a,b\in X$ where $a>b$,
 $(a,b)\sim(\sigma(a),\sigma(b))$.
 Consider a vertex pair $(a,b)$ in $X$. Clearly, $|a-b|=|\sigma(a)-\sigma(b)|$.
 By Theorem~
\ref{e:new}, if $t\ge0$ is the least integer satisfying
$t\ge\log_2{a-b\over \ell+1}$ then there exists an integer $\tilde
m$ such that $a+\tilde{m}2^t$ and $b+\tilde{m}2^t$ where $\tilde m$
is an integer that satisfies
\begin{align}\label{cog-fractal-2}
{1-b\over2^t} \le \tilde{m}\le {n-a\over2^t}.
\end{align}
We claim that if we take $\tilde{m}=m2^{s-t}$ then
$(a,b)\sim(\sigma(a),\sigma(b))$. Indeed, since $j<b<a<i$ it follows
from \eref{cog2} that
$${1-b\over2^t}\le{1-j\over2^t} \le m2^{s-t}\le {n-i\over2^t}\le
{n-a\over2^t}.
$$ Thus $\tilde{m}:=m2^{s-t}$ satisfies \eref{cog-fractal-2} which
implies that $a+\tilde{m}2^t=a+m2^{s}=\sigma(a)$ and
$b+\tilde{m}2^t=b+m2^{s}=\sigma(b)$. Thus
$(a,b)\sim(\sigma(a),\sigma(b))$, which completes the proof.
\hfill{\rule{2mm}{2mm}}

\begin{example} {\rm Consider cognate pairs of $(1,5)$ in the Catalan graph $G=G_n(C,tC)$.
Since the binary $A$-sequence of $G$ is $(1,1,\ldots)$, we have
$\ell=0$ so that $s=2$. By Theorem \ref{e:new}, every vertex pair
with the form of $(1+4m,5+4m)$ is a cognate pair of $(1,5)$ where
$m$ is an integer $0\le m\le {n-5\over4}$. For example, such pairs
are $(1,5),(5,9)$ for $n=9$, and $(1,5),(5,9),(9,13)$ for $n=15$.
Furthermore, if $n=15$ then by Theorem~\ref{e:coro2} we have
\begin{align*}
\left<\{1,2,3,4,5\}\right>\cong\left<\{5,6,7,8,9\}\right>\cong\left<\{9,10,11,12,13\}\right>.
\end{align*}
Fig.\ 10 shows a way to draw $CG_9$ in a ``fractal form''.}
\end{example}

\begin{center}
\scalebox{0.5} 
{
\begin{pspicture}(0,-0.55)(16.284687,8.063)
\psarc[linewidth=0.025999999](7.09,-0.04){3.01}{0.0}{180.0}
\psarc[linewidth=0.025999999](11.13,-0.08){5.03}{0.0}{180.0}
\psarc[linewidth=0.025999999](10.14,0.01){6.04}{0.0}{180.0}
\psarc[linewidth=0.025999999](9.15,-0.02){7.03}{0.0}{180.0}
\psarc[linewidth=0.025999999](8.13,0.0){8.05}{0.0}{180.0}
\psarc[linewidth=0.025999999,linecolor=Red](1.1,0.05){0.98}{0.0}{180.0}
\psarc[linewidth=0.025999999,linecolor=Red](3.1,0.05){0.98}{0.0}{180.0}
\psarc[linewidth=0.025999999,linecolor=Red](5.1,0.05){0.98}{0.0}{180.0}
\psarc[linewidth=0.025999999,linecolor=Red](7.1,0.05){0.98}{0.0}{180.0}
\psarc[linewidth=0.025999999,linecolor=Red](2.1,0.03){2.0}{0.0}{180.0}
\psarc[linewidth=0.025999999,linecolor=Red](6.1,0.03){2.0}{0.0}{180.0}
\psarc[linewidth=0.025999999,linecolor=Red](5.1,-0.09){3.02}{0.0}{180.0}
\psarc[linewidth=0.025999999,linecolor=Red](4.11,0.01){4.03}{0.0}{180.0}
\psarc[linewidth=0.025999999,linecolor=Red](9.12,0.05){0.98}{0.0}{180.0}
\psarc[linewidth=0.025999999,linecolor=Red](11.12,0.05){0.98}{0.0}{180.0}
\psarc[linewidth=0.025999999,linecolor=Red](13.12,0.05){0.98}{0.0}{180.0}
\psarc[linewidth=0.025999999,linecolor=Red](15.12,0.05){0.98}{0.0}{180.0}
\psarc[linewidth=0.025999999,linecolor=Red](10.12,0.03){2.0}{0.0}{180.0}
\psarc[linewidth=0.025999999,linecolor=Red](14.12,0.03){2.0}{0.0}{180.0}
\psarc[linewidth=0.025999999,linecolor=Red](13.12,-0.09){3.02}{0.0}{180.0}
\psarc[linewidth=0.025999999,linecolor=Red](12.13,0.01){4.03}{0.0}{180.0}
\psdots[dotsize=0.24](0.12,-0.01) \psdots[dotsize=0.24](2.12,-0.01)
\psdots[dotsize=0.24](4.12,-0.01) \psdots[dotsize=0.24](6.12,-0.01)
\psdots[dotsize=0.24](8.12,-0.01) \psdots[dotsize=0.24](10.12,-0.01)
\psdots[dotsize=0.24](12.12,-0.01)
\psdots[dotsize=0.24](14.12,-0.01)
\psdots[dotsize=0.24](16.14,-0.01)
\psdots[dotsize=0.08](16.8,-0.01)\psdots[dotsize=0.08](17.1,-0.01)\psdots[dotsize=0.08](17.4,-0.01)
\usefont{T1}{ptm}{m}{n}
\rput(0.15,-0.4){\large$1$}
\usefont{T1}{ptm}{m}{n}
\rput(2.15,-0.4){\large$2$}
\usefont{T1}{ptm}{m}{n}
\rput(4.15,-0.4){\large$3$}
\usefont{T1}{ptm}{m}{n}
\rput(6.15,-0.4){\large$4$}
\usefont{T1}{ptm}{m}{n}
\rput(8.15,-0.4){\large$5$}
\usefont{T1}{ptm}{m}{n}
\rput(10.15,-0.4){\large$6$}
\usefont{T1}{ptm}{m}{n}
\rput(12.15,-0.4){\large$7$}
\usefont{T1}{ptm}{m}{n}
\rput(14.15,-0.4){\large$8$}
\usefont{T1}{ptm}{m}{n}
\rput(16.2,-0.4){\large$9$}
\end{pspicture}
}
\end{center}
\vskip.3pc
 \centerline{Fig.\ 10: The Catalan graph $G_9(C,tC)$}

%
%
%
%

\subsection{Reverse relabeling of Riordan graphs}\label{revRio}

Consider a relabeling of a graph $G$ with $n$ vertices labeled by
$1,2,\ldots, n$. The relabeling can be done by {\em reversing} the
vertices in $[n]$, that is, by replacing a label $i$ by $n+1-i$ for
each $i\in[n]$. Let ${\cal A}(G)=B+B^T$ for some $n\times n$
$(0,1)$-lower triangular matrix with zero diagonal. Then, the
adjacency matrix of the resulting relabeled graph $G^{\prime}$ can
be obtained from $E{\cal A}(G)E=EBE+EB^TE$ where $E$ is the $n\times
n$ backward identity matrix corresponding to the reverse relabeling
given by
\begin{align*}
E=\left(
    \begin{array}{cccc}
      0 & \cdots & 0 & 1 \\
      0 & \cdots & 1 & 0 \\
      \vdots & \iddots & \vdots & \vdots \\
      1 & \cdots & 0 & 0 \\
    \end{array}
  \right).
\end{align*}
Since $EB^TE$ is a lower triangular matrix with zero diagonal, the
adjacency matrix ${\cal A}(G^{\prime})$ can be expressed by the {\it
flip-transpose} \cite{HJ} of $B$ defined as $B^F:=EB^TE$:
\begin{eqnarray}\label{e:labeling}
{\cal A}(G^{\prime})=B^F+(B^F)^T.
\end{eqnarray}

Relabeling a Riordan graph does not necessarily result in a Riordan
graph. In this section, we show that if $G=G_n(g,f)$ is a proper
Riordan graph then the graph resulting in relabeling done by
reversing the vertices in $[n]$ is a proper Riordan graph.

\begin{proposition}[\cite{CJ}]\label{e:CJK}
Let $L=(g,f)_n$ be an $n\times n$ leading principal matrix of a
proper Riordan matrix $(g,f)$. Then, the flip-transpose $L^F=EL^TE$
is a proper Riordan matrix given by
\begin{align*}
L^F=(g(\bar f)\cdot \bar f'\cdot(t/\bar f)^n,\bar f)_n
\end{align*}
\end{proposition}

Proposition~\ref{e:CJK} implies the following theorem.

\begin{theorem}\label{e:reverse}
The reverse relabelling of a Riordan graph $G=G_n(g,f)$ with
$f'(0)=1$ is a Riordan graph given by $G^\prime=G_n(g(\bar f)\cdot
\bar f'\cdot(t/\bar f)^{n-1},\bar f).$
\end{theorem}

From \eref{e:labeling} we have
\begin{align}\label{e:BF}
B^F={\cal B}(tg(\bar f)\cdot \bar f'\cdot(t/\bar f)^{n-1},\bar f)_n.
\end{align}

To give an example for Theorem~\ref{e:reverse} we need the following lemma.

\begin{lemma}{\rm (Lucas Theorem)}\label{e:lucas}
Let $m$ and $n$ be nonnegative integers with the base $p$ (a prime)
expansions:
\begin{itemize}
\item $m=m_kp^k+m_{k-1}p^{k-1}+\cdots+m_1p+m_0$ and
\item $n=n_kp^k+n_{k-1}p^{k-1}+\cdots+n_1p+n_0$.
\end{itemize}
Then, the following congruence relation holds:
$$
{m\choose n}\equiv \prod_{i=0}^k{m_i\choose
n_i}\;\;{(\rm{mod}}\;p).$$
\end{lemma}

\begin{example}{\rm Consider the Catalan graph
$G=G_{n}(C,tC)$ when $n=2^j-1$ with $j\ge2$. Let
$g=C={1-\sqrt{1-4t}\over 2t}$ and $f=tC$. By simple computations we
obtain $\bar f=t+t^2\equiv t-t^2$ and $g(\bar f)=C(\bar f)={1\over
1-t}$. Thus,
\begin{align*}
g(\bar f)\cdot \bar f'\cdot(t/\bar
f)^{2^j-2}&\equiv{1\over 1-t}\cdot(1-2t)\cdot\left({1\over
1-t}\right)^{2^j-2}\equiv\left({1\over 1-t}\right)^{2^j-1}.
\end{align*}
Let ${\tilde g}=\sum_{k=0}^{2^j-3}[t^k]1/(1-t)^{2^j-1}$. We claim
that $\tilde g\equiv 1-t$. By Newton's binomial theorem we obtain
$$
\left({1\over 1-t}\right)^{2^j-1}=\sum_{k=0}^\infty{2^j+k-2\choose
k}t^k\equiv 1-t+\sum_{k=2}^\infty{2^j+k-2\choose k}t^k.
$$
Now, let $k-2=a_i2^i+\cdots+a_22^2+a_12+a_02^0$ where
$a_i\in\{0,1\}$. Since $0\le k\le2^j-3$, we have $j\ge i$. Thus, we
obtain
\begin{eqnarray*}
2^j+k-2&=&2^j+a_i2^i+\cdots+a_22^2+a_12+a_02^0\\
k&=&a_i2^i+\cdots+a_22^2+(a_1+1)2+a_02^0.
\end{eqnarray*}
If $a_1=0$ then ${a_1\choose a_1+1}={0\choose 1}=0$. Let $a_1=1$.
Then, we can assume that $a_{\ell}=0,a_{\ell-1}=\cdots=a_1=1$ where
$i\ge\ell-1$. Then ${a_{\ell}\choose a_{\ell}+1}={0\choose 1}=0.$ By
Lucas theorem, for $k\ge2$,
\begin{align}\label{By-Lucas-Thm1}
{2^j+k-2\choose k}\equiv0.
\end{align}
Thus, $\tilde g=1-t$ and by Theorem \ref{e:reverse}, the reverse
relabeling of the Catalan graph $G_{2^j-1}(C,tC)$ is a Riordan graph
corresponding to the Riordan matrix $(1-t,t-t^2)$. We note that
$(1-t,t-t^2)$ is the inverse matrix of $(C,tC)$. Fig.\ 12
illustrates, as fractals, the adjacency matrices of the two graphs.}
\end{example}

\begin{center}
\epsfig{file=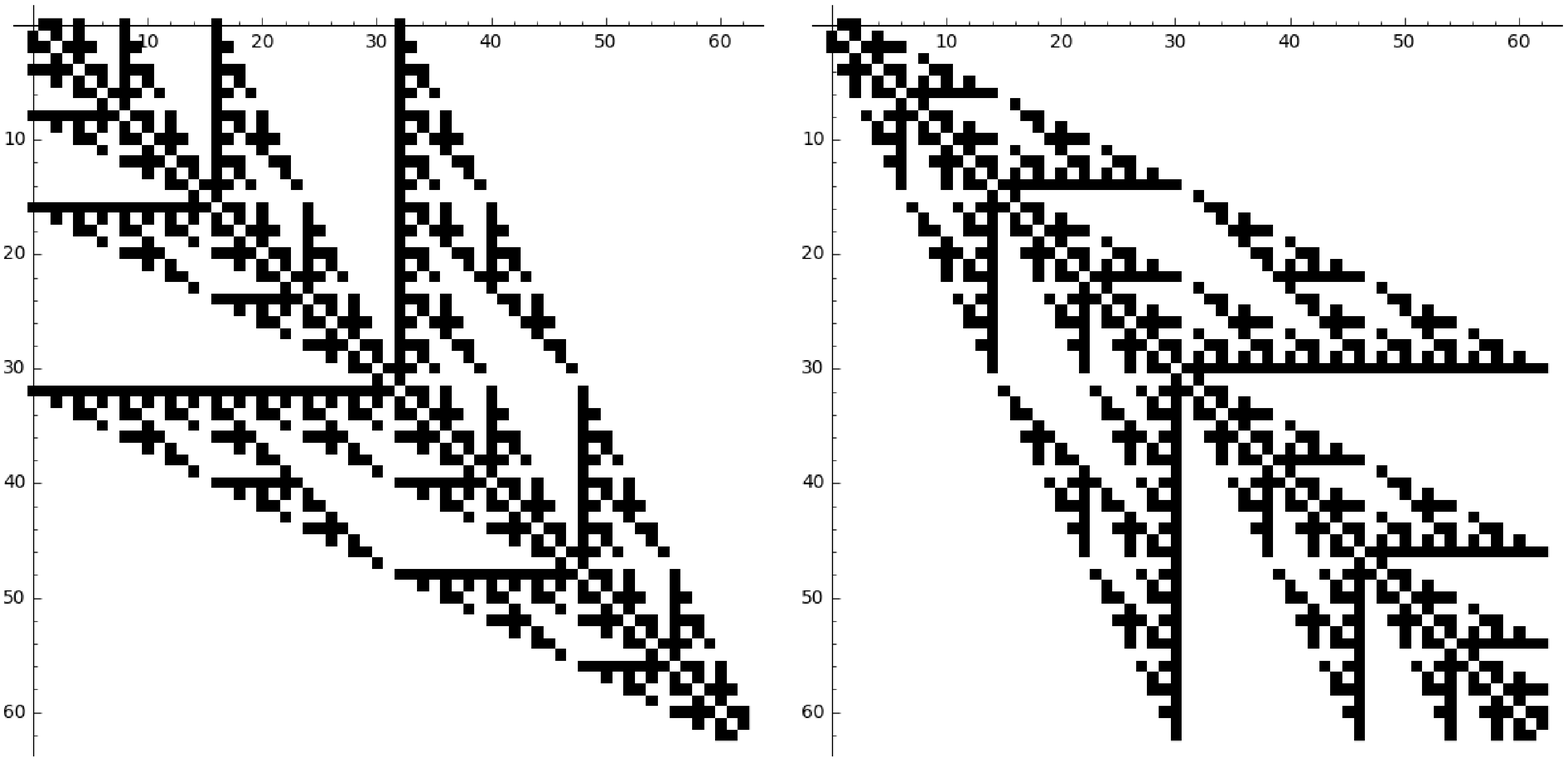,scale=0.3}
\end{center}
\centerline{Fig.\ 12: The fractal natures of a Catalan graph and its
reverse labeling}

%

\subsection{Eulerian and
Hamiltonian Riordan graphs}\label{eh-sec}

In this section we give certain conditions on Riordan graphs to have
an Eulerian trail/cycle or a Hamiltonian cycle. As the result, we
obtain large classes of Riordan graphs which are Eulerian or
Hamiltonian.

\begin{theorem}\label{vertex-degree} Let $G=G_n(g,f)$ be a
Riordan graph with $f'(0)=1$. Then, the degree $d(i)$ of a vertex $i$
is odd (resp., even) if and only if
\begin{align*}
\phi_{i-1}+\varphi_{n-i}\equiv 1\;\;\textrm{(resp., $0$)}
\end{align*}
where for $j=0,\ldots,n-1$
\begin{align*}
\phi_j\equiv[t^j]{tg\over 1-f}\quad\textrm{and}\quad
\varphi_j\equiv[t^j]{tg(\bar f)\cdot(\bar f)'\cdot(t/\bar
f)^{n-1}\over 1-\bar f}.
\end{align*}
\end{theorem}
\noindent\textbf{Proof.} Let $\mathcal{A}(G)=B+B^T$, ${\bold
d}=(d(1),\ldots,d(n))^T$ and
$\mathbf{1}=(1,1,\ldots,1)^T\in\mathbb{Z}^n$ where $B={\cal
B}(tg,f)_n$. Since
$B^T\mathbf{1}=(EB^FE)\mathbf{1}=(EB^F)\mathbf{1}$ where $B^F$ is
given by \eref{e:BF}, applying the fundamental theorem given by
\eref{e:ftbrm} to $B$ and $B^F$, respectively we obtain
\begin{align*}
{\bold
d}&=\mathcal{A}(G)\mathbf{1}=\left(B+B^T\right)\mathbf{1}=B\mathbf{1}+E(B^F\mathbf{1})\\
&\equiv\left(\begin{array}{c}
  \phi_0\\
  \phi_1\\
  \vdots \\
  \phi_{n-1}
\end{array}\right)+\left(\begin{array}{c}
  \varphi_{n-1}\\
  \varphi_{n-2} \\
  \vdots \\
  \varphi_{0}
\end{array}\right)\equiv\left(\begin{array}{c}
  \phi_0+\varphi_{n-1}\\
  \phi_1+\varphi_{n-2}\\
  \vdots \\
  \phi_{n-1}+\varphi_{0}
\end{array}\right).
\end{align*}
Thus, $d(i)$ is odd (resp., even) for
$i=1,\ldots,n$ if and only if $\phi_{i-1}+\varphi_{n-i}\equiv 1$ (resp., $0$),
as required.
 \hfill{\rule{2mm}{2mm}}\\

Let $G$ be a simple connected graph. Since $G$ has an Eulerian trail if and only if $G$
has no odd vertex or exactly two odd vertices, by Theorem~\ref{vertex-degree} we obtain the following corollary.

\begin{corollary}\label{euler-coro}
Let $G=G_n(g,f)$ be a proper Riordan graph with $f'(0)=1$. Then, $G$
has an Eulerian trail if and only if $\phi_{i-1}$~and
$\varphi_{n-i}$ satisfy either for all $i=1,\ldots,n$,
\begin{align}\label{eulerian-cycle}
\phi_{i-1}+\varphi_{n-i}\equiv 0
\end{align}
or for $i_1\neq i_2$
\begin{align*}
\phi_{i-1}+\varphi_{n-i}\equiv\left\{
\begin{array}{ll}1
 & \text{if $i=i_1,i_2\in\{1,\ldots,n\}$} \\
0 & \text{if $i\in\{1,\ldots,n\}\setminus\{i_1,i_2\}$}%
\end{array}
\right..
\end{align*}
In particular, if \eqref{eulerian-cycle} is satisfied then $G$ is
Eulerian.
\end{corollary}

\begin{example}
{\rm The Pascal graph $PG_n=G_n\left({1\over 1-t},{t\over
1-t}\right)$ for $n=2^k+1$ ($k\ge2$) has an Eulerian trail. We note
that $PG_n$ is a proper Riordan graph so that it is a connected
graph. Indeed, since $\bar f={t\over 1+t}$, simple computations show
that
\begin{align*}
\phi_{i-1}&\equiv[t^{i-1}]{t\over 1-2t}\equiv[t^{i-1}]t,\\
\varphi_{n-i}&\equiv[t^{n-i}]t(1+t)^{n-1}={n-1\choose i}.
\end{align*}
By Lucas Theorem, if $n=2^k+1$ and $i=2^k$ then ${n-1\choose
i}\equiv1$, and $0$ otherwise. Thus,
\begin{align*}
\phi_{i-1}+\varphi_{n-i}\equiv\left\{
\begin{array}{ll}
1 & \text{if $i=2$ or $i=2^k$,} \\
0 & \text{if $i\in\{1,\ldots,2^k+1\}\setminus\{2,2^k\}$}.
\end{array}
\right.
\end{align*}
By Corollary \ref{euler-coro}, $PG_{2^k+1}$ has an Eulerian trail.
Similarly, if $n=2^k$ for $k\ge3$ then it is shown that $PG_{n}$
does not have an Eulerian trail.}
\end{example}

A polynomial $P=\sum_{k=0}^na_kt^k$ of degree $n$ or a sequence
$(a_k)_{k=0}^n$ is called {\it palindromic} if $a_k=a_{n-k}$ for
$k=0,1,\ldots,n$.
\begin{corollary}\label{euler-coro2}
Let $G=G_n(g,t)$ be a proper Riordan graph of the Appell type with
$f'(0)=1$ and $n\ge3$ vertices where
$g=\sum_{k=0}^{n-2}g_kt^k\in{\mathbb Z}_2[t]$. If $g$ is palindromic
with $g(1)\equiv 0$ then $G$ is Eulerian.
\end{corollary}
\noindent\textbf{Proof.} Since $f=t$, it can be easily shown from
Theorem~\ref{vertex-degree} that
\begin{align*}
\phi_j=\varphi_j\equiv[t^j]{tg\over 1-t}=\left\{
\begin{array}{ll}0 & \text{if $j=0$} \\
\sum_{k=0}^{j-1}g_k & \text{if $j=1,\ldots,n-1$}.
\end{array}
\right.
\end{align*}
Since $g_0+g_1+\cdots+g_{n-2}\equiv 0$ and $g_k\equiv g_{n-2-k}$ for
$k=0,\ldots,n-2$, we obtain $\phi_0=\phi_{n-1}$ and for
$j=1,\ldots,n-2$,
\begin{align*}
\phi_j\equiv\sum_{k=0}^{j-1}g_k\equiv\sum_{k=j}^{n-2}g_k\equiv\sum_{k=j}^{n-2}g_{n-2-k}\equiv\phi_{n-1-j}.
\end{align*}
Thus, the sequence $\left(\phi_k\right)_{k=0}^{n-1}$ is palindromic.
Since $\phi_k=\varphi_k$ it follows from Corollary~\ref{euler-coro}
that $G$ is Eulerian. \hfill{\rule{2mm}{2mm}\\

Hamiltonian properties of Toeplitz graphs, which are Riordan graphs
of the Appell type, have been studied in~\cite{DTT}. Next, we give
more general results. Recall that every proper Riordan graph
$G=G_n(g,f)$ has the Hamiltonian path
 $1\rightarrow 2 \rightarrow \cdots \rightarrow n$. In particular, if $g=\frac{1}{1-t}$ then $G$ is
 Hamiltonian, i.e.\ $G$ has a Hamiltonian cycle.

\begin{theorem}\label{th1}
 Let $G=G_n(g,f)$ be a proper Riordan graph. If one of the following holds then $G$ is Hamiltonian.

 {\rm(i)} There exists $i\in\{2,\ldots,n-1\}$ such that $[t^{i-1}] g\equiv 1$ and $[t^{n-2}] (gf^{i-1})\equiv 1$.

 {\rm(ii)} $[t]g\equiv 1$ and $[t^2] f \equiv 0$.
 \end{theorem}
\noindent\textbf{Proof.}
 (i) Let ${\cal A}(G)=[a_{ij}]_{1\le i,j\le n}$. Since $G$ is proper, we have the path $1\rightarrow 2 \rightarrow \cdots \rightarrow i$. If $[t^{n-2}] (gf^{i-1})\equiv 1$ for some $i\ne1$,
 i.e. $a_{ni}=1$ then we have the edge $in\in E$. Again, since $G_n$ is proper,
 we have the path $n\rightarrow n-1 \rightarrow \cdots \rightarrow i+1$. Finally, if $[t^{i-1}] g\equiv 1$ for some $i\ne1$, i.e. $a_{i+1,1}=1$ then $(i+1)1\in E$.
 Thus we have the following Hamiltonian cycle in $G$:
 $$1\rightarrow 2 \rightarrow \cdots \rightarrow i\rightarrow n \rightarrow n-1 \rightarrow \cdots \rightarrow i+1 \rightarrow 1$$
 as required.

 (ii) Since $G$ is proper and $[t]g\equiv1$, by Theorem~\ref{e:th} we obtain that ${\cal A}(\langle V_o\rangle)$ is proper, i.e.\
 we have the path $1 \rightarrow 3 \rightarrow \cdots \rightarrow 2\left\lceil\frac{n}{2}\right\rceil-1$. Our assumption that $[t^2] f\equiv 0$ implies $[t^2] (gf)\equiv 1$. By Theorem~\ref{e:th} we obtain that ${\cal A}(\langle V_e\rangle)$ is proper, i.e.\
 we have the path $2\rightarrow 4 \rightarrow  \cdots \rightarrow 2\left\lfloor\frac{n}{2}\right\rfloor$. Finally, we obtain the following Hamiltonian cycle in
 $G$:
  $$1\rightarrow 3 \rightarrow \cdots \rightarrow 2\left\lceil\frac{n}{2}\right\rceil-1 \rightarrow 2\left\lfloor\frac{n}{2}\right\rfloor\rightarrow\cdots \rightarrow 4 \rightarrow 2 \rightarrow 1$$
  as required.
 \hfill{\rule{2mm}{2mm}}

\begin{example}
{\rm The Catalan graph $G=G_n(C,tC)$ when $n=2^k+1$ is Hamiltonian.
Indeed, it is known \cite{DS} that $C_n=[t^n] C\equiv 1$ if and only
if $n=2^k-1$ for $k\ge0$:
$$ C=\sum_{n=0}^\infty{1\over n+1}{2n\choose n}t^n\equiv 1+t+t^3+t^7+t^{15}+\cdots.
$$
 Take $i=n-1$ in (i) of Theorem~\ref{th1}.
If $n=2^k+1$ then we obtain
$$
[t^{n-2}]C=[t^{2^k-1}]C \equiv 1\;{\rm and}\;[t^{n-2}]
C(tC)^{n-2}=[t^0]C^{n-1}\equiv 1.
$$
 Thus the Catalan graph $G_{2^k+1}(C,tC)$ is Hamiltonian.}
 \end{example}

The following result is obtained by Theorem~3.1 in \cite{DS} and
Proposition~1 in \cite{AAB}.
\begin{lemma} \label{lem4}
 The Motzkin number $M_n$ is even if and only if either $n\in S_1$ or $n\in S_2$ where
 $S_1=\{4^i(2j-1)-2\;|\;i,j\ge1\}$ and $S_2=\{4^i(2j-1)-1\;|\;i,j\ge1\}$.
\end{lemma}

\begin{theorem}  \label{th1-1}
 If $n\neq 4^i(2j-1)$ for $i,j\ge 1$ then the Motzkin graph
$G=G_n(M,tM)$ is Hamiltonian for $n\ge 3$.
\end{theorem}
\noindent\textbf{Proof.} Consider the generating function $M$ of the
Motzkin numbers $M_n$:
$$M=\frac{1-t-\sqrt{1-2t-3t^2}}{2t^2}=\sum_{n\ge0}M_nt^n=1+t+2t^2+4t^3+9t^4+21t^5+51t^6+\cdots.$$
 Take $i=n-1$
in (i) of Theorem~\ref{th1}. Note that $S_1\cap S_2=\emptyset$ in
Lemma~\ref{lem4}. By Lemma~\ref{lem4} we obtain
$$
[t^{n-2}]M=M_{n-2}\equiv 1\;{\text{ if and only if} }\; n\notin S_3
\cup S_4
$$
where $S_3:=\{ 4^i(2j-1)\;|\;i,j\ge1\}$ and
$S_4:=\{4^i(2j-1)+1\;|\;i,j\ge1\}$. In addition,
$[t^{n-2}]M(tM)^{n-2}=[t^0]M^{n-1}\equiv 1$. Thus by
Theorem~\ref{th1}, if $n\notin S_3 \cup S_4$ then $MG_n$ is
Hamiltonian. It remains to prove that if $n\in S_4$ then $MG_n$ is
Hamiltonian. To show this, consider the case of $i=2$ in (i) of
Theorem~\ref{th1}. Since
$$
[t^{n-2}]tM^2(t)\equiv [t^{n-2}]tM(t^2)=[t^{n-2}]\sum_{k\ge 0} M_k
t^{2k+1}
$$
it follows from $n-2=2k+1$ that if $n\in S_4$ then
\begin{eqnarray}\label{e:kk}
k\in\{4^{i-1}(4j-2)-1\;|\;i,j\ge1\}=\{1,5,7,9,23,\ldots\}.
\end{eqnarray}
 Since $k\notin S_1\cup S_2$ for any $k$ in \eref{e:kk}, we have $M_k\equiv 1$ by Lemma \ref{lem4} so that
 $[t^{n-2}]tM^2(t)\equiv 1$. In addition,
$[t]M=M_1\equiv1$. Thus by Theorem~\ref{th1}, if $n\in S_4$ then $G$
is Hamiltonian as required. \hfill{\rule{2mm}{2mm}}
\bigskip

A complete split graph $CS_{m,\,n-m}$, $m\leq n$, is a graph on $n$
vertices consisting of a clique $K_{m}$ on $m$ vertices and a stable
set (i.e.\ independent set) on the remaining $n-m$ vertices, such
that any vertex in the clique is adjacent to each vertex in the
stable set. The bipartite graph $G(m,n-m)$ obtained from
$CS_{m,\,n-m}$ by deleting all edges in $K_{m}$ is referred to as
the bipartite graph corresponding to $CS_{m,\,n-m}$.

\begin{lemma}[\cite{BH}]\label{lemma1}
 In a split graph $CS_{m,\,n-m}$, if $m<n-m$ then
 $CS_{m,\,n-m}$ is not Hamiltonian. If
 $m=n-m$ then $CS_{m,\,n-m}$ is Hamiltonian if
 and only if
 the corresponding bipartite graph $G({m,\,n-m})$ is Hamiltonian.  \end{lemma}

 \begin{theorem}\label{th3}
If $G=G_n(g,f)$ is an improper Riordan graph with $[t]f\equiv 0$
then $G_n$ is  not Hamiltonian.
 \end{theorem}
\noindent\textbf{Proof.} If $f$ satisfies $[t]f\equiv 0$ then the
Riordan graph $G_n(g,f)$ for any $g$ is a subgraph of
$G_n(\frac{1}{1-t}, t^2)$. Thus it is sufficient to prove the claim
for $G_n(\frac{1}{1-t}, t^2)$. Let ${\cal A}(G)=[a_{i,j}]_{1\le
i,j\le n}$ and let $C_j(t)$ be the $j$th column generating function
of ${\cal B}(tg,f)$ where $g=\frac{1}{1-t}$ and $f=t^2$.
 We may assume that $i>j$ and let $n=n_1+n_2$ where $n_1=\lceil
n/2\rceil$ and $n_2=\lfloor n/2\rfloor$. Then $C_{n_2+1}(t)=g\cdot
t^{2n_2}$, i.e.\
 $a_{i,j}=0$ if $j\ge n_2+1$. Thus the subgraph of $G$ induced by $\langle \{n_2+1, \ldots, n
 \}\rangle$ is a null graph. Clearly, $G$ is a subgraph of the complete split graph $CS_{n_2, n_1}$. If $n$ is odd then by Lemma~\ref{lemma1} $CS_{n_2, n_1}$ is not Hamiltonian.
 It follows that $G$ is not Hamiltonian. Otherwise, $n$ is even. Again, by Lemma~\ref{lemma1}, $CS_{n_2,n_1}$ is Hamiltonian if the corresponding bipartite graph $G({n_2,n_1})$ is
 Hamiltonian. But in this case
 $C_{n_2}(t)=g\, t^{n-2}$. Now if  $g_0\equiv 1$, then $a_{n/2,(n-2)}\equiv
 1$ and otherwise if $g_0\equiv 0$ then $a_{n/2,(n-2)}\equiv 0$. Thus the vertex $n/2$ has maximum degree 1 in the graph $G(n_2,n_1)$, i.e.\ $G(n_2,n_1)$ is not Hamiltonian. Again, by Lemma~\ref{lemma1}, $CS_{\lfloor n/2 \rfloor, \lceil n/2\rceil
 }$ has no Hamiltonian cycle so that $G_n$ is not Hamiltonian.
 \hfill{\rule{2mm}{2mm}}\\

\begin{proposition}
 Let $G=G_n(g,f)$ be an e-decomposable/o-decomposable Riordan graph and $n$ be even. Then $G$ is Hamiltonian if and only if the bipartite graph with partitions $V_o$ and $V_e$ is Hamiltonian.
 Moreover, if $G$ is an e-decomposable Riordan graph and $n$ is odd then $G$ is not Hamiltonian.
 \end{proposition}
\noindent\textbf{Proof.}
 Let $G$ be an e-decomposable/o-decomposable Riordan graph and $n$ be even. In this case, $G$ contains an independent set on $n/2$ vertices and so $G$ is a subgraph of the complete split graph $CS_{n/2,n/2}$.
 By Lemma~\ref{lemma1},  $CS_{n/2,n/2}$ is Hamiltonian if and only if the corresponding bipartite graph $G({n/2},\,{n/2})$ is Hamiltonian.

If $G$ is an e-decomposable Riordan graph and $n$ is odd, then $G$ is a subgraph of the complete split graph $CS_{(n-1)/2,\,(n+1)/2}$. By Lemma~\ref{lemma1} we obtain the desired result. \hfill{\rule{2mm}{2mm}}\\

We end this section by observing that every Riordan graph
$G=G_n(g,f)$ of the checkerboard type of odd order $n$ is not
Hamiltonian since $G$ is bipartite of odd order by (iv) in Theorem
\ref{e:th13}.

\section{Four families of Riordan graphs}\label{Four-families}

In this section, we define io-decomposable and ie-decomposable
Riordan graphs, and also consider Riordan graphs of the Bell type
and of the derivative type.

\subsection{io-decomposable and ie-decomposable Riordan graphs}\label{ioie-sec}

\begin{definition}\label{faithful-def} {\rm Let $G=G_{n}(g,f)$ be a proper Riordan graph
 with the odd and even vertex sets $V_{o}$ and $V_{e}$, respectively.
\begin{itemize}
\item If $\left< V_{o}\right> \cong G_{\lceil n/2\rceil}(g,f)$ and $\left<
V_{e}\right>$ is a null graph then $G_n$ is {\em io-decomposable}.
\item If $\left< V_{o}\right>$ is a null graph and $\left< V_{e}\right> \cong
G_{\lfloor n/2\rfloor}(g,f)$ then $G_n$ is {\em ie-decomposable}.
\end{itemize}
``io'' and ``ie'' stand for ``isomorphically odd'' and
``isomorphically even'', respectively.}\end{definition}

\begin{theorem}\label{e:th12}
Let $G=G_{n}(g,f)$ be a proper Riordan graph.
\begin{itemize}
\item[{\rm(i)}] $G$ is io-decomposable if and only if $g'\equiv g^2$ and $gf\equiv
t\cdot(f/t)'.$
\item[{\rm(ii)}] $G$ is ie-decomposable if and only if $g^\prime\equiv 0$ and $t^2g\equiv
tf'+f.$
\end{itemize}
\end{theorem}

\textbf{Proof.} (i) Since $g(t^2)\equiv g^2(t)$ by
Lemma~\ref{e:lem2}, by the definitions and Theorem \ref{e:th}, $G$
is io-decomposable if and only if the matrix $X$ in \eref{e:bm} is
given by
\begin{align*}
\mathcal{A}\left(G_{\lceil n/2\rceil}(g^{\prime}(\sqrt{t})
,f(t))\right)=\mathcal{A}\left(G_{\lceil
n/2\rceil}(g(t),f(t))\right)\quad\Leftrightarrow\quad
g'(t)=g(t^2)\equiv g^2(t)
\end{align*}
and
\begin{align*}
G_{\lfloor n/2\rfloor
}\left(\left(\frac{gf}{t}\right)'(\sqrt{t}),f(t)\right)\cong
N_{\lfloor n/2\rfloor }\quad&\Leftrightarrow\quad 0\equiv
\left(\frac{gf}{t}\right)'\equiv g^2{f\over
t}+g\cdot\left({f\over t}\right)'\\
&\Leftrightarrow\quad gf\equiv t\cdot(f/t)'.
\end{align*}

(ii) Similarly, $G$ is ie-decomposable if and only if the matrix $X$
in \eref{e:bm} is given by
\begin{align*}
G_{\lceil n/2\rceil}(g^{\prime}(\sqrt{t}) ,f(t))\cong N_{\lceil
n/2\rceil}\quad\Leftrightarrow\quad g'(t)\equiv 0
\end{align*}
and
\begin{align*}
&\mathcal{A}\left(G_{\lfloor n/2\rfloor
}\left(\left(\frac{gf}{t}\right)'(\sqrt{t}),f(t)\right)\right)=\mathcal{A}\left(G_{\lfloor
n/2\rfloor }(g(t),f(t))\right)\\
&\Leftrightarrow\quad g^2\equiv \left(\frac{gf}{t}\right)'\equiv
g'{f\over t}+g\cdot\left({f\over t}\right)'\equiv
g\cdot\left({f\over t}\right)'\\
&\Leftrightarrow\quad t^2g\equiv tf'+f.
\end{align*}
\hfill{\rule{2mm}{2mm}}

\subsection{Riordan graphs of the Bell type}\label{Bell-type-sec}

In this section, we consider some properties of a Riordan graph
$G_n(g,tg)$ of Bell type with odd and even vertex sets $V_{o}$ and $
V_{e}$, respectively.

\begin{theorem}\label{Bell-is-o-decomp}
Any Riordan graph $G=G_n(g,tg)$ of the Bell type is $o$-decomposable
in which its adjacency matrix is permutational equivalent to
\begin{eqnarray*}
\left(
\begin{array}{cc}
X & \tilde{B} \\
\tilde{B}^{T} & O%
\end{array}
\right)
\end{eqnarray*}
where $X$ is the adjacency matrix of $\left< V_{o}\right>=G_{\lceil
n/2\rceil}(g^{\prime}(\sqrt{t}),tg(t))$ and
\begin{eqnarray}\label{e:bm1}
\tilde{B}=\mathcal{B}(tg,tg)_{\lceil n/2\rceil ,\lfloor
n/2\rfloor}+\mathcal{B}((tg)^{\prime}(\sqrt{t}),tg)) _{\lfloor
n/2\rfloor ,\lceil n/2\rceil }^{T}.
\end{eqnarray}.
\end{theorem}
\noindent\textbf{Proof.} Since for $m\ge1$
\begin{align*}
[t^{2m}]t\cdot g^2(t)=[t^{2m-1}]g^2(t)\equiv [t^{2m-1}]g(t^2)\equiv
0,
\end{align*}
by (iii) of Theorem \ref{e:th13} $G_n$  is o-decomposable. Its
adjacency matrix immediately follows from by Theorem~\ref{e:th}.
\hfill{\rule{2mm}{2mm}} \vskip.3pc

Now, we consider io-decomposable Riordan graphs of the Bell type. We
first derive conditions on the $A$-sequences of such graphs.

\begin{lemma}\label{e:lem3}
A Riordan graph $G=G_{n}(g,tg)$ is io-decomposable if and only if
\begin{eqnarray*}
g^{2}\equiv g^{\prime}, \text{ i.e.\ $[t^j]g\equiv[t^{2j+1}]g$.}
\end{eqnarray*}
\end{lemma}
\textbf{Proof.} Since $\left<V_e\right>$ is a null graph, by
Theorem~\ref{e:th12} $G$ is io-decomposable if and only if $g'\equiv
g^2$. \hfill{\rule{2mm}{2mm}}

\begin{remark}\label{remark-io} {\rm Let $[t^j]g=g_j$ with $g_0=1$. By Lemma
\ref{e:lem3}, a Riordan graph $G_{n}(g,tg)$ is io-decomposable if
and only if $g_{2m}\equiv g_{(2m+1)2^k-1}$ for each $m\ge0$ and
$k\ge1$, i.e.
\begin{align*}
g=\sum_{n\ge0}g_{2n}\left(\sum_{k\ge0} t^{(2n+1)2^k-1}\right).
\end{align*}
Since $g$ depends  only  on its even coefficients, the number of
io-decomposable Riordan graphs of the Bell type is equal to
$2^{\lfloor n/2\rfloor-1}$.}
\end{remark}

\begin{theorem}\label{e:th11}
A Riordan graph $G=G_n(g,tg)$ is io-decomposable if and only if its
binary $A$-sequence~is of the form
$\left(1,1,a_2,a_2,a_4,a_4,\ldots\right)$, $a_{2j}\in\{0,1\}$ whose
generating function is
\begin{align*}
A(t)\equiv(1+t)+(1+t)\sum_{j\ge1}a_{2j}t^{2j}.
\end{align*}
\end{theorem}
\textbf{Proof.} Let $G$ be io-decomposable. Since there is a unique
binary generating function $A(t)=\sum_{i\ge0}a_it^i\in{\mathbb
Z}_2[t]$ such that $g\equiv A(tg)$, by applying derivative to both
sides, we obtain $g'\equiv (g+tg')\cdot A'(tg).$ Since $g^{2}\equiv
g^{\prime}$ by Lemma \ref{e:lem3}, it implies
\begin{align}\label{e:eq9}
g\equiv(1+tg)\cdot A'(tg),\ \text{i.e.}\ A(tg)\equiv(1+tg)\cdot
A'(tg).
\end{align}
Since $A'(t)\equiv \sum_{i\ge 0}a_{2i+1}t^{2i}$, the equation
\eref{e:eq9} is equivalent to
\begin{align*}
\sum_{j\ge0}a_j(tg)^j&\equiv(1+tg)\left(\sum_{i\ge
0}a_{2i+1}(tg)^{2i}\right)\\
&= a_1+a_1tg+a_3(tg)^2+a_3(tg)^3+a_5(tg)^4+a_5(tg)^5+\cdots.
\end{align*}
Thus $a_{2i}\equiv a_{2i+1}$  for $i\ge0$ as desired.
\hfill{\rule{2mm}{2mm}}

\begin{corollary}\label{e:th1}
Let $A(t)$ be a generating function of the binary $A$-sequence for a
proper Riordan graph $G=G_n(g,tg)$. Then we have:
\begin{itemize}
\item[{\rm(i)}] If $A(t)\equiv \sum_{j=0}^{2s+1}t^{n}$ or $A(t)\equiv \left( 1+t\right) ^{2t-1}$ then $G$ is
io-decomposable.
\item[{\rm(ii)}] If $A(t)\equiv \sum_{j=0}^{2s}t^{n}$ or $A(t)\equiv \left( 1+t\right) ^{2t}$ then $G$ is not io-decomposable.
\end{itemize}
\end{corollary}
\textbf{Proof.} (i) If $A(t)\equiv \sum_{j=0}^{2s+1}t^{n}$ then
clearly $G$ is io-decomposable by Theorem \ref{e:th11}. Now let
$A(t)=(1+t)^{2t-1}.$ Since $(g,tg)$ is of the Bell type, we obtain
\begin{eqnarray*}
g\equiv A(tg)=(1+tg)^{2t-1}.
\end{eqnarray*}
Applying derivative to both sides and using $2t-1\equiv1$, we obtain
\begin{eqnarray*}
g^{\prime } \equiv (g+tg^{\prime }) (1+tg)^{2t-2}=\frac{g+tg^{\prime
}}{1+tg}(1+tg)^{2t-1}\equiv\frac{g+tg^{\prime }}{1+tg}g.
\end{eqnarray*}
It implies that
\begin{eqnarray*}
g^{\prime }(1+tg) \equiv (g+tg^{\prime }) g,\ \text{i.e.}\ g^{\prime
}\equiv g^{2}.
\end{eqnarray*}
By Lemma \ref{e:lem3}, $G$ is io-decomposable.

(ii) If $A\equiv \sum_{j=0}^{2s}t^{n}$ then clearly $G$ is not
io-decomposable by Theorem~\ref{e:th11}. Now let $A=(1+t)^{2t}$.
Using $g= A(tg) \equiv (1+tg)^{2t}$ we obtain $g^{\prime }\equiv 0$.
Since $g^{2}(t) \equiv g(t^{2}) \not\equiv 0$ we have $g^{\prime
}\not\equiv g^{2}$. Thus, $G$ is not io-decomposable.
\hfill{\rule{2mm}{2mm}}

\begin{example}\label{ex3} {\rm Note that the Pascal graph $PG_{n}$, the\
Catalan graph $CG_{n}$ and the Mottkin graph $MG_{n}$ have the
generating functions $1+t$, $(1+t)^{-1}$ and $1+t+t^{2}$ for binary
$A$-sequences, respectively.
By Corollary~\ref{e:th1} we see that $PG_{n}$ and $CG_{n}$ are
io-decomposable, but $MG_{n}$ is not io-decomposable.}
\end{example}

\begin{lemma}\label{e:lem}
If $G=G_{n}(g,tg)$ is io-decomposable then the number $m(G)$ of
edges of $G$ is given by
\begin{eqnarray*}
m(G) =2m(G_{\lceil n/2\rceil}) +m(H_{\lfloor n/2\rfloor+1})
\end{eqnarray*}
where $H_{\lfloor n/2\rfloor +1}\cong G_{\lfloor n/2\rfloor
+1}\left(( tg)^{\prime }(\sqrt{t}) ,tg\right).$
\end{lemma}

\noindent \textbf{Proof.} Let $\sigma(A_n)$ denote the number of
$1$s in the adjacency matrix $A_n=\mathcal{A}\left(G\right)$. Since
$G$ is io-decomposable, by~\eref{e:bm1} we obtain
\begin{eqnarray}\label{e:m1}
\sigma(A_n)=\sigma(A_{\lceil n/2\rceil })+ 2\sigma(B_1)+2\sigma(B_2)
\end{eqnarray}
where $B_1=\mathcal{B}((tg)^{\prime }(\sqrt{t}),tg) _{\lfloor
n/2\rfloor\times\lceil n/2\rceil}$ and
$B_2=\mathcal{B}(tg,tg)_{\lceil n/2\rceil\times\lfloor n/2\rfloor}$.
We can see that
\begin{eqnarray*}
2\sigma(B_1)=\sigma\left(\mathcal{A}\left(H_{\lfloor n/2\rfloor
+1}\right)\right)\;\;\textrm{and}\;\;2\sigma(B_2)=\sigma(A_{\lceil
n/2\rceil}).
\end{eqnarray*}
Applying this to \eref{e:m1}, we obtain
\begin{eqnarray*}
\sigma(A_n)=2\sigma(A_{\lceil n/2\rceil })+
\sigma\left(\mathcal{A}\left(H_{\lfloor n/2\rfloor +1}\right)\right)
\end{eqnarray*}
which implies the desired result.
\hfill{\rule{2mm}{2mm}}\\

Since $d_{G_n}(n)=m(G_{n})-m(G_{n-1})$, by Lemma \ref{e:lem} we obtain the
following lemma.
\begin{lemma}\label{e:lem1}
If a Riordan graph $G_{n}=G_{n}\left( g,tg\right) $ is
io-decomposable then
\begin{itemize}
\item[{\rm(i)}] $d_{G_{n}}(n)=2\left\{ m\left( G_{k}\right) -m\left(
G_{k-1}\right)\right\}=2\;d_{G_{k}}(k)$ if $n=2k-1$
\item[{\rm(ii)}]
$m\left(H_{k+1}\right)-m\left(H_{k}\right)=d_{H_{k+1}}\left(k+1\right)$
if $n=2k$
\end{itemize}
where $m(G_{0})=0$ and $H_{j}\cong G_{j}((tg)^{\prime }(
\sqrt{t}),tg).$
\end{lemma}

\begin{example} {\rm
(a) Consider the Catalan graph $G_n=G_n(C,tC)$ which is
io-decomposable. Since $C=1+tC^2$ we obtain
\begin{align*}
(tC)'={d\over dt}\left(t+(tC)^2\right)\equiv
1+2(tC)(tC)^\prime\equiv1.
\end{align*}
Thus it follows from Lemmas~\ref{e:lem} and~\ref{e:lem1} that
\begin{align*}
m(G_{n}) =2\;m\left(G_{\left\lceil n/2\right\rceil }\right) +m\left(
G_{\left\lfloor n/2\right\rfloor }\right) +1.
\end{align*}
Equivalently, we obtain
\begin{align*}
m\left( G_{n}\right)=\left\{
\begin{array}{ll}
2\;m\left(G_{k}\right)+m\left(G_{k-1}\right) +1 &
\text{if $n=2k-1$} \\
3\;m\left(G_{k}\right) +1 & \text{if $n=2k$.}%
\end{array}%
\right.
\end{align*}
This recurrence relation gives
\begin{align*}
m\left(G_{2^{k}}\right) =\frac{3^{k}-1}{2}\quad\text{and}\quad
m\left(G_{2^{k}+1}\right) =\frac{3^{k}-1}{2}+2^{k}.
\end{align*}
(b) Let $G_n=G_n\left({1\over 1-t},{t\over 1-t}\right)$ be the
Pascal graph which is io-decomposable. Since $(tg)'=\left({t\over
1-t}\right)'\equiv {1\over 1-t^2}$, by Lemma \ref{e:lem} we obtain
\begin{align*}
m(PG_{n}) =2m(PG_{\lceil n/2\rceil}) +m(PG_{\lfloor n/2\rfloor+1}).
\end{align*}
This recurrence relation gives
\begin{align*}
m(PG_{2^k+1})=3^k\quad\textrm{and}\quad m(PG_{2^k})=3^k-2^k
\end{align*}
which also known as in \cite{SGB}.}
\end{example}

\begin{definition} {\rm A vertex in a graph $G$ is {\em universal} (or {\em apex}, or {\em dominating vertex}) if it is adjacent to all other vertices in $G$.} \end{definition}

It is known \cite{DQ} that the Pascal graphs $PG_n$ which are
io-decomposable have two or three universal vertices for $n\ge2$.
The following theorem gives a result in this direction for any
io-decomposable Riordan graphs of the Bell type with $2^i+1$ or
$2^i+2$ vertices.

\begin{theorem}\label{e:th4}
Let $G_n=G_{n}(g,tg) $ be io-decomposable. If $n=2^{i}+1$ for $i\geq
0$, then $G_{n}$ and $G_{n+1}$ have at least one universal vertex,
namely the vertex $2^i+1$.
\end{theorem}
\noindent \textbf{Proof.} Let $n=2^{i}+1$. It is enough to show that
$d_{G_{n}}(n)=2^{i}.$ We prove this by induction on $i\geq 0.$ Let
$i=0.$ Since $g(0)=1,$ $d_{G_{2}}(2)=1.$ Thus it holds for $i=0.$
Let $i\geq 1.$ Then,
\begin{eqnarray*}
d_{G_{n}}(n)&=&2\{ m(G_{2^{i-1}+1})-m(
G_{2^{i-1}})\}\quad\quad\textrm{(by
Lemma \ref{e:lem1})} \\
&=&2\{2m(G_{2^{i-2}+1})+m(H_{2^{i-2}+1})-2m(G_{2^{i-2}})-m(
H_{2^{i-2}+1})\}\quad\textrm{(by Lemma \ref{e:lem})} \\
&=&2^2\{m(G_{2^{i-2}+1})-m(G_{2^{i-2}})\} \\
&=&2\;d_{G_{2^{i-1}+1}}(2^{i-1}+1)\quad\quad\textrm{(by Lemma \ref{e:lem1})} \\
&=&2^{i}\quad\quad\textrm{(by the induction hypothesis)}.
\end{eqnarray*}
Thus if $n=2^{i}+1$ then the vertex $n$ is a universal vertex of
$G_{n}$. In addition, two vertices $n$ and $n+1$ are adjacent in
$G_{n+1}$ since $G_{n+1}$ is proper. Thus the vertex $n$ is also a
universal vertex of $G_{n+1}$ if $n=2^{i}+1$.
\hfill{\rule{2mm}{2mm}}

\begin{theorem}\label{e:th5}
An io-decomposable Riordan graph $G_{n}(g,tg)$ is
$(\lceil\log_2n\rceil+1)$-partite.
\end{theorem}
\textbf{Proof.} We proceed by induction on $n\ge2$. Let $n=2$. Since
$G_{2}(g,tg)$ is clearly bi-partite, the theorem holds for $n=2$.
Let $n\ge3$. Since $\left< V_{o}\right> \cong G_{\lceil
n/2\rceil}(g,tg)$ is io-decomposable, and $\left<V_{e}\right>$ is a
null graph,  by the induction hypothesis, $G_{n}(g,tg)$ is the
$(\lceil\log_2{\lceil n/2\rceil}\rceil+2)$-partite graph. Now it is
enough to show that $\lceil\log_2{\lceil
n/2\rceil}\rceil=\lceil\log_2n\rceil-1$. For all $k\ge0$, when
$2^k<n\le2^{k+1}$ we have
\begin{eqnarray*}
\lceil\log_2{\lceil n/2\rceil}\rceil=k=\lceil\log_2n\rceil-1.
\end{eqnarray*}
Hence we obtain the desired result. \hfill{\rule{2mm}{2mm}}

\begin{remark}
{\rm We note that if $G_{n}(g,tg)$ is io-decomposable then it is
$(\lceil\log_2n\rceil+1)$-partite with partitions
$V_1,V_2,\ldots,V_{\lceil\log_2n\rceil+1}$ such that
\begin{align*}
V_j=\left\{2^{j-1}+1+i2^j\;|\;0\le i\le\left\lfloor
{n-1-2^{j-1}\over 2^j}\right\rfloor\right\} \quad \text{if $1\le
j\le\lceil\log_2n\rceil$}
\end{align*}
and  $V_{\lceil\log_2n\rceil+1}=\{1\}$.}
\end{remark}

Riordan arrays frequently arising in combinatorics are of Bell type.
So it would be interesting to study graph invariants for Riordan
graphs of the Bell type.

\begin{definition}{\rm A {\em clique} is a subset of vertices of a graph $G$ such that its
induced subgraph is a complete graph. The {\em clique number} of $G$
is the number of vertices in a maximum clique in $G$, and it is
denoted by $\omega(G)$.} \end{definition}

\begin{theorem}\label{e:th6}
For $n\geq 1$, if $G=G_n(g,tg)$ is io-decomposable then
\begin{eqnarray*}
\omega(G)=\lceil\log_2n\rceil+1.
\end{eqnarray*}
\end{theorem}
\textbf{Proof.} It follows from Theorem \ref{e:th5} that
\begin{eqnarray}\label{e:eq3}
\omega(G)\le\lceil\log_2n\rceil+1.
\end{eqnarray}
Let $\tilde{V}=\{1\}\cup\{2^i+1\;|\;0\le
i\le\lceil\log_2n\rceil-1\}\subseteq V$. By Theorem \ref{e:th4}, for
every $i$, $0\le i\le\lceil\log_2n\rceil$, the vertex
$2^{i}+1\in\tilde{V}$ is adjacent to all vertices in
$\{1\}\cup\{2^{j}+1\;|\;0\le j\le i-1\}$. Thus the induced subgraph
of $\tilde{V}$ is the complete graph $K_{\lceil\log_2n\rceil+1}$. By
\eref{e:eq3} we obtain the desired result.
\hfill{\rule{2mm}{2mm}}\\

Since the complete graph $K_5$ is not planar, from Theorem
\ref{e:th6} we immediately obtain the following corollary.
\begin{corollary}
An io-decomposable graph $G_n(g,tg)$ is not planar for all $n\ge9$.
\end{corollary}

It is known \cite{DQ} that the Pascal graph $PG_n$ is planar for
$n\le7$ but it is not for $n\ge8$. Also we may check that the
Catalan graph $CG_n$ is planar for $n\le8$ but it is not for
$n\ge9$.

\begin{definition} {\rm The {\it chromatic number} of a graph $G$ is the smallest number of
colors needed to color the vertices of $G$ so that no two adjacent
vertices share the same color, and it is denoted by $\chi(G)$.}
\end{definition}

\begin{theorem}\label{chromatic-thm}
For $n\geq 1$, if a Riordan graph $G=G_n(g,tg)$ is io-decomposable
then
\begin{eqnarray*}
\chi(G)=\lceil\log_2n\rceil+1.
\end{eqnarray*}
\end{theorem}
\textbf{Proof.} Since a $k$-partite graph is $k$-colorable and
$\chi(G)\ge\omega(G)$, we obtain the desired result by
Theorems~\ref{e:th5} and~\ref{e:th6}. \hfill{\rule{2mm}{2mm}}

\begin{definition} {\rm The {\em distance} between two vertices $u,v$ in a graph $G$ is the
number of edges in a shortest path between $u$ and $v$, and it is
denoted by $d(u,v)$. The {\em diameter} of $G$ is the maximum
distance between all pairs of vertices, and it is denoted by
diam$(G)$.} \end{definition}

It is obvious that if $G$ has a universal vertex then
$\mathrm{diam}(G)=1$ or $2$.

\begin{theorem}\label{e:th9}
If a Riordan graph $G=G_n(g,tg)$ is io-decomposable then
\begin{align}\label{e:dia}
\mathrm{diam}(G)\le \lfloor{\rm log}_2n\rfloor.
\end{align}
In particular, if $n=2^k+2$ or $2^{k+1}+1$, for $k\ge1$, then
$\mathrm{diam}(G)=2.$
\end{theorem}
\textbf{Proof.} We proceed by induction on $n\ge2$. Let $n=2$. Since
clearly $\mathrm{diam}(G)=1\le \lfloor{\rm log}_22\rfloor$, the
statement is true for $n=2$. Suppose $n\ge 3$. Let $V_1=\{i\;|\;1\le
i\le 2^b+1\}$ and $V_2=\{i\;|\; 2^b+1\le i\le n\}$ where $b$ is an
integer such that $2^{b}< n\le2^{b+1}$. By Theorem \ref{e:th4},
$2^b+1$ is a universal vertex in the induced subgraph
$\left<V_1\right>$. Thus, $\mathrm{diam}(\left<V_1\right>)\le 2$.
From Theorem~\ref{e:coro2} and Lemma \ref{e:th11}, we obtain
$\left<V_2\right>\cong G_{n-2^b}(g,tg).$ By induction hypothesis, we
have
\begin{align}\label{e:V2}
\mathrm{diam}(\left<V_2\right>)\le \lfloor{\rm
log}_2(n-2^b)\rfloor\le \lfloor{\rm log}_2n\rfloor-1.
\end{align}
Let $u\in V_1\backslash \{2^b+1\}$ and $v\in V_2\backslash
\{2^b+1\}$. Now it is enough to show that $d(u,v)\le \lfloor{\rm
log}_2n\rfloor$. Since $2^b+1$ is an universal vertex of
$\left<V_1\right>$, it follows from \eref{e:V2} that
\begin{align*}
d(u,v)\le d(u,2^b+1)+d(2^b+1,v)\le\lfloor{\rm log}_2n\rfloor
\end{align*}
which proves \eref{e:dia}. Now let $n=2^k+2$ or $2^{k+1}+1$ for
$k\ge1$. By Theorem \ref{e:th4}, every io-decomposable Riordan graph
$G=G_n(g,tg)$ has at least one universal vertex. Thus
$\mathrm{diam}(G)$ is 1 or 2. Now it is enough to show that $G$ is
not a complete graph for $n\ge4$. Let ${\cal A}(G)=[a_{i,j}]_{1\le
i,j\le n}$. Since $a_{4,2}\equiv[t^2]tg^2(t)\equiv[t]g(t^2)=0$, $G$
cannot be $K_n$ for $n\ge4$. Hence we obtain the desired result.
 \hfill{\rule{2mm}{2mm}}

\begin{corollary}
Let $n\ge6$ and a Riordan graph $G=G_n(g,tg)$ be io-decomposable. If
$2^k+1<n<2^{k+1}$ then
$$\mathrm{diam}(G)\le \lfloor{\rm log}_2(n-2^k)\rfloor+1.$$
\end{corollary}
\textbf{Proof.} Let $V_1=\{i\;|\;1\le i\le 2^k+1\}$ and
$V_2=\{i\;|\; 2^k+1\le i\le n\}$. By Theorem~\ref{e:coro2} and
Lemma~\ref{e:th11}, we obtain $\left<V_2\right>\cong
G_{n-2^k}(g,tg).$ Thus, by \eref{e:dia} we have
\begin{align}\label{e:V2-1}
\mathrm{diam}(G_{n-2^k})\le \lfloor{\rm log}_2(n-2^k)\rfloor.
\end{align}
Let $u\in V_1\backslash \{2^k+1\}$ and $v\in V_2\backslash
\{2^k+1\}$. Since the vertex $2^k+1$ is a universal vertex in the
induced subgraph $\left<V_1\right>$, it follows from \eref{e:V2-1}
that
\begin{align*}
d(u,v)\le d(u,2^k+1)+d(2^k+1,v)\le\lfloor{\rm
log}_2(n-2^k)\rfloor+1,
\end{align*}
which completes the proof. \hfill{\rule{2mm}{2mm}}

\subsection{Riordan graphs of the derivative type}\label{derivative-sec}

We now consider some properties of a Riordan graph
$G_n(f',f)$ of the derivative type.

\begin{theorem}\label{deriv-is-e-decomp}
Any Riordan graph $G=G_n(f',f)$ of the derivative type is
$e$-decomposable in which its adjacency matrix is permutational
equivalent to
\begin{align*}
\left(
\begin{array}{cc}
O & \tilde{B} \\
\tilde{B}^{T} & Y%
\end{array}%
\right)
\end{align*}
where $Y$ is the adjacency matrix of $\left< V_{e}\right>=G_{\lfloor
n/2\rfloor }\left(\left(f'f/t\right)^{\prime}(\sqrt{t}),f(t)\right)$
and
\begin{align*}
\tilde{B}=\mathcal{B}(tf'(t),f(t))_{\lceil n/2\rceil\times\lfloor
n/2\rfloor}+\mathcal{B}(f'(\sqrt{t}),f(t))_{\lfloor
n/2\rfloor\times\lceil n/2\rceil}^{T}.\end{align*}
\end{theorem}
\noindent\textbf{Proof.} Let $f=\sum_{i\ge1}f_it^i$. Since
$[t^{2m-1}]f'=2mf_{2m}\equiv0$ for $m\ge1$, it follows from (ii) in
Theorem \ref{e:th13} that $G_n$ is o-decomposable. The adjacency
matrix of $G$ immediately follows from by Theorem~\ref{e:th}.
\hfill{\rule{2mm}{2mm}} \vskip.3pc

Now, we turn our attention to ie-decomposable Riordan graphs of the
derivative type.

\begin{lemma}\label{e:th2}
 A Riordan graph $G=G_{n}(f',f)$ is
ie-decomposable if and only if
\begin{eqnarray*}
(t+t^2)f'\equiv f,\ \text{i.e.\ $[t^{2m-1}]f\equiv[t^{2m}]f$ for all
$m\ge1$.}
\end{eqnarray*}
\end{lemma}
\textbf{Proof.} Since $\left<V_o\right>$ in $G$ is a null graph and
$f''\equiv 0$ for all $f\in\mathbb{Z}[[t]]$, by Theorem~\ref{e:th12}
$G$ is ie-decomposable if and only if $(t+t^2)f'\equiv f$.
\hfill{\rule{2mm}{2mm}}

\begin{theorem}\label{e:th3}
A Riordan graph $G=G_n(f',f)$ is ie-decomposable if and only if its
binary $A$-sequence is of the form $(1,1,a_2,0,a_4,0,a_6,0,\ldots)$
where $a_{2i}$ is $0$ or $1$ for $i\ge1$, i.e.\ $A'(t)\equiv1$.
\end{theorem}
\textbf{Proof.} Since $f\equiv tA(f)$ it follows from Lemma
\ref{e:th2} that $G$ is ie-decomposable if and only if
$(1+t)f'\equiv f/t\equiv A(f)$. Applying derivative to both sides of
$f\equiv tA(f)$ we have
\begin{eqnarray*}
f'\equiv A(f)+tf'A'(f)\equiv(1+t)f'+tf'A'(f).
\end{eqnarray*}
After simplification we obtain $A'(f)\equiv1$, i.e.\ $A'(t)\equiv1$.
Since $G$ is proper, it follows that $G$ is ie-decomposable if and
only if $A(t)\equiv 1+t+\sum_{i\ge1} a_{2i}t^{2i}$ where
$a_{2i}\in\{0,1\}$. \vspace{-5mm}\hfill{\rule{2mm}{2mm}}

\begin{example} {\rm Consider $G=G_n\left({1\over
1+t^2},{t\over 1+t}\right)$. Then $A(t)=1+t$. Since $A'(t)=1$, it
follows from Theorem \ref{e:th3} that $G$ is ie-decomposable. For
instance, if  $n=9$ then}
\begin{eqnarray*}
P\left(\begin{array}{cccc|ccccc}
0 & 1 & 0 & 1 & 0 & 1 & 0 & 1 & 0 \\
1 & 0 & 1 & 1 & 0 & 0 & 1 & 1 & 0 \\
0 & 1 & 0 & 1 & 0 & 0 & 0 & 1 & 0 \\
1 & 1 & 1 & 0 & 1 & 1 & 1 & 1 & 0 \\
\hline
0 & 0 & 0 & 1 & 0 & 1 & 0 & 1 & 0 \\
1 & 0 & 0 & 1 & 1 & 0 & 1 & 1 & 0 \\
0 & 1 & 0 & 1 & 0 & 1 & 0 & 1 & 0 \\
1 & 1 & 1 & 1 & 1 & 1 & 1 & 0 & 1 \\
0 & 0 & 0 & 0 & 0 & 0 & 0 & 1 & 0
\end{array}\right)P^T=\left(\begin{array}{ccccc|cccc}
 &  &  &  &  & 1 & 1 & 1 & 1 \\
 &  &  &  &  & 1 & 1 & 0 & 1 \\
 &  & O &  &  & 0 & 1 & 1 & 1 \\
 &  &  &  &  & 1 & 1 & 1 & 1 \\
 &  &  &  &  & 0 & 0 & 0 & 1 \\
\hline
1 & 1 & 0 & 1 & 0 & 0 & 1 & 0 & 1 \\
1 & 1 & 1 & 1 & 0 & 1 & 0 & 1 & 1 \\
1 & 0 & 1 & 1 & 0 & 0 & 1 & 0 & 1 \\
1 & 1 & 1 & 1 & 1 & 1 & 1 & 1 & 0
\end{array}\right)
\end{eqnarray*}
\end{example}
where $P=\left[e_{1}\;|\; e_{3}\;|\; \cdots\;|\; e_{9}\;|\;
e_{2}\;|\; e_{4} \;|\; \cdots \;|\; e_{8}\right] ^{T}$.

\begin{theorem}
Let $G=G_{n}(f',f)$ be ie-decomposable. Then, $G$ is
$(\lfloor\log_2n\rfloor+1)$-partite for $n\geq 1$.
\end{theorem}
\textbf{Proof.} We proceed by induction on $n$. Since $G_1(f',f)$ is
a single vertex and $G_{2}(f',f)$ is bi-partite, the theorem holds
for $n=1,2$. Assume $n\ge3$. Since $\left< V_{e}\right> \cong
G_{\lfloor n/2\rfloor}$ is ie-decomposable and $\left<V_{o}\right>$
is a null graph, it follows from induction hypothesis that $G_{n}$
is $(\lfloor\log_2{\lfloor n/2\rfloor}\rfloor+2)$-partite. Now it is
enough to show that $\lfloor\log_2{\lfloor
n/2\rfloor}\rfloor=\lfloor\log_2n\rfloor-1$. If $2^k\le n<2^{k+1}$
for all $k\ge1$ then $\lfloor\log_2{\lfloor
n/2\rfloor}\rfloor=k-1=\lfloor\log_2n\rfloor-1.$ Thus we obtain the
desired result. \hfill{\rule{2mm}{2mm}}

\section{Concluding remarks and open problems}\label{last-sec}

In this paper, we use the notion of a Riordan matrix to introduce
the notion of a Riordan graph, and based on it, to introduce the
notion of an {\em unlabelled} Riordan graph. The studies conducted
by us are aimed at {\em structural properties} of (various classes
of) Riordan graphs; {\em spectral properties} of Riordan graphs are
studied by us in the follow up paper~\cite{CJKM}.

Even though our paper establishes a number of fundamental structural
results, many more such results are yet to be discovered. In
particular, we would like to extend our results on graph properties
for Riordan graphs of the Bell type to other family of Riordan
graphs. Other specific problems we would like to be solved are as
follows.

\begin{problem}\label{prob3}  {\rm Characterize {\em unlabelled} Riordan graphs.} \end{problem}

\begin{problem}\label{prob4} {\rm Enumerate {\em unlabelled} Riordan graphs.} \end{problem}

\begin{problem}\label{prob2-5} {\rm Characterize Riordan graphs whose complements are Riordan in labelled and unlabelled cases. See Section~\ref{complement-sec} for relevant observations.}\end{problem}

\begin{problem}\label{prob2} {\rm What is the complexity of recognizing labelled/unlabelled Riordan graphs?} \end{problem}

\begin{problem}\label{prob1} {\rm Characterize Riordan graphs in terms of forbidden subgraphs, or otherwise.} \end{problem}

\begin{problem}\label{prob5} {\rm Find graph invariants not considered in this paper for io-decomposable Riordan graphs of the Bell type,
e.g. the {\em independence number}, {\em Wiener index}, {\em average
path length}, and so on.} \end{problem}

Let $G_n$ be an io-decomposable Riordan graph of the Bell type. Then,
one can check that diam$(G_1)=0$ and diam$(G_n)=1$ for $n=2,3$. Based on computer experiments, in the ArXiv version of our paper \cite{CJKM0}, we have published the following conjecture.

\begin{conjecture}\label{conj1} {\rm Let $G_n$ be an
io-decomposable Riordan graph of the Bell type. Then,
\begin{align*}
2={\rm diam}(PG_n)\le {\rm diam}(G_n) \le {\rm diam}(CG_n)
\end{align*}
for $n\geq 4$. Moreover, $PG_n$ is the only graph in the class of
io-decomposable graphs of the Bell type  whose diameter is $2$ for
{\em all} $n\ge4$.}
\end{conjecture}

However, it was shown in \cite{Jung} that Conjecture~\ref{conj1} is false, because its statement does not work for certain graphs $G_n$ and values of $n$. See \cite{Jung} for the updated version of Conjecture~\ref{conj1}. In either case, we conclude our paper with the following conjecture.

\begin{conjecture}\label{conj2} {\rm We have that {\em diam}$(CG_{2^k})=k$ and there are no io-decomposable Riordan graphs $G_{2^k}\not\cong CG_{2^k}$ of the Bell type satisfying {\em diam}$(G_{2^k})=k$ for all $k\ge
1$.} \end{conjecture}

\section*{Acknowledgements}
The authors are thankful to the anonymous referee for reading carefully our manuscript and providing many useful suggestions. Also, the authors are grateful to Jeff Remmel for the useful comments on
our paper during his stay at the Applied Algebra and Optimization
Research Center, Sungkyunkwan University, in South Korea, on
September 10--19, 2017. This paper is dedicated to the memory of Jeff Remmel.


\begin{thebibliography}{99}

\bibitem{AAB} J.-P. Allouche, A. Arnold, J. Berstel, S. Brlek, W. Jockusch, S. Plouffe, B.E. Sagan, \textit{A relative of the Thue-Morse sequence}, Discrete Math. 139 (1995), 455--461.

\bibitem{BM} J. A. Bondy, U. S. R. Murty, \textit{Graph
Theory}, Graduate Texts in Mathematics, Springer 2008, ISBN
978-1-84628-969-9.

 \bibitem{BD} R. C. Brigham, R. D. Dutton, \textit{Node connectivity equals minimum degree for Pascal graphs}, Ars Combinatoria 35-A (1993),
 143--154.

\bibitem{BH} R. E. Burkard, P. L. Hammer, \textit{A Note on Hamiltonian Split Graphs}, Journal of Combinatorial Theory Series B
28 (1980), 245--248.

\bibitem{CH} G.-S. Cheon, H. Kim, \textit{The elements of finite order in the Riordan group over the complex
field}, Linear Algebra Appl. 439 (2013), 4032--4046.

\bibitem{CJ} G.-S. Cheon, S.-T. Jin, \textit{Structural properties of Riordan matrices and extending the matrices}, Linear Algebra Appl. 435 (2011), 2019--2032.

\bibitem{CJKM0} G.-S. Cheon, J.-H. Jung, S. Kitaev, S. A. Mojallal, \textit{Riordan graphs I: Structural properties}, arXiv:1710.04604.

\bibitem{CJKM} G.-S. Cheon, J.-H. Jung, S. Kitaev, S. A. Mojallal, \textit{Riordan graphs II: Spectral properties}, Linear Algebra Appl., to appear.

\bibitem{CKM} G.-S. Cheon, J. S. Kim, S. A. Mojallal, \textit{Spectral properties of Pascal graphs}, Linear and Multilinear Algebra 66 (7) (2018), 1403--1417.

\bibitem{CLMPS} G.-S. Cheon, A. Luz\'{o}n, M. A. Mor\'{o}n, L. F. Prieto-Martinez, M. Song,  \textit{Finite and infinite dimensional Lie group structures on Riordan groups},
Advances in Mathematics 319 (2017), 522--566.

\bibitem{DQ} N. Deo, M. J. Quinn, \textit{Pascal Graphs and Their Properties}, The
Fibonacci Quarterly 21 (1983), 203--214.

\bibitem{DS} E. Deutsch,  B. E. Sagan, \textit{Congruences for Catalan and Motzkin numbers and related sequences}, J. Number Theory 117 (2006), 191--215.

\bibitem{DBD} R. D. Dutton, R. C. Brigham, N. Den,  \textit{On Pascal graphs}, Presented at 23rd Southeastern Conference on Combinatorics, Graph Theory and Computing, Boca Raton, FL; 1992.

\bibitem{ES} A. Edelman, G. Strang,  \textit{Pascal Matrices}, The American Mathematical Monthly 111 (3) (2004), 189--1972.


\bibitem{DTTVZZ} R. van Dal, G. Tijssen, Z. Tuza, J. van der Veen, C. Zamfirescu, T.
Zamfirescu, \textit{ Hamiltonian properties of Toeplitz graphs},
Discrete Math. 159 (1996), 69--81.

 \bibitem{Eul}  R. Euler, \textit{Characterizing bipartite
Toeplitz graphs}, Theor. Comput. Sci 263 (2001), 47--58.

 \bibitem{Eul3} R. Euler, \textit{Coloring planar Toeplitz graphs and the stable set polytope}, Discret. Math. 276 (2004), 183--200.

 \bibitem{EZ} R. Euler, T. Zamfirescu, \textit{On Planar Toeplitz Graphs}, Graphs and Combinatorics 29 (2013), 1311--1327.


\bibitem{Ghorban} S. H. Ghorban, \textit{Toeplitz graph decomposition}, Transactions on Combinatorics 1(4) (2012), 35--41.

\bibitem{He}T.-X. He, \textit{Matrix characterizations of Riordan arrays}, Linear Algebra
Appl. 465 (2015), 15--42.

\bibitem{He1} T.-X. He, L. W. Shapiro, \textit{ Row sums and alternating sums of Riordan arrays}, Linear Algebra Appl. 507 (2016), 77--95.

\bibitem{HJ} R. A. Horn and C. R. Johnson, {\it Matrix Analysis}, Cambridge Univ. Press, 2012.

\bibitem{IK} W. Imrich, S. Klavzar, \textit{Product Graphs: Structure and Recognition}, Wiley 2000, ISBN 0-471-37039-8


\bibitem{IB} A. Ili{\'c}, Milan Ba{\v s}i{\'c}, \textit{ On the chromatic number of integral circulant graphs},
Computers \& Math. with Applications 60 (2010), 144--150.

\bibitem{Alg} C. Jean-Louis, A. Nkwanta, \textit{Some algebraic structure
 of the Riordan group}, Linear Algebra Appl. 438 (2013), 2018--2035.

\bibitem{Jung} J.-H. Jung, \textit{Diameter of io-decomposable Riordan graphs of the Bell type}, arXiv:1901.11156.

\bibitem{MRSV} D. Merlini, D.G. Rogers, R. Sprugnoli, M.C. Verri, \textit{On some
alternative characterizations of Riordan arrays}, Canad. J. Math. 49
(1997), 301--320.

\bibitem{MerSpr} D. Merlini, R. Sprugnoli, \textit{Algebraic aspects of some Riordan arrays related to binary words avoiding a pattern}, Theoretical Computer Science 412(27) (2011), 2988--3001.

\bibitem{NP} S. Nicoloso, Ugo Piestropaoli, \textit{ On the chromatic number of Toeplitz graphs},
Discrete Appl. Math. 164 (2014), 286--296.

\bibitem{S2} L.W. Shapiro, \textit{Bijections and the Riordan
group}, Theoretical Computer Science 307 (2003), 403--413.

\bibitem{SGWW} L.W. Shapiro, S. Getu, W.-J. Woan, L. Woodson, \textit{The Riordan group}, Discrete Appl. Math. 34 (1991), 229--239.

\bibitem{SGB} B. P. Sinha, S. Ghose, B. B. Bhattacharya, \textit{A further note on Pascal graphs}, The Fibonacci Quarterly 24 (1986), 251--257.

\bibitem{oeis} N. J. A. Sloane, The On-line Encyclopedia of Integer Sequences, published electronically at  https://oeis.org/

\bibitem{Sprugnoli} R. Sprugnoli, \textit{Riordan arrays and combinatorial sums}, Discrete Math. 132 (1994), 267--290.

\bibitem{DTT} R. van Dal, G. Tijssen, Z. Tuza, J. van der Veen, C. H. Zamfirescu, T. Zamfirescu, \textit{Hamiltonian properties of Toeplitz graphs}, Discret. Math. 159 (1996), 69--81.

\end{thebibliography}
\end{document}